\documentclass[a4paper,nonatbib]{article}
\pdfoutput=1
\usepackage{simplewick}
\usepackage[nonatbib]{customstyle}
\usepackage{enumitem,graphicx,xurl,xcolor}
\usepackage{booktabs,nicefrac,microtype,wrapfig,multirow,xspace}
\usepackage[titletoc]{appendix}
\usepackage{grffile}
\def\bi{\begin{itemize}[leftmargin=2em]}
\def\ei{\end{itemize}}
\def\bnum{\begin{enumerate}[leftmargin=2em,label={(\arabic*)}]}
\def\enum{\end{enumerate}}
\def\bmat{\begin{pmatrix}}
\def\emat{\end{pmatrix}}
\def\EE{\mathbb{E}}
\def\CC{\mathbb{C}}
\def\RR{\mathbb{R}}

\def\GG{\mathcal{G}}
\def\PP{\mathbb{P}}
\def\WW{\hat{W}}
\def\HH{\hat{H}}
\def\MM{\mathcal{M}}
\def\rmd{\mathrm{d}}
\DeclareMathOperator{\rme}{e}
\def\1{\mathbb{I}}
\def\xx{\mathbf{x}}
\def\uu{\mathbf{u}}
\def\vv{\mathbf{v}}
\def\CN{\mathcal{N}_{\CC}}
\def\RN{\mathcal{N}_{\RR}}
\DeclareMathOperator{\He}{He}
\def\Var{\mathsf{Var}}
\newcommand{\al}[1]{%
  \begin{align}
    #1
  \end{align}
}
\def\mkakko#1{\left(#1\right)}
\def\ckakko#1{\left\{#1\right\}}
\def\kkakko#1{\left[#1\right]}
\renewcommand{\bar}[1]{\overline{#1}}
\renewcommand{\tilde}[1]{\widetilde{#1}}
\renewcommand{\epsilon}{\varepsilon}
\let\Re\relax
\let\Im\relax
\DeclareMathOperator{\Re}{Re}
\DeclareMathOperator{\Im}{Im}
\DeclareMathOperator{\Tr}{Tr}
\newcommand{\MP}{Marchenko--Pastur\xspace}
\newcommand{\moment}{\mathfrak{m}}
\newtheorem{theorem}{Theorem}
\theoremstyle{lemma}
\newtheorem{lemma}{Lemma}
\theoremstyle{proposition}
\newtheorem{proposition}{Proposition}
\allowdisplaybreaks

\title{Widely Linear Gaussian Equivalence for Complex Random Features}

\author{%
  Takuya Kanazawa\\
  Department of Business Administration, Kobe Gakuin University\\
  1-1-3 Minatojima, Chuo-ku, Kobe 650-8586, Japan
}

\begin{document}

\maketitle
\thispagestyle{plain}

\begin{abstract}
Entrywise nonlinearity makes the spectral analysis of high-dimensional random feature matrices substantially more difficult than that of the underlying matrix products. We establish a spectral Gaussian equivalence for complex random feature matrices with general, possibly non-holomorphic activation functions. The widely linear surrogate combines the preactivation matrix and its complex conjugate with independent Gaussian residual noise that may be improper. Under a condition ensuring proper preactivations, we prove convergence of the expected empirical spectral distributions of the nonlinear feature covariance matrix and its surrogate to a common limit in the proportional high-dimensional regime. The proof combines a complex Hermite expansion with a graphical moment method. We also characterize broad classes of activation functions for which one or both response channels vanish. Numerical experiments confirm the predicted spectral equivalence.
\end{abstract}

\clearpage
\hrule\vspace{-0.6em}
\tableofcontents
\vspace{1.6em}
\hrule
\vspace{1.6em}

\section{Introduction}

Feature matrices obtained by applying entrywise nonlinearities to large random matrix products arise naturally in random feature methods, neural networks, and signal processing. Random features were introduced as explicit, low-dimensional randomized representations for approximating kernel methods at reduced computational cost \autocite{Rahimi2007}. A representative high-dimensional model is $f(H)$, where $H=XW/\sqrt{d}$ is the preactivation matrix, $X$ is a data matrix, $W$ is a random weight matrix, and the activation function $f$ is applied entrywise. Although the entries of $X$ and $W$ are independent in standard models, the entries of $H$ and $f(H)$ are strongly dependent through the shared factors. Describing the spectral behavior of the resulting feature covariance matrix is a nontrivial problem in random matrix theory (RMT) and has motivated a substantial line of work, including the foundational analysis by \textcite{Pennington2017}.

The Gaussian equivalence principle (GEP) addresses this difficulty by replacing a nonlinear random feature matrix with a simpler surrogate assembled from Gaussian random matrices and low-order response coefficients. Depending on the formulation, a GEP may describe limiting spectra, normalized resolvents, local fields, or training and test errors in high-dimensional supervised learning \autocite{Louart2018,Goldt2021,Hastie2022,Mei-Montanari-2022}. Its central premise is that the effect of the nonlinearity on macroscopic observables can often be captured by a small number of response coefficients and residual moments, while retaining the underlying random matrix structure.

The complex setting introduces structures absent from the conventional real-valued GEP. Complex-valued neural networks (CVNNs), and complex random feature models more generally, often employ non-holomorphic activation functions \autocite{Hirose2012book}. The linear projection of such an activation generally responds to both a complex preactivation and its conjugate. Moreover, for improper complex variables, the covariance alone does not determine the second-order statistics, and the pseudo-covariance must also be retained. A general complex surrogate should therefore adopt the widely linear structure familiar from complex-valued statistical signal processing \autocite{Schreier2010book}, while allowing its Gaussian residual to be improper.

Complex extensions were briefly mentioned in earlier spectral analyses. \textcite[Remark~1.3]{Peche2019} considered matrices whose real and imaginary parts are independent and identically distributed, while \textcite[Remark~E.1]{Piccolo2021} described an analogous extension for real-valued activations applied to complex feature and weight matrices. These earlier extensions are naturally consistent with the proper, circular setting and a single-channel surrogate. However, proper input statistics alone do not ensure such a reduction: a general non-holomorphic activation may produce a conjugate-linear response and an improper residual even from a proper preactivation. For improper inputs, these structures arise in an even more general form.

Motivated by this distinction, we formulate a general complex GEP with the widely linear surrogate $\mu 1_n 1_p^\top + c_1 H + c_2 \bar{H} + Z$. Here, $\mu$ is the mean response, $c_1$ and $c_2$ describe the linear and conjugate-linear response channels, and $Z$ is independent Gaussian noise whose variance and pseudo-variance match those of the nonlinear residual. The formulation reduces to the conventional real GEP when the model degenerates to the real axis, while allowing richer behavior in genuinely complex settings. In particular, a proper preactivation may yield an improper output under a non-holomorphic activation.

The general complex GEP is proposed as a conjectural principle for a broad class of macroscopic observables. Our main rigorous result, Theorem~\ref{thm:complex_gep}, proves its spectral version under additional assumptions: the expected empirical spectral distributions of the nonlinear feature covariance matrix and its widely linear surrogate converge weakly to the same limit. The limiting spectrum of the surrogate is derived separately in the companion paper \autocite{kanazawa_rmt_2026}. Thus, the present paper establishes the equivalence between the nonlinear model and its surrogate, whereas the companion paper solves the spectral problem for the resulting surrogate.

The proof combines a complex Hermite expansion with a graphical moment method. Hermite expansions have played a central role in Gaussian-equivalence analyses of real nonlinear random matrices \autocite{Cheng_Singer_2013,Mei-Montanari-2022}; here, their complex extension separates the linear, conjugate-linear, and higher-order components of the activation. After a finite Hermite truncation, the dependence structure is encoded by contraction graphs. Power counting identifies the configurations that survive at leading order, and their contributions match the Wick contractions of the widely linear surrogate. The remaining approximation steps remove the truncation and auxiliary normalizations.

We also characterize broad classes of complex activation functions for which one or both response coefficients vanish. For instance, the conjugate-linear coefficient vanishes for Split-ReLU, a popular activation in the CVNN literature \autocite{Trabelsi2018,Gao2019}, even though Split-ReLU is non-holomorphic. Such reductions clarify when the complex surrogate genuinely requires two response channels and when it collapses to a simpler form. Numerical experiments support the spectral complex GEP for diverse activation functions and covariance structures, including several settings outside the assumptions of Theorem~\ref{thm:complex_gep}.

The remainder of the paper is organized as follows. Section~\ref{sc:relwork} reviews Gaussian equivalence, nonlinear random matrix models, CVNNs, and widely linear methods. Section~\ref{sc:notation} introduces the complex random feature model and the notation used throughout. Section~\ref{sc:surrogate} derives the scalar widely linear decomposition, formulates the general complex GEP, examines its reduction to the real case, and studies conditions under which one or both response channels vanish. Section~\ref{sc:gep_rigorous} proves Theorem~\ref{thm:complex_gep}. Section~\ref{sc:numerics} presents numerical experiments, and Section~\ref{sc:conclusions} concludes the paper. The appendices collect supplementary proofs and technical details used in the main argument.

\section{Related Work}
\label{sc:relwork}
\subsection{Gaussian Equivalence}

A classical universality phenomenon in RMT is that spectral observables can be insensitive to the detailed distributions of the matrix entries, depending instead on coarse properties such as their low-order moments and the symmetry and dimensional structure of the matrix ensemble \autocite{Tao2012book}. Related universality phenomena have emerged for nonlinear random matrices arising in machine learning and high-dimensional statistics, notably through the GEP.

Precursors to the spectral GEP appeared in the study of inner-product kernel random matrices. For smooth kernels, \textcite{El_Karoui_2010} showed that the kernel matrix is asymptotically equivalent to a linear transformation of the sample covariance matrix, up to a scalar shift and a finite-rank term. \textcite{Cheng_Singer_2013} extended the spectral analysis to nonsmooth and discontinuous kernels and showed that the limiting spectral distribution depends only on the linear coefficient and total variance of a Hermite-type expansion of the kernel.

A closely related line of work concerns real-valued random feature matrices of the form $Y=f(WX)$. Random features were originally introduced as explicit randomized approximations to kernel methods \autocite{Rahimi2007}. Their high-dimensional spectral behavior was later investigated through nonlinear random matrix models. For independent real Gaussian matrices $W$ and $X$, \textcite{Pennington2017} characterized the asymptotic spectrum of $YY^\top$ by a quartic equation whose coefficients depend on $f$ through a small number of Gaussian averages. \textcite{Louart2018} considered a similar model with a deterministic data matrix $X$ and derived a deterministic equivalent for the resolvent of the nonlinear random feature Gram matrix, governed by its expected covariance structure. \textcite{Peche2019} identified the limiting distribution obtained by \textcite{Pennington2017} with that of an information-plus-noise model consisting of a linear component and independent Gaussian noise. The bulk result was subsequently extended to sub-Gaussian matrices, and the corresponding extremal eigenvalues were also studied \autocite{Benigni2021,Benigni_largest_2022}. For real-valued random features with an additive bias term, \textcite{Piccolo2021} developed a resolvent-based analysis.

Gaussian equivalence has also become an important tool for analyzing high-dimensional learning problems, where the observables of interest are training and test errors rather than the feature spectrum alone. \textcite{Hastie2022} analyzed the asymptotic prediction risk of ridgeless regression for both linear and nonlinear random feature models, while \textcite{Mei-Montanari-2022} derived the full asymptotic prediction risk and related nonlinear random feature regression to a Gaussian covariates model with matching low-order statistics. \textcite{HuLu2023} subsequently proved a general universality theorem showing that, for training and generalization errors under convex losses, nonlinear random features can be replaced by a linear Gaussian model with matching covariance. Further equivalence results cover classification and max-margin problems \autocite{Montanari2019maxmargin,Montanari2023universalitymaxmarginclassifiers} and nonconvex empirical risk minimization \autocite{Montanari_Saeed_2022}. A parallel line of work has developed a local-field formulation of the GEP for structured data and applied it to high-dimensional learning problems \autocite{Goldt2020,Goldt2021,Gerace2021}.

The precise meaning of GEP is therefore not uniform across the literature: depending on the formulation, the equivalence may concern spectral observables, training and test errors, or finite collections of local fields.

Subsequent extensions of GEP have addressed polynomial scaling regimes \autocite{HLM2024}, self-attention layers \autocite{Hayase2025attention}, and deep random feature models \autocite{Schroder2023,Bosch2023}, alongside traffic-based analyses of profiled Pennington--Worah matrices \autocite{Dabo2024traffic} and entrywise nonlinear functions of orthogonally invariant matrices \autocite{Speicher2024GOE}. A simple form of the GEP may fail in certain polynomial scaling regimes \autocite{Wen2025}, while heavy-tailed nonlinear models exhibit further distinct behavior \autocite{Guionnet2026kernel}.

Building on these developments in the real-valued setting, the present work extends the spectral GEP to complex random feature matrices with general non-holomorphic activations, leading to a surrogate with both linear and conjugate-linear response channels.

\subsection{CVNNs and Widely Linear Models}

Complex-valued signals arise naturally in a wide range of scientific and engineering applications, including biomedical imaging, acoustics, optics and wireless communications. This has motivated a substantial body of work on statistical signal processing for complex-valued data \autocite{Mandic2009book,Schreier2010book,Ollila2011, Adali2011,Bluhm_Mota_Cavalcante_2025}. CVNNs extend real-valued networks by allowing inputs, activations, outputs, and trainable parameters to take values in the complex domain, thereby retaining amplitude and phase information throughout the network \autocite{Hirose2012book,Bassey2021review,Lee2022review,Barrachina2023,Hammad2024}. CVNNs commonly use non-holomorphic activations, whose linear projections generally respond to both a complex preactivation $h$ and its conjugate $\bar h$, naturally giving rise to a form such as $c_1h+c_2\bar h$. Such dependence on a variable and its conjugate is known in the statistical signal processing literature as a \emph{widely linear} structure \autocite{Brown1969,Picinbono_Chevalier_1995}. 

Classical complex statistical models often assume \emph{proper} signals \autocite{Neeser1993,Picinbono1994}, whose second-order structure is characterized by the mean and covariance alone. For \emph{improper} signals, however, the pseudo-covariance carries additional second-order information. Since improper statistics arise in many real-world datasets, a general modeling framework should retain this information. 

The present work connects these developments in the setting of high-dimensional complex random feature models. Non-holomorphic activations may produce conjugate-linear response channels, while improper complex statistics require the pseudo-covariance to be retained alongside the covariance. Together, these features lead naturally to a widely linear Gaussian surrogate rather than the conventional single-channel model.

\section{Notation and Definitions}
\label{sc:notation}

For a positive integer $K$, let $[K]$ denote the set $\{1,2,\dots,K\}$. For the dimension parameters $n$, $p$, and $d$ introduced below, we use $i,j,\dots$ to index $[n]$, $\alpha,\beta,\dots$ to index $[p]$, and $\mu,\nu,\dots$ to index $[d]$.

The sets of real and complex numbers are denoted by $\RR$ and $\CC$, respectively. $\CC_+$ denotes the set of complex numbers with positive imaginary part. For $z\in\CC$, its real and imaginary parts are denoted by $\Re z$ and $\Im z$, respectively. The complex conjugate of $z$ is denoted by $\bar{z}$. The $d$-dimensional vector with all elements equal to $1$ is denoted by $1_d$. The $\ell_2$-norm of a vector $\xx$ is denoted by $\|\xx\|_2$ and the $\ell_\infty$-norm by $\|\xx\|_\infty$. 

For a complex matrix $A\in\CC^{n \times p}$, $\bar{A}$ denotes the complex conjugate of $A$, $A^\top$ the transpose of $A$, and $A^*$ the conjugate transpose of $A$, i.e., $A^*=\bar{A^{\top}}$. The operator norm of $A$ is denoted by $\|A\|_{\rm op}$, the nuclear norm of $A$ by $\|A\|_{\star}$, and the Frobenius norm of $A$ by $\|A\|_{\rm F}$. The $d$-dimensional identity matrix is denoted by $\1_d$. For $z\in\CC_+$, we denote by
\al{
	m^{}_A(z) \equiv \frac{1}{p} 
	\Tr\bigg[ \left( \frac{1}{n}A^{*}A - z\1_p \right)^{-1} \bigg] 
}
the normalized Stieltjes transform of the empirical spectral distribution of $A^{*}A/n$.

The expectation and variance of a random variable $z$ are denoted by $\EE[z]$ and $\Var(z)\equiv \EE[\left|z-\EE[z]\right|^2]$, respectively. A real-valued Gaussian distribution with mean $\mu\in\RR$ and variance $v\geq 0$ is denoted by $\RN(\mu,v)$. A complex-valued Gaussian distribution with mean $\mu\in\CC$, variance $v \geq 0$, and \emph{pseudo-variance} $\tau\in\CC$ with $|\tau|\leq v$, is denoted by $\CN(\mu, v, \tau)$. It is called nondegenerate when $|\tau|<v$. In the nondegenerate case, its probability density function $p(z)$ is given by \autocite{Schreier2010book}
\al{
	p(z) & = \frac{1}{\pi \sqrt{v^2 - |\tau|^2}} 
	\exp\left(
		- \frac{(v - \Re \tau) x^2 
		- 2 (\Im \tau) xy 
		+ (v + \Re \tau) y^2}
		{v^2 - |\tau|^2}
	\right),
}
where we denoted $z-\mu\equiv x+iy$ with $x,y\in\RR$. One can verify $v = \EE[|z-\mu|^2]$ and $\tau = \EE[(z-\mu)^2]$. The variable $z$ is called \emph{proper} when $\tau=0$, and \emph{improper} otherwise. For $v>0$ and $\tau=0$, its distribution is circularly symmetric%
\footnote{Propriety does not mean circular symmetry for general non-Gaussian distributions.} with density $\exp(-|z-\mu|^2/v)/(\pi v)$. The pseudo-variance $\tau$ determines the anisotropy and orientation of the complex Gaussian distribution. When $|\tau|=v>0$, the distribution is degenerate and is supported on a line in the complex plane; when $v=\tau=0$, it is concentrated at its mean.

We consider complex random features generated by nonlinear transformations of random linear projections. Let $X\in\CC^{n\times d}$ denote the input data matrix and $W\in\CC^{d\times p}$ the random weight matrix. The entries of $X$ and $W$ are assumed to be independent Gaussian random variables satisfying 
\al{
	X_{i\mu}\sim\CN(0,v_x,\tau_x) 
	\qquad\text{and}\qquad 
	W_{\mu\alpha}\sim\CN(0,v_w,\tau_w)\,, 
}
where $v_x, v_w>0$ with $|\tau_x|\leq v_x$ and $|\tau_w|\leq v_w$. Define the column-normalized weight matrix $\WW$ as
\al{
	\WW \equiv WS \in\CC^{d\times p}\,, \qquad 
	S \equiv \operatorname{diag}
	(s_1,\dots,s_p)
	\in\RR^{p\times p}\,, \qquad 
	s_\alpha \equiv \sqrt{\frac{dv_w}
	{\sum_{\mu=1}^{d}|W_{\mu\alpha}|^2}}\,.
	\label{eq:column-normalization}
}
The norm of each column of $\WW$ is normalized to $\sqrt{dv_w}$. Note that $S$ is close to $\1_p$ when $d\gg 1$.

We define the feature matrices by
\al{
	H \equiv \frac{XW}{\sqrt{d}} 
	\in\CC^{n\times p}\,, \qquad 
	f(H) \in\CC^{n\times p}\,, 
	\qquad\text{and}\qquad 
	\HH \equiv \frac{X\WW}{\sqrt{d}} 
	= HS \in\CC^{n\times p}\,.
	\label{eq:feature-matrices}
}
Here and throughout, the activation function $f:\CC\to\CC$ is applied entrywise to matrices.%
\footnote{In complex analysis, the notation $f(z)$ often suggests holomorphic dependence, while a general non-holomorphic function may be written as $f(z,\bar{z})$. For brevity, we use the notation $f(z)$ throughout this paper regardless of whether $f$ is holomorphic, anti-holomorphic, or non-holomorphic.} Our primary object of interest is the uncentered empirical covariance matrix $\frac{1}{n}f(H)^*f(H)$. The matrix $\HH$ will be used as an intermediate building block.

We assume throughout that there exist constants $C_f>0$ and $q\geq0$ such that
\al{
	|f(z)-f(w)| \leq C_f \bigl( 1+|z|^q+|w|^q \bigr) |z-w|\,, 
	\qquad z,w\in\CC\,.
	\label{eq:activation_increment_bound}
}
This condition implies that $f$ is continuous and locally Lipschitz, and hence differentiable almost everywhere in the real sense. Its Wirtinger derivatives $\partial_z f$ and $\partial_{\bar z}f$ are therefore defined almost everywhere and grow at most polynomially. The condition also implies that $f$ itself has at most polynomial growth, so that all its Gaussian moments, as well as those of its Wirtinger derivatives, are finite. This regularity condition includes globally Lipschitz activations and arbitrary polynomials in $z$ and $\bar z$.

Our analysis targets the proportional asymptotic regime
\al{
	n,d,p \to \infty \qquad\text{with}\qquad 
	\frac{n}{d} \to \phi \in (0,\infty), \quad 
	\frac{p}{d} \to \psi \in (0,\infty). 
	\label{eq:lmlm}
}
For each fixed pair $(i,\alpha)$, the random variables $X_{i\mu}W_{\mu\alpha}$ are independent and identically distributed over $\mu$, with zero mean, variance $v_xv_w$, and pseudo-variance $\tau_x\tau_w$. Consequently, for every finite $d$,
\al{
	\EE[|H_{i\alpha}|^2] = v_xv_w \equiv v_h > 0
	\qquad \text{and} \qquad 
	\EE[H_{i\alpha}^2] = \tau_x\tau_w \equiv \tau_h\in\CC\,.
}
As $d\to\infty$, the bivariate real central limit theorem gives, for each fixed $(i,\alpha)$,
\al{
	H_{i\alpha} \xrightarrow{~\mathrm{d}~} 
	\CN(0,v_h,\tau_h)\,. 
}
This is an entrywise convergence statement and does not imply that the full matrix $H$ converges to a matrix with independent Gaussian entries. Its entries remain dependent because they share factors from $X$ and $W$.

\section{Complex Gaussian Surrogate}
\label{sc:surrogate}

In this section, we derive the complex Gaussian surrogate and formulate a general complex GEP conjecture. We subsequently examine its reduction to the real GEP and identify conditions under which one or both response coefficients vanish. A rigorous spectral version of the complex GEP is formulated and proved in Section~\ref{sc:gep_rigorous}.

\subsection{Widely Linear Decomposition}
\label{sc:widely_linear}

Consider a nondegenerate zero-mean complex Gaussian variable $h\sim\CN(0,v_h,\tau_h)$ with $v_h>|\tau_h|$. Then $f(h)$ admits the unique decomposition
\al{
	f(h) = \mu + c_1h + c_2\bar{h} + z\,,
	\label{eq:ansatz}
}
where $\mu+c_1h+c_2\bar{h}$ is the orthogonal projection of $f(h)$ onto $\mathsf{span}\{1,h,\bar{h}\}$ and $z$ is the corresponding residual. In particular, $\EE[z]=\EE[\bar{h}z]=\EE[hz]=0$. We now derive explicit formulas for these parameters and determine the second moments of $z$.%
\footnote{The construction below parallels the Bussgang decomposition and the theory of widely linear estimation in signal processing \autocite{Bussgang1952,Minkoff1985,Ermolova2004,Schreier2010book,Demir2021}.}
We stress that this is a single-variable problem and does not involve matrices.

Matching the mean on both sides of \eqref{eq:ansatz} gives%
\footnote{The mean parameter $\mu$ should not be confused with the data index $\mu\in[d]$ used elsewhere.}
\al{
	\mu = \EE[f(h)]\,.
}
To determine $c_1$ and $c_2$, define the cross-moments
\al{
	\kappa_1\equiv\EE[(f(h)-\mu)\bar{h}]
	\qquad\text{and}\qquad
	\kappa_2\equiv\EE[(f(h)-\mu)h]\,.
	\label{eq:kkr4sx}
}
The conditions $\EE[\bar{h}z]=\EE[hz]=0$ then yield
\al{
	\bmat
		\kappa_1 \\ \kappa_2
	\emat
	= \bmat
		v_h & \bar{\tau_h}
		\\
		\tau_h & v_h
	\emat
	\bmat
		c_1 \\ c_2
	\emat .
	\label{eq:kaka}
}

The coefficients can also be expressed in terms of the expected Wirtinger derivatives of $f$. To see this, write $h=x+iy$, where $x,y\in\RR$. The real vector $\xx=(x,y)^\top$ is a zero-mean bivariate Gaussian vector with covariance matrix 
\al{
	\Sigma = \frac{1}{2}
	\bmat
		v_h+\Re(\tau_h) & \Im(\tau_h)
		\\
		\Im(\tau_h) & v_h-\Re(\tau_h)
	\emat .
}
The real multivariate Stein identity \autocite{Stein1981,Couillet_Liao_2022}, applied separately to the real and imaginary parts of $f$, gives
\al{
	\bmat
		\EE[f(h)x] \\ \EE[f(h)y]
	\emat
	= \Sigma
	\bmat
		\EE[\partial_x f(h)]
		\\
		\EE[\partial_y f(h)]
	\emat .
	\label{eq:real_stein_h}
}
Using $\partial_h=(\partial_x-i\partial_y)/2$ and $\partial_{\bar{h}}=(\partial_x+i\partial_y)/2$, or equivalently $\partial_x=\partial_h+\partial_{\bar{h}}$ and $\partial_y=i(\partial_h-\partial_{\bar{h}})$, \eqref{eq:real_stein_h} becomes
\al{
	\bmat \kappa_1 \\ \kappa_2 \emat
	= \bmat
		v_h & \bar{\tau_h}
		\\
		\tau_h & v_h
	\emat
	\bmat
		\EE[\partial_h f(h)] 
		\\ 
		\EE[\partial_{\bar{h}}f(h)]
	\emat .
	\label{eq:complex_stein_h}
}
Comparison of \eqref{eq:complex_stein_h} with \eqref{eq:kaka}, together with $v_h>|\tau_h|$, gives
\al{
	c_1 = \frac{v_h\kappa_1-\bar{\tau_h}\kappa_2}
	{v_h^2-|\tau_h|^2} =
	\EE[\partial_h f(h)]\,,
	\label{eq:stein_gain_1}
	\\
	c_2 = \frac{v_h\kappa_2-\tau_h\kappa_1}
	{v_h^2-|\tau_h|^2}
	= \EE[\partial_{\bar{h}}f(h)]\,.
	\label{eq:stein_gain_2}
}

Let $v_z\in\RR_{\geq0}$ and $\tau_z\in\CC$ denote the variance and pseudo-variance of $z$, respectively. Using \eqref{eq:kkr4sx} and \eqref{eq:kaka}, we obtain
\al{
	v_z & = \EE[|f(h)-\mu-c_1h-c_2\bar{h}|^2]
	\\
	& =
	\EE[|f(h)-\mu|^2] 
	- v_h(|c_1|^2+|c_2|^2) 
	- 2 \Re(\tau_h c_1 \bar{c_2}) \,,
	\label{eq:vc87vc9}
	\\
	\tau_z & =
	\EE[(f(h)-\mu-c_1h-c_2\bar{h})^2]
	\\
	& = \EE[(f(h)-\mu)^2]
	- \tau_h c_1^2 - \bar{\tau_h}c_2^2
	- 2 v_h c_1c_2 \,.
	\label{eq:RTEGD1F}
}

An important point is that $z$ is generally non-Gaussian and depends on $h$. Although $z$ has no correlation with either $h$ or $\bar{h}$, this does not imply independence. As a simple example, let $h\sim\CN(0,1,0)$ and consider $f(h) = h + a|h|^2 \bar{h}$ with $a\in\CC\setminus\{0\}$. The above formulas give $\mu=0$, $c_1=1$ and $c_2=2a$. For the residual $z\equiv f(h)-h - 2a\bar{h}$, we have
\al{
	\EE[f(h)^2] = 4a\ne 0 \qquad \text{and}\qquad 
	\EE[z h^2\bar{h}] = 2a
	\ne \EE[z]\EE[h^2\bar{h}]=0\,.
}
The first relation shows that a non-holomorphic activation can produce an improper output even when the input is proper and circularly symmetric. The second shows that $z$ remains dependent on $h$ after its linear correlations have been removed. Equations~\eqref{eq:vc87vc9} and \eqref{eq:RTEGD1F} give $v_z=2|a|^2$ and $\tau_z=0$, so the residual itself is proper. If $z$ were Gaussian, its propriety would imply $\EE[|z|^4]=2v_z^2$. Instead, direct calculation gives
\al{
	\EE[|z|^4] = 176|a|^4 \ne 2v_z^2=8|a|^4.
}
Thus, $z$ is non-Gaussian.

\subsection{A General Complex GEP Conjecture}
\label{sc:gengep}

While the decomposition \eqref{eq:ansatz} operates at the single-variable level, we speculate that an analogous principle holds at the matrix level, providing a natural extension of the real GEP to complex data $X$, complex weights $W$, and complex activations $f$.

\paragraph{Conjecture.} \emph{Under the setup of Section~\ref{sc:notation} and the nondegeneracy condition $v_h>|\tau_h|$, the nonlinear feature matrix $f(H)$ is asymptotically equivalent, for a broad class of high-dimensional macroscopic observables, to the widely linear Gaussian model}
\al{
    f(H) \approx \mu 1_n1_p^{\top}
    + c_1H + c_2 \bar{H} + Z 
    \label{eq:sdf9}
}
\emph{in the proportional limit \eqref{eq:lmlm}, where $Z\in\CC^{n\times p}$ is independent of $H$ and has \mbox{i.i.d.} entries $Z_{i\alpha}\sim\CN(0,v_z,\tau_z)$. The parameters $\{\mu,c_1,c_2,v_z,\tau_z\}$ are defined in Section~\ref{sc:widely_linear}.}

The equivalence in \eqref{eq:sdf9} is intended as a macroscopic high-dimensional principle rather than an entrywise approximation. It encompasses spectral observables such as normalized resolvents and the bulk spectral density, and is expected to extend to broader quantities arising in supervised learning problems such as training and test errors. For example, a more specific spectral formulation would compare normalized traces of the form
\al{
    \frac{1}{p}\Tr\left[
    \left( \frac{1}{n}f(H)^*f(H)-z\1_p \right)^{-1}A
    \right]
}
with their counterparts under the Gaussian model, where $A\in\CC^{p\times p}$ is deterministic with bounded operator norm and $z\in\CC_+$. Likewise, for projected fields such as $\uu^\top f(H)\vv$, the equivalence is expected when the deterministic unit vectors $\uu$ and $\vv$ are sufficiently delocalized, for instance when
\al{
    \|\uu\|_\infty  =  O \big(n^{-1/2}\big)
    \qquad\text{and}\qquad  \|\vv\|_\infty
    =  O \big(p^{-1/2} \big)\,.
}
Such conditions prevent the observable from isolating individual non-Gaussian entries. The conjecture does not assert entrywise convergence, and microscopic quantities such as local eigenvalue statistics and the exact joint distribution of eigenvectors are not included in its intended scope.

The difference from the real GEP is twofold. First, \eqref{eq:sdf9} contains both the response to $H$ and the response to its complex conjugate $\bar{H}$. Second, the Gaussian noise $Z$ may be improper, in the sense that circular symmetry can be broken and $\tau_z$ can be nonzero. To the best of our knowledge, an explicit Gaussian-equivalence formulation incorporating both response channels together with possibly improper residual noise has not appeared in the literature. The complex GEP is characterized by the five surrogate parameters $\{\mu,c_1,c_2,v_z,\tau_z\}$, comprising nine real degrees of freedom; the input variance $v_h$ and pseudo-variance $\tau_h$ are regarded as fixed ensemble parameters and are not included in this count. By contrast, the real GEP involves only three real parameters associated with the activation, namely $\EE[f(x)]$, $\EE[f(x)x]$, and $\EE[f(x)^2]$ for a real Gaussian variable $x$ and a real-valued activation $f$.

The rigorous result proved in this work is a more restrictive spectral formulation of \eqref{eq:sdf9}, stated as Theorem~\ref{thm:complex_gep} in Section~\ref{sc:gep_rigorous} and established by a spectral-moment expansion. Before turning to this result, the following subsections examine several special cases and structural consequences of the proposed equivalence.

\subsection{Reduction to the Real GEP}

Consider a real activation $f:\mathbb R\to\mathbb R$ and take the real-axis limit $\tau_x\to v_x$ and $\tau_w\to v_w$, so that $X$, $W$, and $H$ become real and $\tau_h\to v_h$. This recovers the setting of the conventional real GEP. In this limit, $h\sim\RN(0,v_h)$ and \eqref{eq:kkr4sx} gives
\al{
	\kappa_1,\kappa_2 \to
	\kappa \equiv \EE[f(h)h] \in\RR\,.
	\label{eq:9Rfge3}
}
The two response channels also merge because $H=\bar{H}$, leaving only the combination $c_1+c_2$ in the surrogate. Although the individual expressions in \eqref{eq:stein_gain_1} and \eqref{eq:stein_gain_2} become singular as $\tau_h\to v_h$, their sum has a finite limit. Along the real-axis degeneration, $\tau_h$ is real, and hence
\al{
	c_1+c_2 &=
	\frac{v_h\kappa_1-\tau_h\kappa_2
	+v_h\kappa_2-\tau_h\kappa_1}
	{v_h^2-\tau_h^2}
	\\
	& =
	\frac{\kappa_1+\kappa_2}{v_h+\tau_h}
	\to
	\frac{\kappa}{v_h} \in\RR\,.
	\label{eq:real_gain_limit}
}
Thus, the apparent $0/0$ singularity reflects only the loss of distinction between the two complex response channels; their sum reduces smoothly to the single real response coefficient.

Using \eqref{eq:9Rfge3} and \eqref{eq:real_gain_limit} in \eqref{eq:vc87vc9} and \eqref{eq:RTEGD1F}, we obtain
\al{
	v_z,\tau_z & \to \EE[(f(h)-\mu)^2]
	- v_h(c_1+c_2)^2
	\\
	& = \Var(f(h))
	- \frac{\kappa^2}{v_h}\,.
}
The equality $\tau_z=v_z$ shows that the Gaussian residual becomes real in this limit. Consequently, the complex surrogate \eqref{eq:sdf9} reduces to
\al{
	f(H) \approx \EE[f(h)]1_n1_p^{\top}
	+ \frac{\EE[f(h)h]}{v_h}H
	+ \mkakko{
		\Var(f(h)) - \frac{\EE[f(h)h]^2}{v_h}
	}^{1/2}
	Z_{\RR}\,,
	\label{eq:real_gep_reduction}
}
where $h\sim\RN(0,v_h)$ and $(Z_{\RR})_{i\alpha}\sim\RN(0,1)$. Up to notational differences, \eqref{eq:real_gep_reduction} agrees exactly with the standard real GEP formulas established in the literature \autocite{Pennington2017,Peche2019,Dhifallah2020,Goldt2020,Benigni2021,Gerace2021,Goldt2021,Mei-Montanari-2022,HuLu2023}.

We note that restricting $f$ to a real-output complex activation $f:\mathbb C\to\mathbb R$ alone does not reduce the complex GEP to the real one. If $f:\CC\to\RR$, then $\mu\in\RR$, while $c_1=c$ and $c_2=\bar{c}$ for some $c\in\CC$. Equations \eqref{eq:vc87vc9} and \eqref{eq:RTEGD1F} then give
\al{
	v_z = \tau_z = \EE[(f(h)-\mu)^2] - 2v_h |c|^2
	- 2 \Re(\tau_h c^2) \in \RR\,.
}
The noise matrix $Z$ is therefore real-valued, and the surrogate takes the form 
\al{
	f(H) \approx \EE[f(h)] 1_n 1_p^\top
	+ 2\Re(c H) + \sqrt{v_z} Z _\RR\,.
	\label{eq:realf}
}
This remains a two-channel problem, since the conjugate-mixing structure is retained in the second term. Indeed, as shown in our companion paper \autocite[Section~3.4]{kanazawa_rmt_2026}, the limiting Stieltjes transform of the covariance matrix associated with \eqref{eq:realf} generally satisfies a sextic equation, rather than the quartic equation found for the real GEP \autocite{Pennington2017,Benigni2021}.

\subsection{Conditions for Vanishing GEP Coefficients}
\label{sc:vavasd}

The general complex GEP contains two response channels governed by $c_1$ and $c_2$. These channels need not both survive: symmetries and other structural properties of the activation may force either coefficient to vanish, leading to a reduced surrogate model. We explore these reductions below for both improper and proper inputs, together with representative examples.

Throughout this subsection, $f$ is assumed to be fixed independently of the Gaussian parameters $v_h$ and $\tau_h$. Derivatives of locally Lipschitz functions and identities involving them are understood to hold almost everywhere unless otherwise stated. We freely use the standard fact that a locally Lipschitz function whose weak gradient vanishes almost everywhere is constant on a connected domain.

\subsubsection{General inputs}

If $c_1=c_2=0$, the surrogate \eqref{eq:sdf9} reduces to the Gaussian residual $Z$. Therefore, when $v_z>0$, the GEP predicts that the bulk spectrum of $f(H)^*f(H)/n$ follows a scaled \MP law \autocite{Marchenko1967}. We now characterize when both coefficients vanish.

\begin{lemma}\label{lem:bothczero}
One has $c_1=c_2=0$ for every $v_h>0$ and $\tau_h\in\CC$ satisfying $|\tau_h|<v_h$ if and only if $f(h)=f(-h)$ for every $h\in\CC$.
\end{lemma}

\begin{proof}
Let $h = x + iy$. We decompose the activation function into its even and odd components: $f(h) = f_{\text{e}}(h) + f_{\text{o}}(h)$, where $f_{\text{e}}(h) = (f(h)+f(-h))/2$ and $f_{\text{o}}(h) = (f(h)-f(-h))/2$. 
From \eqref{eq:stein_gain_1} and \eqref{eq:stein_gain_2}, the condition $c_1 = c_2 = 0$ is equivalent to $\EE[\partial_x f] = \EE[\partial_y f] = 0$.

First, observe that the partial derivatives $\partial_x f_{\text{e}}$ and $\partial_y f_{\text{e}}$ of $f_{\rm e}$ are odd functions. Since the zero-mean Gaussian density is an even function of $\xx = (x, y)^\top$, the integrals $\EE[\partial_x f_{\text{e}}]$ and $\EE[\partial_y f_{\text{e}}]$ vanish identically. Therefore, the requirement $c_1 = c_2 = 0$ necessitates $\EE[\partial_x f_{\text{o}}] = \EE[\partial_y f_{\text{o}}] = 0$ for all covariance matrices $\Sigma \succ 0$.

Let $g(\xx) = \partial_x f_{\text{o}}(x, y)$. Note that $g(\xx)$ is even. The condition $\EE[g(\xx)] = 0$ implies that, for any $A = \Sigma^{-1}$,
\al{
	\int_{\RR^2}\rmd\xx\; g(\xx) 
	\exp\left(-\frac{1}{2} \xx^\top A \xx\right) = 0\,.
}
By \eqref{eq:activation_increment_bound}, $g$ has at most polynomial growth. Hence, in a sufficiently small neighborhood of $A=\1_2$, the integrand and all its derivatives with respect to $A_{11}$, $A_{22}$, and $A_{12}$ are dominated by integrable functions. We may therefore differentiate under the integral sign arbitrarily many times and then set $A=\1_2$, which yields
\al{
	\int_{\RR^2}\rmd\xx\; g(\xx) x^{2k} y^{2m} (xy)^\ell 
	\rme^{- (x^2 + y^2)/2} = 0
	\label{eq:moments_vanishing}
}
for all nonnegative integers $k, m, \ell \geq 0$. Equation~\eqref{eq:moments_vanishing} implies that $g(\xx)$ is orthogonal, under the standard bivariate Gaussian measure $\gamma$, to every monomial $x^a y^b$ of even total degree $a+b$. Furthermore, since $g(\xx)$ itself is an even function, the integral of $g(\xx) x^a y^b$ is trivially zero whenever $a+b$ is odd. Consequently, $g(\xx)$ is orthogonal to all polynomials in the Hilbert space $L^2(\RR^2, \gamma)$. Since the tensor products of Hermite polynomials form a complete orthogonal basis for $L^2(\RR^2, \gamma)$, polynomials are dense in this space. Thus, $g(\xx) = \partial_x f_{\text{o}}(x, y) = 0$. 

By an identical argument, $\partial_y f_{\text{o}}(x, y) = 0$. This implies that $f_{\text{o}}(h)$ must be a constant. Since $f_{\text{o}}(h)$ is odd by definition, this constant is zero. We thus conclude that $f(h) = f_{\text{e}}(h)$, which establishes that $f(h) = f(-h)$.
\end{proof}
	
We next ask when the conjugate-linear response $c_2$ vanishes identically over all admissible variances and pseudo-variances. Although holomorphicity of $f$ is sufficient, the following lemma shows that it is not necessary.
\begin{lemma}
\label{lem:c2zero}
One has $c_2=0$ for every $v_h>0$ and $\tau_h\in\CC$ $(|\tau_h| < v_h)$ if and only if the activation function can be decomposed as
\al{
	f(h) = f_{\rm odd,\, holo}(h) + f_{\rm even}(h)\,, 
	\label{eq:fdf9fds}
}
where $f_{\rm odd,\, holo}(h)$ is a holomorphic odd function $(f_{\rm odd,\, holo}(-h) = -f_{\rm odd,\, holo}(h))$ and $f_{\rm even}(h)$ is an even function $(f_{\rm even}(-h) = f_{\rm even}(h))$.
\end{lemma}
\begin{proof}
Let $c_2[f]$ denote the coefficient $c_2$ viewed as a linear functional of $f$.

\emph{Sufficiency:} By the linearity of the functional $c_2$, we have $c_2[f_{\text{odd,\,holo}} + f_{\text{even}}] = c_2[f_{\text{odd,\,holo}}] + c_2[f_{\text{even}}]$. Recalling \eqref{eq:stein_gain_2} and using $\partial_{\bar{h}} f_{\text{odd,\,holo}}(h) = 0$, we trivially obtain $c_2[f_{\text{odd,\,holo}}] = 0$. Meanwhile, Lemma~\ref{lem:bothczero} guarantees that $c_2[f_{\text{even}}] = 0$. Thus, $c_2[f] = 0$.

\emph{Necessity:} Suppose $c_2[f] = 0$ for all valid covariance matrices. We uniquely decompose $f$ into its even and odd components: $f(h) = f_{\text{e}}(h) + f_{\text{o}}(h)$. By the linearity of $c_2$, we have
\al{
	c_2[f] = c_2[f_{\text{e}}] + c_2[f_{\text{o}}] 
	= c_2[f_{\text{o}}] \,,
}
where Lemma~\ref{lem:bothczero} was utilized. Thus the condition $c_2[f] = 0$ implies $c_2[f_{\text{o}}] = \EE[\partial_{\bar{h}} f_{\text{o}}] = 0$.

Since $f_{\text{o}}(h)$ is an odd function, its derivative $\partial_{\bar{h}} f_{\text{o}}(h)$ is an even function. As established in the proof of Lemma~\ref{lem:bothczero}, if the expectation of an even function uniformly vanishes under any zero-mean Gaussian measure, the function itself must be zero. Consequently, $\partial_{\bar{h}} f_{\text{o}}(h) = 0$. Local Lipschitz continuity then implies that $f_{\mathrm{o}}$ is holomorphic.

Setting $f_{\text{odd,\,holo}} = f_{\text{o}}$ and $f_{\text{even}} = f_{\text{e}}$ immediately yields the required decomposition \eqref{eq:fdf9fds}.
\end{proof}
 
By the same token, $c_1$ uniformly vanishes for arbitrary $v_h$ and $\tau_h$ if and only if $f$ can be decomposed as $f(h) = f_{\text{odd,\,anti-holo}}(h) + f_{\text{even}}(h)$, where $f_{\text{odd,\,anti-holo}}(h)$ is an anti-holomorphic odd function. 

As an example for Lemma~\ref{lem:c2zero}, consider $f(h)=h^3$, which is an entire function. One can verify that the surrogate becomes $f(H)\approx 3\tau_h H + Z$ with no $\bar{H}$ component, where $Z_{ij}\sim \CN\big(0,6v_h^3,6\tau_h^3 \big)$. 
A more nontrivial example is the complex cardioid activation \autocite{Virtue2017}:
\al{
	f(h) = \frac{1+\cos(\angle h)}{2}h \,,
}
where $\angle h$ denotes the argument of $h$. This $f$ can be decomposed as
\al{
	f(h) = \frac{h}{2} + \frac{h+\bar{h}}{4|h|} h \,.
}
Since this matches the decomposition \eqref{eq:fdf9fds}, it follows from Lemma~\ref{lem:c2zero} that $c_2=0$.

Next, consider a split-type activation function $f(h)=g(x)+ ig(y)$, in which the same real function $g:\RR\to\RR$ is applied to both the real and imaginary parts of the input. This class of activations, which belongs to Type~A activations in the terminology of \textcite{Kuroe2003}, is widely used in the literature \autocite{Bassey2021review,Lee2022review,Barrachina2023}.
\begin{lemma}\label{lem:spg}
For a split activation $f(h)=g(x)+ ig(y)$, one has $c_2=0$ for all $v_h >0$ and $\tau_h\in\CC$ $(|\tau_h| < v_h)$ if and only if $g$ can be decomposed as
\al{
	g(x) = ax + g_{\rm even}(x)\,, 
	\label{eq:1gdecompmm}
}
where $a\in\RR$ is a constant and $g_{\rm even}(x)$ is an even function of $x$.
\end{lemma}
\begin{proof}
\emph{Sufficiency}: Assuming that $g$ is of the form \eqref{eq:1gdecompmm}, we have $f(h)=ah + g_{\rm even}(x) + ig_{\rm even}(y)$. Since this is a sum of a holomorphic odd function and an even function, it readily follows from Lemma~\ref{lem:c2zero} that $c_2=0$.

\emph{Necessity:} Assume that $c_2 = 0$ uniformly holds for any valid covariance structure. According to Lemma~\ref{lem:c2zero}, the activation function must be of the form
\al{
	f(h) & = f_{\rm odd,\,holo}(h) + f_{\rm even}(h)\,.
} 
To connect this with the split structure, we decompose $g(x)$ into its odd and even parts as $g(x) = g_{\rm o}(x) + g_{\rm e}(x)$. Then
\al{
	f(h) = [g_{\rm o}(x) + ig_{\rm o}(y)] 
	+ [ g_{\rm e}(x) + ig_{\rm e}(y)]\,.
}
By uniqueness of the parity decomposition, we must have $f_{\rm odd,\,holo}(h) = g_{\rm o}(x) + ig_{\rm o}(y)$, which means that $g_{\rm o}(x) + ig_{\rm o}(y)$ is holomorphic. Applying the Cauchy--Riemann equations yields $g_{\rm o}'(x) = g_{\rm o}'(y)$. Since $x$ and $y$ are independent variables, this equality implies that the derivative must be a constant, say $g_{\rm o}'(x) = a$ for some $a \in \RR$. Integrating this equation yields $g_{\rm o}(x) = ax$, where the constant of integration was put to zero due to the oddness of $g_{\rm o}$. Consequently, we established the required decomposition $g(x) = ax + g_{\rm e}(x)$.
\end{proof}

An analogous argument shows that, for an activation $f(h)=g(x)-ig(y)$ with $g:\RR\to\RR$, the coefficient $c_1$ vanishes uniformly for arbitrary $v_h$ and $\tau_h$ if and only if $g$ takes the form \eqref{eq:1gdecompmm}.

An important example for Lemma~\ref{lem:spg} is given by the split-ReLU activation \autocite{Trabelsi2018,Gao2019}, defined by $g(x)=\max(x, 0)$. Since the latter can be decomposed as $g(x)=x/2 + |x|/2$ in line with \eqref{eq:1gdecompmm}, it follows from Lemma~\ref{lem:spg} that $c_2=0$. It is also easy to verify that $c_1=1/2$. 

More generally, $c_2=0$ can be proven not only for Split-ReLU but also for the split versions of Leaky-ReLU \autocite{Maas2013}, GELU \autocite{Hendrycks2016}, and SiLU \autocite{Elfwing2018} (also known as Swish \autocite{Ramachandran2017}). Indeed, it is easy to decompose them as in \eqref{eq:1gdecompmm}:
\al{
	\operatorname{Leaky-ReLU}(x) 
	& = \max(x,0) - a\max(-x,0) = \frac{1+a}{2} x 
	+ \frac{1-a}{2}|x|\,,
	\\
	\operatorname{GELU}(x) & = x\Phi(x)
	= \frac{x}{2} + x \kkakko{\Phi(x)-\frac{1}{2}} \,,
	\\
	\operatorname{SiLU}(x) & = x\sigma(x)
	= \frac{x}{2} + x \kkakko{\sigma(x)-\frac{1}{2}} \,,
}
where $\Phi$ and $\sigma$ denote the standard normal cumulative distribution function and the sigmoid function, respectively. Note that both $\Phi(x)-1/2$ and $\sigma(x)-1/2$ are odd functions.

\subsubsection{Proper inputs}

We now turn to the proper case $(\tau_h=0)$, where the input $h$ obeys the circularly symmetric complex Gaussian distribution $p(h)=\frac{1}{\pi v_h} \exp\left(-|h|^2 /v_h\right)$.

\begin{lemma}\label{lem:propercczero}
Suppose $\tau_h = 0$. One has $c_1=0$ for every $v_h>0$ if and only if the Fourier coefficient of order $1$ of $f$ vanishes. Equivalently, writing $h = r \rme^{i\theta}$,
\al{
	\int_0^{2\pi} \!\!\!\!\!\rmd \theta ~
	f(r \rme^{i\theta}) \rme^{-i\theta} = 0\,,
	\qquad \text{$\forall r > 0$\,.}
}
Similarly, one has $c_2=0$ for every $v_h>0$ if and only if the Fourier coefficient of order $-1$ of $f$ vanishes. Equivalently,
\al{
	\int_0^{2\pi} \!\!\!\!\!\rmd \theta ~
	f(r \rme^{i\theta}) \rme^{i\theta} = 0\,,
	\qquad \text{$\forall r > 0$\,.}
	\label{eq:sd789879fgdf}
}
\end{lemma}
\begin{proof}
From \eqref{eq:kaka} evaluated at $\tau_h = 0$, the condition $c_1 = 0$ reduces to $\kappa_1 = \EE[f(h)\bar{h}] = 0$.
In polar coordinates, this condition is expressed as
\al{
	\kappa_1 = \frac{1}{\pi v_h} 
	\int_0^\infty \!\!\! \rmd r\, r 
	\int_0^{2\pi}\!\!\!\!\! \rmd \theta~
	f(r \rme^{i\theta})\; r \rme^{-i\theta} 
	\exp\left(-\frac{r^2}{v_h}\right) = 0 \,.
}
Let $U_1(r) = \frac{1}{2\pi} \int_0^{2\pi} \!\! \rmd \theta\;
f(r \rme^{i\theta}) \rme^{-i\theta}$ denote the Fourier coefficient of order $1$ of $f$ with respect to $\theta$. By substituting $U_1(r)$, the condition becomes
\al{
	\int_0^\infty \!\!\! \rmd r\; r^2 \, 
	U_1(r) \exp\left(-\frac{r^2}{v_h}\right) = 0\,,
	\qquad \forall v_h > 0 \,.
}
Applying the change of variables $t = r^2$ and defining $s = 1/v_h > 0$, the integral can be recast as
\al{
	\int_0^\infty \!\!\!\rmd t\; 
	\sqrt{t}\; U_1(\sqrt{t}) \rme^{-st} = 0 \,, 
	\qquad \forall s > 0 \,.
}
This is the Laplace transform of the function $t\mapsto\sqrt{t}\,U_1(\sqrt{t})$. Since it vanishes for every $s>0$, Lerch's theorem \autocite{Cohen2007} implies that $U_1(r)=0$ for almost every $r>0$. The continuity of $f$ implies that $U_1$ is continuous, and hence $U_1(r)=0$ for every $r>0$.

By a completely symmetric argument, the condition $c_2 = 0$ reduces to $\kappa_2 = \EE[f(h)h] = 0$. 
Extracting the angular integral for $\rme^{i\theta}$, we obtain the requirement
\al{
	U_{-1}(r) = \frac{1}{2\pi} \int_0^{2\pi} \!\!\!\!\! 
	\rmd \theta ~f(r \rme^{i\theta}) \rme^{i\theta} = 0\,, 
	\qquad \forall r>0\,.
}
This completes the proof.
\end{proof}

It is worth noting that the condition for $c_1 = c_2 = 0$ is significantly relaxed under the proper Gaussian assumption compared to the general improper case in Lemma~\ref{lem:bothczero}. While Lemma~\ref{lem:bothczero} strictly requires the activation function to be even, meaning that all odd-order Fourier coefficients must vanish, the proper case merely demands the absence of the $\pm 1$-order coefficients. For example, all activation functions of the form
\al{
	f(h) & = a_0(r) + \sum_{k=2}^{\infty}\ckakko{
		a_k(r) h^k + \tilde{a}_k(r) \bar{h}^k
	}
}
for arbitrary functions $a_0$, $a_k$, $\tilde{a}_k$ of the modulus $r=|h|$, satisfy the conditions of Lemma~\ref{lem:propercczero} and yield $c_1=c_2=0$. In particular, monomials $f(h)=h^k$ and $\bar{h}^k$ with $k\geq 2$ belong to this class.

A popular class of non-split-type activations is given by $f(h) = u(r)\rme^{i\theta}$ with $h \equiv r\rme^{i\theta}$ \autocite{Hirose1992,Georgiou1992}, which belongs to Type~B activation functions in the terminology of \textcite{Kuroe2003}. The widely used modReLU activation \autocite{Arjovsky2016}, defined by $u(r)=\max(r+b, 0)$, belongs to this class. Under phase rotations, $f(h)$ transforms in the same way as $h$, which is a favorable property for some applications.%
\footnote{A shallow and sufficiently wide CVNN using such an activation is guaranteed to retain the universal approximation property for the class of all $\mathrm{U}(1)$-equivariant functions \autocite[Proposition 2.4]{Yarotsky2022}.} 
\begin{lemma}\label{lem:c2_0_for_f_with_phase}
Suppose $\tau_h=0$. For an activation of the form $f(h) = u(r)\rme^{i\theta}$, one has $c_2=0$ and $\tau_z = 0$ for every $v_h>0$.
\end{lemma}
\begin{proof}
Since $f$ does not possess the Fourier coefficient of order $-1$, Lemma~\ref{lem:propercczero} implies that the coefficient $c_2$ uniformly vanishes for every $v_h>0$. 

Plugging $\tau_h=c_2=0$ into \eqref{eq:kaka} gives $\kappa_2=0$. Equation \eqref{eq:RTEGD1F} evaluated at $c_2=\kappa_2=0$ yields $\tau_z = \EE[(f(h)-\mu)^2]=\EE[f(h)^2]$, where we used $\mu=\EE[f(h)]\propto \int_0^{2\pi} \!\!\rmd \theta \; \rme^{i\theta} = 0$. Moreover, $\EE[f(h)^2]\propto \int_0^{2\pi} \!\!\rmd \theta \; \rme^{2i\theta}=0$. This proves $\tau_z=0$, i.e., the residual noise $Z$ is proper.
\end{proof}

For split-type activations, one can show the following strong result for proper inputs.

\begin{lemma}\label{lem:ooqwQOWc}
Suppose $\tau_h=0$. For a split activation $f(h) = g(x) + ig(y)$, $c_2=0$ for every $v_h>0$.
\end{lemma}

\begin{proof}
Employing the Wirtinger derivative with respect to $\bar{h}$ yields
\al{
	\partial_{\bar{h}} f(h) 
	= \frac{1}{2}(\partial_x + i\partial_y)\big(g(x) + ig(y)\big) 
	= \frac{1}{2}\big(g'(x) - g'(y)\big) \,.
}
Recalling \eqref{eq:stein_gain_2}, we obtain
\al{
	c_2 = \EE[\partial_{\bar{h}}f(h) ]
	= \frac{1}{2}\big( \EE[g'(x)] - \EE[g'(y)] \big) \,.
	\label{eq:345cscERTRE}
}
When $\tau_h = 0$, the real and imaginary parts $x$ and $y$ are uncorrelated and share the exact same marginal zero-mean Gaussian distribution, $\RN(0, v_h/2)$. Therefore, exact cancellation occurs in the right-hand side of \eqref{eq:345cscERTRE}, yielding $c_2 = 0$.

Alternatively, the same result can be obtained by directly verifying the condition \eqref{eq:sd789879fgdf} in Lemma~\ref{lem:propercczero} for the split activation, though the explicit integration is omitted here for brevity.
\end{proof}

Unlike Lemma~\ref{lem:spg}, Lemma~\ref{lem:ooqwQOWc} imposes no additional structural condition on $g$.

\section{Spectral Complex GEP}
\label{sc:gep_rigorous}

In this section, we formulate and prove a spectral version of the general complex GEP proposed in Section~\ref{sc:gengep}.

\begin{theorem}[Spectral complex Gaussian equivalence]
\label{thm:complex_gep}
Under the assumptions of Section~\ref{sc:notation} and the proportional scaling \eqref{eq:lmlm}, suppose in addition that at least one of $X$ and $W$ is proper, that is, $\tau_x=0$ or $\tau_w=0$. Define the Gaussian surrogate
\al{
	Q \equiv c_1 H + c_2 \bar{H} + Z\,,
	\label{eq:Q_rigorous}
}
where $Z$ is independent of $H$ and has \mbox{i.i.d.} entries distributed as $\CN(0,v_z,\tau_z)$. The parameters $c_1,c_2,v_z$ and $\tau_z$ are defined in Section~\ref{sc:widely_linear}.

Let $\bar\nu_f^{(n,d,p)}$ and $\bar\nu_Q^{(n,d,p)}$ denote the expected empirical spectral distributions of $f(H)^*f(H)/n$ and $Q^*Q/n$, respectively. Then there exists a compactly supported probability measure $\nu$ on $\RR_{\geq0}$ such that
\al{
	\bar\nu_f^{(n,d,p)} \longrightarrow\nu\,, \qquad
	\bar\nu_Q^{(n,d,p)} \longrightarrow\nu
	\label{eq:expected_esd_equivalence}
}
weakly as $n,d,p\to\infty$.
\end{theorem}

We regard the requirement that at least one factor be proper as a technical condition of the present proof rather than an intrinsic limitation of the complex GEP. The limiting spectrum of the Gaussian surrogate \eqref{eq:Q_rigorous} is studied separately in our companion paper \autocite{kanazawa_rmt_2026} without assuming propriety of either $X$ or $W$. There, a system of algebraic self-consistent equations for the Stieltjes transform of $Q^*Q/n$ is derived.

The remainder of this section first introduces the intermediate models and the proof strategy (Section~\ref{sc:spectral_gep_overview}). We then establish the limiting spectral properties of the Gaussian surrogate (Lemma~\ref{lem:surrogate_limit}), compare the surrogate with the nonlinear model (Propositions~\ref{prop:approximation_steps} and~\ref{prop:finite_degree_equivalence}), and finally complete the proof of Theorem~\ref{thm:complex_gep} (Section~\ref{sc:spectral_gep_proof}).

\subsection{Outline of the Proof}
\label{sc:spectral_gep_overview}

The proof is written for the case in which $X$ is proper, without loss of generality. If only $W$ is proper, we transpose the model: $H^\top = W^\top X^\top /\sqrt{d}$. The transposed model has sample size $p$ and feature dimension $n$, with $W^\top$ now playing the role of the proper data matrix. Since the nonlinearity is applied entrywise, $f(H^\top) = f(H)^\top$. For any $A\in\CC^{n\times p}$, the nonzero eigenvalues of $A^*A/n$ are obtained by multiplying those of $\bar A A^\top/p$ by $p/n$. The remaining eigenvalues are zero, with multiplicities determined only by $n$ and $p$. The same relation holds for the nonlinear matrix and its surrogate. Therefore, spectral equivalence for the transposed model implies spectral equivalence for the original model.

We first remove the mean component from the nonlinear feature matrix. Subtracting $\mu1_n1_p^\top$ from $f(H)$ changes its empirical covariance matrix by a perturbation of rank at most two and therefore does not affect its limiting expected empirical spectral distribution (see \textcite[Theorem~A.44]{Bai2010book} and \textcite[Lemma~2.9]{Couillet_Liao_2022}). We may thus prove the theorem for the centered activation $f-\mu$. For notational simplicity, we retain the notation $f$ for this centered activation throughout the proof.

The next step is to replace the original preactivation matrix $H$ by its column-normalized version $\HH=HS$ (recall \eqref{eq:column-normalization} and \eqref{eq:feature-matrices} for the definitions of $S$ and $\HH$). When $X$ is proper, conditioning on $W$ gives
\al{
	\HH_{i\alpha}\mid W \sim \CN(0,v_h,0)
	\label{eq:HH_dist_exact}
}
exactly at finite dimensions for every $(i,\alpha)$. This common Gaussian marginal allows us to use a single complex Hermite expansion for all entries of $\HH$.

Writing $h\sim\CN(0,v_h,0)$, we expand the centered activation as
\al{
	f(h) = \sum_{s,t\geq0}
	a_{s,t}\, \He_{s,t}(h,\bar{h})\,,
	\qquad a_{0,0} = 0 \,,
}
where $\{\He_{s,t}\}$ denotes the complex Hermite polynomials \autocite{Ito1952}; see Appendix~\ref{app:hermite} for details. For a general nonlinear activation, this expansion contains infinitely many terms. For technical transparency, we first establish the spectral equivalence for a finite Hermite truncation and then remove the truncation in a separate step. Specifically, for an integer $M\geq1$, define
\al{
	f_M(h) \equiv \sum_{1\leq s+t \leq M}
	a_{s,t}\, \He_{s,t}(h,\bar{h})\,.
	\label{eq:complex_hermite_truncation}
}
Let $c_{1,M}$, $c_{2,M}$, $v_{z,M}$, and $\tau_{z,M}$ be the surrogate parameters associated with $f_M$ through the formulas of Section~\ref{sc:widely_linear}. We introduce a Gaussian matrix $Z_M$, independent of $\HH$ and $H$, with \mbox{i.i.d.} entries
\al{
	(Z_M)_{i\alpha} \sim
	\CN(0,v_{z,M},\tau_{z,M})
}
and define
\al{
	\hat{Q}_M \equiv c_{1,M}\HH
	+ c_{2,M}\bar{\HH} + Z_M\,.
	\label{eq:Qhat_M_definition}
}
The column-normalized surrogate associated with the full activation is
\al{
	\hat{Q} \equiv c_1\HH + c_2\bar{\HH} + Z\,.
}

The complete comparison is summarized by the chain
\al{
	Q \longleftrightarrow
	\hat Q \longleftrightarrow
	\hat Q_M \longleftrightarrow
	f_M(\HH) \longleftrightarrow
	f(\HH) \longleftrightarrow f(H)\,.
	\label{eq:proof_chain}
}
The proof combines the central finite-degree comparison with four approximation steps connecting the truncated models to the original surrogate and nonlinear feature matrices. The four approximation steps are controlled using the Frobenius-norm bound for normalized Stieltjes transforms established in Lemma~\ref{lem:stieltjes_comparison} of Appendix~\ref{app:frobenius}, whereas the finite-degree comparison is proved by a graphical moment argument. The main steps are summarized below.

\begin{enumerate}
	\item
	We first establish the compactly supported limiting expected spectral measure $\nu$ of $Q^*Q/n$.
	
	\item
	The comparison $Q\longleftrightarrow\hat Q$ introduces the column normalization on the surrogate side. Since the difference is linear in $\HH-H$, the result follows from the concentration of $S$ around $\1_p$ and Lemma~\ref{lem:stieltjes_comparison}.
	
	\item
	The comparison $\hat Q\longleftrightarrow\hat Q_M$ replaces the full surrogate by the one associated with the finite-degree Hermite truncation. It follows from the convergence of the scalar surrogate parameters as $M\to\infty$, together with the construction of $Z_M$ and $Z$ from the same underlying Gaussian variables.
	
	\item
	The central comparison $\hat Q_M\longleftrightarrow f_M(\HH)$ is the combinatorial core of the proof. For fixed $M$, $f_M$ is a finite-degree polynomial in $h$ and $\bar h$, and a graphical expansion shows that the two models have the same limiting expected spectral moments. Moment determinacy then transfers the compactly supported limiting measure of $\hat Q_M^*\hat Q_M/n$ to
	$f_M(\HH)^*f_M(\HH)/n$.
	
	\item
	The comparison $f_M(\HH)\longleftrightarrow f(\HH)$ removes the Hermite cutoff. It follows from the Gaussian $L^2$ convergence of $f_M$ to $f$ and Lemma~\ref{lem:stieltjes_comparison}.
	
	\item
	Finally, the comparison $f(\HH)\longleftrightarrow f(H)$ removes the column normalization on the nonlinear side. This comparison follows from the increment bound on $f$ in \eqref{eq:activation_increment_bound}, the concentration of the column norms of $W$ and Lemma~\ref{lem:stieltjes_comparison}.
\end{enumerate}
The high-dimensional limit is always taken first with $M$ fixed, followed by removal of the Hermite cutoff:
\al{
	M\ \text{fixed} \qquad\longrightarrow\qquad
	n,d,p\to\infty \qquad\longrightarrow\qquad M\to\infty\,.
}
Thus, spectral moments are used only for the comparison of the two finite-degree models in the center of \eqref{eq:proof_chain}; all remaining arrows are handled at the level of normalized Stieltjes transforms.

\subsection{Auxiliary Results I}

\begin{lemma}[Limiting spectral measure of the Gaussian surrogate]
	\label{lem:surrogate_limit}
	Let $Q=c_1H+c_2\bar H+Z$. Then the expected empirical spectral distribution of $Q^*Q/n$ converges weakly to a compactly supported probability measure $\nu$ on $\RR_{\geq0}$. Moreover, for every fixed integer $k\geq1$,
	\al{
		\frac{1}{p}\EE\bigg[\Tr
		\bigg\{\left(\frac{Q^*Q}{n}\right)^k\bigg\}\bigg]
		\longrightarrow
		\int_{\RR_{\geq0}}\lambda^k\,\rmd\nu(\lambda)
		\equiv \moment_k\,.
	}
	In particular, $\nu$ is uniquely determined by the moment sequence $\{\moment_k\}_{k\geq0}$, with $\moment_0=1$.
\end{lemma}

\begin{proof}
Fix $k\geq1$. Expanding $(Q^*Q)^k$ and substituting $H=XW/\sqrt d$ gives a finite sum of Gaussian moments. By Wick's rule, each Gaussian moment is a finite sum over pairings, with the resulting Kronecker deltas identifying the corresponding indices. For example, one contraction contributing to the pure-noise word for $k=3$ is
\al{
	\frac{1}{p}\EE\bigg[
		\Tr\bigg\{\mkakko{\frac{Z^*Z}{n}}^3\bigg\}
	\bigg]
	&=\frac{1}{pn^3}
	\sum_{\substack{i,j,m,\\ \alpha,\beta,\gamma}}
	\EE\bigg[
		\contraction{}{
			\bar{Z_{i\alpha}}
		}{
			Z_{i\beta}\bar{Z_{j\beta}}
		}{
			Z_{j\gamma}
		}
		\contraction[2.5ex]{
			\bar{Z_{i\alpha}}
		}{
			Z_{i\beta}
		}{}{
			\bar{Z_{j\beta}}
		}
		\contraction{
			\bar{Z_{i\alpha}}Z_{i\beta} \bar{Z_{j\beta}}Z_{j\gamma}
		}{
			\bar{Z_{m\gamma}}
		}{}{
			Z_{m\alpha}
		}
		\bar{Z_{i\alpha}} Z_{i\beta} \bar{Z_{j\beta}} Z_{j\gamma}
		\bar{Z_{m\gamma}} Z_{m\alpha}
	\bigg] + \cdots
	\\
	&=\frac{1}{pn^3} 	\sum_{\substack{i,m,\alpha,\beta}}
	\EE[\bar{Z_{i\alpha}} Z_{i\alpha}]
	\EE[\bar{Z_{i\beta}} Z_{i\beta}]
	\EE[\bar{Z_{m\alpha}} Z_{m\alpha}]
	+ \cdots.
}
The displayed contraction, shown in Figure~\ref{fg:zzzz}, leaves the four free indices $i,m,\alpha,\beta$. Its index volume is therefore of order $n^2p^2$, which exactly balances the normalization $pn^3$ in the proportional regime.

\begin{figure}[tb]
	\centering
	\includegraphics[width=14em]{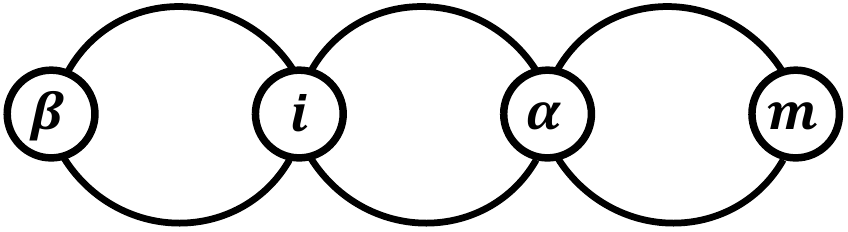}
	\caption{\label{fg:zzzz}A Wick contraction of a product of six Gaussian variables. The four nodes represent the free indices $i$, $m$, $\alpha$ and $\beta$, and each edge represents one Gaussian factor. Each pair of parallel edges represents one Wick pair. The three Wick pairs connect the index pairs $(i,\alpha)$, $(i,\beta)$ and $(m,\alpha)$.}
\end{figure}

We next extend this counting argument to a general contraction. Consider a word in the expansion of $(Q^*Q)^k$ containing $r$ factors of $H$, $\bar H$, $H^*$, or $H^{\mathsf T}$. A nonzero Gaussian moment requires the remaining number $2k-r$ of $Z$-type factors to be even; hence $r=2\ell$ for some $0\leq\ell\leq k$. After inserting $H=XW/\sqrt d$ and its conjugates, the word carries the normalization $d^{-\ell}$ and contains $2r+(2k-r)=2k+2\ell$ Gaussian factors.

Fix a Wick contraction $\pi$, and let $V_\pi$ be the total number of free indices remaining after all index identifications. Represent each Gaussian factor by an edge joining its two index nodes. The resulting contraction graph has $V_\pi$ nodes and $2k+2\ell$ edges. It is connected: before contraction, the factors belong to a single trace and form a connected cyclic graph, and the index identifications introduced by the contraction cannot disconnect it.

Each Wick pair produces two parallel edges joining the same pair of identified nodes. We construct a simple graph by replacing each nonempty bundle of parallel edges by a single edge, and let $q_\pi$ denote its number of edges. This graph remains connected and has the same $V_\pi$ nodes. Since each of its edges represents at least two of the original $2k+2\ell$ edges,
\al{
	2q_\pi\leq 2k+2\ell,
	\qquad\text{so that}\qquad
	q_\pi\leq k+\ell.
}
On the other hand, a connected graph with $V_\pi$ nodes must have at least $V_\pi-1$ edges. Hence
\al{
	V_\pi-1\leq q_\pi\leq k+\ell
	\qquad \Longrightarrow \qquad 
	V_\pi\leq k+\ell+1.
	\label{eq:free_index_bound}
}

Let $N$ denote a common dimensional scale, so that $n\asymp p\asymp d\asymp N$. The normalization factor of a word with $2\ell$ $H$-type factors is $(pn^k)^{-1}(1/\sqrt{d})^{2\ell}\asymp N^{-k -\ell -1}$, whereas a contraction with $V_\pi$ free indices has index volume of order $N^{V_\pi}$. Hence its contribution is of order $N^{V_\pi-(k+\ell+1)}$. By \eqref{eq:free_index_bound}, no contraction grows with the dimension. Contractions satisfying the strict inequality vanish, whereas the leading contractions converge to finite constants. Since the number of words and Wick contractions is finite for fixed $k$, the normalized expected trace converges to a finite limit, denoted by $\moment_k$.

Let $\bar\nu_Q^{(n,d,p)}$ denote the expected empirical spectral distribution of $Q^*Q/n$. We next use an operator-norm bound to establish tightness and compactness of its possible weak limits. Since $H=XW/\sqrt d$ and $\|\bar H\|_{\rm op}=\|H\|_{\rm op}$,
\al{
	\frac{\|Q\|_{\rm op}}{\sqrt n}
	&\leq (|c_1|+|c_2|)
	\frac{\|H\|_{\rm op}}{\sqrt n}
	+\frac{\|Z\|_{\rm op}}{\sqrt n}
	\\
	&\leq (|c_1|+|c_2|)
	\frac{\|X\|_{\rm op}}{\sqrt n}
	\frac{\|W\|_{\rm op}}{\sqrt d}
	+\frac{\|Z\|_{\rm op}}{\sqrt n}\,.
}
The Bai--Yin theorem (see \autocite{YinBaiKrishnaiah1988}, \autocite[Theorem~5.8]{Bai2010book}, and \autocite[Theorem~2.11]{Couillet_Liao_2022}) shows that the three normalized norms on the right-hand side converge almost surely to finite deterministic limits in the proportional regime. Hence there exists a deterministic constant $\Gamma<\infty$, independent of $n,d,p$, such that
\al{
	\PP\bigg(
	\operatorname{spec}\mkakko{\frac{Q^*Q}{n}}
	\not\subseteq[0, \Gamma]
	\bigg)\longrightarrow 0
	\qquad \text{as~~$n,d,p\to\infty$.}
}
It follows that
\al{
	\bar\nu_Q^{(n,d,p)}([0, \Gamma]^c)
	\leq \PP\bigg(
	\operatorname{spec}\mkakko{\frac{Q^*Q}{n}}
	\not\subseteq[0, \Gamma]
	\bigg)
	\longrightarrow 0
	\qquad \text{as~~$n,d,p\to\infty$.}
}
This establishes that $\big\{\bar\nu_Q^{(n,d,p)}\big\}$ is tight. By Prokhorov's theorem \autocite{Billingsley1999}, every subsequence contains a further subsequence that converges weakly to some probability measure $\tilde\nu$. Moreover, every such subsequential limit is supported on $[0,\Gamma]$. Indeed, since
\al{
	\bar\nu_Q^{(n,d,p)}\bigl([0, \Gamma]^c\bigr)
	\longrightarrow 0 \,,
}
the Portmanteau theorem \autocite{Billingsley1999} gives
\al{
	\tilde\nu\bigl([0, \Gamma]^c\bigr) 
	\leq \liminf_{j\to\infty}\bar\nu_Q^{(n_j, d_j, p_j)}
	\bigl([0, \Gamma]^c\bigr) = 0 \,.
}
It remains to identify $\tilde\nu$. The Wick estimates above apply at every fixed order. In particular, the estimate at order $k+1$ gives 
\al{
	\sup_j \int_{\RR_{\geq 0}} \lambda^{k+1}\, 
	\rmd\bar\nu_Q^{(n_j, d_j, p_j)}(\lambda) < \infty \,. 
}
It follows that, for every $R>0$,
\al{
	\sup_j \int_{(R,\infty)} \lambda^k\, 
	\rmd\bar\nu_Q^{(n_j, d_j, p_j)}(\lambda) 
	\leq \frac{1}{R} \sup_j 
	\int_{\RR_{\geq 0}} \lambda^{k+1}\, 
	\rmd\bar\nu_Q^{(n_j, d_j, p_j)}(\lambda)
	= O\big(R^{-1}\big) \,.
}
The right-hand side can be made arbitrarily small by taking $R$ sufficiently large. Thus, the functions $\lambda^k$ are uniformly integrable along the subsequence. Weak convergence then gives
\al{
	\int_{\RR_{\geq 0}} \lambda^k\,
	\rmd \tilde\nu(\lambda) = \lim_{j\to\infty}
	\int_{\RR_{\geq 0}} \lambda^k\, 
	\rmd\bar\nu_Q^{(n_j, d_j, p_j)}(\lambda) = \moment_k \,.
}
Therefore, every weak subsequential limit has the moment sequence $\{\moment_k\}_{k\geq0}$, where $\moment_0=1$.

Since every such limit is supported on the common compact interval $[0, \Gamma]$, it is uniquely determined by its moments \autocite{Lin2017review}. Hence all weak subsequential limits coincide with a single compactly supported probability measure $\nu$. It follows that the full sequence $\big\{\bar\nu_Q^{(n,d,p)}\big\}$ converges weakly to $\nu$, whose moments are precisely $\{\moment_k\}_{k \geq 0}$.
\end{proof}

\begin{proposition}[Truncation and normalization errors]
\label{prop:approximation_steps}
For every fixed $z\in\CC_+$, the following four comparisons hold, with all dimensional limits taken under \eqref{eq:lmlm}:
\begin{gather}
	\EE\big[\big|m_Q(z)-m_{\hat{Q}}(z)\big|\big] \longrightarrow 0 \,,
	\label{eq:closure_Q_Qhat}
	\\
	\lim_{M\to\infty} \limsup_{n,d,p\to\infty}
	\EE\big[\big| m_{\hat{Q}}(z)-m_{\hat{Q}_M}(z) \big| \big]= 0\,,
	\label{eq:closure_Q_QhatM}
	\\
	\lim_{M\to\infty} \limsup_{n,d,p\to\infty}
	\EE\big[\big| m_{f_M(\HH)}(z)-m_{f(\HH)}(z) \big| \big] = 0\,,
	\label{eq:closure_fM_f}
	\\
	\EE\big[\big| m_{f(\HH)}(z)-m_{f(H)}(z) \big|\big] \longrightarrow 0\,.
	\label{eq:closure_Hhat_H}
\end{gather}
\end{proposition}

\begin{proof}
We apply Lemma~\ref{lem:stieltjes_comparison} to each of the four comparisons. It therefore suffices to show, for each relevant pair of matrices $A$ and $B$, that
\al{
	\frac{1}{np}\EE \big[\|A-B\|_{\rm F}^2 \big] \longrightarrow 0
}
with the order of limits stated in the proposition, and that
\al{
	\frac{1}{np}\EE\big[\|A\|_{\rm F}^2+\|B\|_{\rm F}^2\big]
}
remains bounded. We verify these conditions in turn.

\paragraph{Step 1: \boldmath $Q\longleftrightarrow\hat Q$.}
Recall that $\HH=HS$, where
\al{
	S = \operatorname{diag}(s_1,\ldots,s_p)\,,
	\qquad
	s_\alpha = \frac{\sqrt{dv_w}}{\|w_\alpha\|_2}\,,
}
and $w_\alpha$ is the $\alpha$-th column of $W$. By construction, $Q$ and $\hat Q$ are defined from the same matrices $H$ and $Z$. Hence
\al{
	Q - \hat Q
	= c_1(H-\HH)+c_2(\bar H-\bar{\HH})\,,
}
so that
\al{
	\EE \big[\|Q-\hat Q\|_{\rm F}^2\big]
	\leq (|c_1|+|c_2|)^2 \;
	\EE \big[\|H-\HH\|_{\rm F}^2 \big] \,.
	\label{eq:prop2_Q_Qhat_reduction}
}
Conditionally on $W$, each row of $X$ is a proper Gaussian vector with covariance $v_x \1_d$. Therefore,
\al{
	\EE_X\big[\|H-\HH\|_{\rm F}^2\mid W\big]
	& = \frac{nv_x}{d}\|W(\1_p-S)\|_{\rm F}^2
	\\
	& = \frac{nv_x}{d}\sum_{\alpha=1}^p
	\big(\|w_\alpha\|_2-\sqrt{dv_w}\big)^2\,.
}
Since the columns of $W$ are identically distributed,
\al{
	\frac{1}{np}\EE \big[\|H-\HH\|_{\rm F}^2 \big]
	& = \frac{v_x}{d}\EE\big[
		\big( \|w_1\|_2-\sqrt{dv_w} \big)^2
	\big]
	\\
	& = v_h\,\EE\bigg[
	\mkakko{\frac{\|w_1\|_2}{\sqrt{dv_w}}-1}^2
	\bigg] \longrightarrow 0\,,
	\label{eq:prop2_H_Hhat_frobenius}
}
where we used the $L^2$ convergence of $\|w_1\|_2/\sqrt{dv_w}$ to $1$. Moreover,
\al{
	\frac{1}{np}\EE\big[\|H\|_{\rm F}^2\big] 
	= \frac{1}{np}\EE\big[\|\HH\|_{\rm F}^2\big] 
	= v_h\,, 
	\qquad \frac{1}{np}\EE\big[\|Z\|_{\rm F}^2\big] =v_z\,.
	\label{eq:23pDS}
}
Since
\al{
	\|Q\|_{\rm F} \leq (|c_1|+|c_2|)\|H\|_{\rm F}+\|Z\|_{\rm F}\,,
}
and the same bound holds for $\hat Q$ with $H$ replaced by $\HH$, it follows that
\al{
	\frac{1}{np}\EE\big[ \|Q\|_{\rm F}^2
	+ \|\hat Q\|_{\rm F}^2 \big] = O(1)\,.
	\label{eq:Q0Qh}
}
Together with \eqref{eq:prop2_Q_Qhat_reduction}--\eqref{eq:prop2_H_Hhat_frobenius} and \eqref{eq:Q0Qh}, Lemma~\ref{lem:stieltjes_comparison} proves \eqref{eq:closure_Q_Qhat}.

\paragraph{Step 2: \boldmath $\hat Q\longleftrightarrow\hat Q_M$.}
Throughout this step, all unlabeled limits are taken as $M\to\infty$. The Hermite truncations satisfy
\al{
	f_M \longrightarrow f \qquad
	\text{in~~}L^2\big(\CN(0,v_h,0)\big)\,.
}
The formulas in Section~\ref{sc:widely_linear} then imply
\al{
	c_{1,M}\longrightarrow c_1\,, \qquad
	c_{2,M}\longrightarrow c_2\,, \qquad
	v_{z,M}\longrightarrow v_z\,, \qquad
	\tau_{z,M}\longrightarrow\tau_z\,.
	\label{eq:prop2_surrogate_parameter_convergence}
}
Indeed, the first two limits follow from the continuity of the linear projections onto $h$ and $\bar h$ in Gaussian $L^2$, while the last two follow from the convergence of the variance and pseudo-variance of the corresponding residual.

We construct $Z_M$ and $Z$ from the same underlying standard Gaussian variables. Identify a complex Gaussian variable with its two-dimensional real representation, and let $\Sigma_M$ and $\Sigma$ be the real covariance matrices determined by $(v_{z,M},\tau_{z,M})$ and $(v_z,\tau_z)$, respectively. For \mbox{i.i.d.} standard real Gaussian vectors $\xx_{i\alpha}\in\RR^2$, define
\al{
	\bmat
		\Re (Z_M)_{i\alpha}\\
		\Im (Z_M)_{i\alpha}
	\emat
	= \Sigma_M^{1/2} \xx_{i\alpha}, \qquad
	\bmat
		\Re Z_{i\alpha}\\
		\Im Z_{i\alpha}
	\emat
	= \Sigma^{1/2} \xx_{i\alpha}\,.
}
By \eqref{eq:prop2_surrogate_parameter_convergence} and the continuity of the positive-semidefinite square root,
\al{
	\delta_M \equiv \EE
	\big[|(Z_M)_{i\alpha}-Z_{i\alpha}|^2\big] 
	\longrightarrow0\,,
}
where $\delta_M$ is independent of $i$ and $\alpha$. Consequently,
\al{
	\EE\big[\|Z_M-Z\|_{\rm F}^2\big]=np\delta_M\,.
	\label{eq:ll29Sdd}
}
Under this common construction,
\al{
	\hat Q_M-\hat Q = (c_{1,M}-c_1)\HH
	+(c_{2,M}-c_2)\bar{\HH} +(Z_M-Z)\,.
}
Applying $\|A+B+C\|_{\rm F}^2\leq 3\big(\|A\|_{\rm F}^2+\|B\|_{\rm F}^2+\|C\|_{\rm F}^2\big)$ and using \eqref{eq:23pDS} and \eqref{eq:ll29Sdd}, we obtain
\al{
	\frac{1}{np}\EE \big[\|\hat Q_M-\hat Q\|_{\rm F}^2\big]
	\leq{} & 3v_h|c_{1,M} - c_1|^2
	+ 3 v_h |c_{2,M}-c_2|^2 + 3 \delta_M
	\longrightarrow 0\,.
}
This upper bound is independent of $n,d,p$. Hence the left-hand side converges to zero uniformly in the matrix dimensions as $M\to\infty$. Moreover, \eqref{eq:prop2_surrogate_parameter_convergence} and the same argument as in \eqref{eq:Q0Qh} show that
\al{
	\frac{1}{np}\EE\big[\|\hat Q_M\|_{\rm F}^2\big] 
	= O(1)
}
uniformly in $n,d,p$ for all sufficiently large $M$. The corresponding bound for $\hat Q$ follows from \eqref{eq:Q0Qh}. Applying Lemma~\ref{lem:stieltjes_comparison} and then taking the high-dimensional limit superior, we obtain
\al{
	\lim_{M\to\infty}\limsup_{n,d,p\to\infty}
	\EE\big[\big|m_{\hat Q}(z)-m_{\hat Q_M}(z)\big|\big] =0\,,
}
which proves \eqref{eq:closure_Q_QhatM}.

\paragraph{Step 3: \boldmath $f_M(\HH)\longleftrightarrow f(\HH)$.}
Conditionally on $W$, every entry of $\HH$ has the common marginal distribution $\CN(0,v_h,0)$. Hence, for $h\sim\CN(0,v_h,0)$,
\al{
	\frac{1}{np}\EE \big[\|f_M(\HH)-f(\HH)\|_{\rm F}^2 \big]
	= \EE \big[|f_M(h)-f(h)|^2 \big] \longrightarrow 0
}
as $M\to\infty$. Similarly,
\al{
	\frac{1}{np}\EE\big[\|f_M(\HH)\|_{\rm F}^2\big]
	= \EE \big[|f_M(h)|^2\big] \,,
	\qquad
	\frac{1}{np}\EE\big[\|f(\HH)\|_{\rm F}^2\big]
	=\EE\big[|f(h)|^2\big]\,.
	\label{eq:24wtfsd}
}
The first quantity is bounded uniformly for all sufficiently large $M$ because $f_M\to f$ in Gaussian $L^2$. Lemma~\ref{lem:stieltjes_comparison} therefore gives
\al{
	\lim_{M\to\infty}\limsup_{n,d,p\to\infty}
	\EE\big[\big| m_{f_M(\HH)}(z)-m_{f(\HH)}(z)
	\big|\big]=0\,,
}
proving \eqref{eq:closure_fM_f}.

\paragraph{Step 4: \boldmath $f(\HH)\longleftrightarrow f(H)$.}
Set
\al{
	r_\alpha\equiv
	\frac{\|w_\alpha\|_2}{\sqrt{dv_w}}
	= s_\alpha^{-1}\,.
}
Since $\HH=HS$, we have, conditionally on $W$,
\al{
	H_{i\alpha} = r_\alpha \HH_{i\alpha}\,,
	\qquad
	\HH_{i\alpha}\mid W \sim\CN(0,v_h,0)\,.
}
The increment condition \eqref{eq:activation_increment_bound} therefore gives
\al{
	|f(\HH_{i\alpha})-f(H_{i\alpha})|^2
	\leq{} & C_f^2
	(r_\alpha - 1)^2|\HH_{i\alpha}|^2
	\mkakko{ 1 + |\HH_{i\alpha}|^q
	+ r_\alpha^q|\HH_{i\alpha}|^q }^2.
	\label{eq:prop2_increment_application}
}
Using $(a+b+c)^2\leq3(a^2+b^2+c^2)$ in \eqref{eq:prop2_increment_application}, and recalling that $\HH_{i\alpha}\mid W\sim\CN(0,v_h,0)$, we obtain
\al{
	\EE_X\big[
	|f(\HH_{i\alpha})-f(H_{i\alpha})|^2
	\mid W\big]
	& \leq 3 C_f^2 (r_\alpha - 1)^2
	\left\{ v_h + \big(1 + r_\alpha^{2q}\big)
	\EE\kkakko{|h|^{2q+2}} \right\}
	\\
	& \leq 3 C_f^2 (r_\alpha - 1)^2
	\mkakko{v_h + \EE\kkakko{|h|^{2q+2}}}
	\big(1 + r_\alpha^{2q}\big)\,,
}
where $h\sim\CN(0,v_h,0)$. Thus, there exists a constant $C<\infty$, independent of $n,d,p,i,\alpha$, and $W$, such that
\al{
	\EE_X\big[
	|f(\HH_{i\alpha})-f(H_{i\alpha})|^2
	\mid W\big]
	\leq
	C\, (r_\alpha - 1)^2 \big(1 + r_\alpha^{2q}\big)\,.
	\label{eq:prop2_increment_conditional}
}
Let $R_d\equiv\|w_\alpha\|_2^2/(dv_w)$. Then $r_\alpha = \sqrt{R_d}$ and $|r_\alpha - 1|\leq |R_d-1|$. Since $R_d$ is the average of $d$ \mbox{i.i.d.}\ nonnegative random variables with mean $1$ and finite moments of every order, independence and centering give
\al{
	\EE\bigl[ (R_d - 1)^4\bigr] 
	&= \frac{1}{d^4} \left\{ d\,\EE\left[
		\left( \frac{|W_{1 \alpha}|^2}{v_w} - 1 \right)^4
	\right] + 3d(d-1) \EE\left[
		\left( \frac{|W_{1 \alpha}|^2}{v_w} - 1 \right)^2
	\right]^2 \right\}
	\\
	& = O\bigl(d^{-2}\bigr)\,.
}
Moreover, $\sup_d\EE\big[ R_d^{2q} \big] < \infty$. Hence, by the Cauchy--Schwarz inequality,
\al{
	\EE\bigl[ (r_\alpha - 1)^2 \bigl(1 + r_\alpha^{2q}\bigr) \bigr]
	& \leq 
	\EE\bigl[(r_\alpha - 1)^2\bigr]
	+ \left\{ \EE\bigl[(r_\alpha - 1)^4\bigr]
	\EE\bigl[r_\alpha^{4q}\bigr] \right\}^{1/2}
	\\
	& \leq \EE\bigl[ (R_d - 1)^2 \bigr]
	+ \left\{ \EE\bigl[ (R_d - 1)^4 \bigr]
	\EE\bigl[ R_d^{2q} \bigr] \right\}^{1/2}
	\longrightarrow 0 \,.
}
Summing \eqref{eq:prop2_increment_conditional} over $i$ and $\alpha$ and taking expectations therefore yields
\al{
	\frac{1}{np}\EE\big[
	\|f(\HH)-f(H)\|_{\rm F}^2
	\big]
	&\leq
	C\,\EE\big[
		(r_1 - 1)^2
		\bigl(1 + r_1^{2q}\bigr)
	\big]
	\longrightarrow 0\,.
	\label{eq:prop2_destandardization_frobenius}
}
Using $f(H)=f(\HH)-\{f(\HH)-f(H)\}$, we obtain
\al{
	\|f(H)\|_{\rm F}^2 \leq
	2\|f(\HH)\|_{\rm F}^2 +
	2\|f(\HH)-f(H)\|_{\rm F}^2\,.
}
Thus, by \eqref{eq:24wtfsd} and \eqref{eq:prop2_destandardization_frobenius},
\al{
	\frac{1}{np}\EE\big[
	\|f(\HH)\|_{\rm F}^2+\|f(H)\|_{\rm F}^2
	\big]
	& \leq \frac{3}{np}\EE\big[\|f(\HH)\|_{\rm F}^2\big]
	+ \frac{2}{np}\EE\big[
	\|f(\HH)-f(H)\|_{\rm F}^2
	\big]
	\\
	& = 3\,\EE\big[|f(h)|^2\big]+o(1)
	= O(1)\,.
}
Together with \eqref{eq:prop2_destandardization_frobenius}, this verifies the two conditions of Lemma~\ref{lem:stieltjes_comparison} and proves \eqref{eq:closure_Hhat_H}.
\end{proof}

\subsection{Auxiliary Results II}

\begin{proposition}[Spectral equivalence for finite-degree truncated activations]
	\label{prop:finite_degree_equivalence}
	Fix an integer $M\geq1$. Then the expected empirical spectral distribution of $\hat Q_M^*\hat Q_M/n$ converges weakly to a compactly supported probability measure $\nu_M$ on $\RR_{\geq0}$. Moreover, for every fixed integer $k\geq1$,
	\al{
		\frac{1}{p}\EE\bigg[ \Tr\bigg\{ \bigg(
		\frac{1}{n}f_M(\HH)^*f_M(\HH)
		\bigg)^k \bigg\} \bigg]
		- \frac{1}{p}\EE\bigg[\Tr\bigg\{ \bigg(
		\frac{1}{n}\hat Q_M^*\hat Q_M
		\bigg)^k \bigg\} \bigg]
		\longrightarrow 0\,.
		\label{eq:finite_degree_moment_equivalence}
	}
	Consequently, the expected empirical spectral distribution of $f_M(\HH)^*f_M(\HH)/n$ also converges weakly to $\nu_M$.
\end{proposition}

It is instructive to first verify the moment matching for $k=1$, which in fact \emph{holds exactly at finite dimensions.} Let $h\sim \CN(0,v_h,0)$. As noted before in \eqref{eq:HH_dist_exact}, conditionally on $\WW$, every entry $\HH_{i\alpha}$ has the same distribution as $h$ and it does not depend on $\WW$. Consequently,
\al{
	\frac{1}{p}\EE\bigg[\Tr
	\left( \frac{1}{n}f_M(\HH)^*f_M(\HH) \right)\bigg]
	= \EE\big[|f_M(h)|^2\big]\,.
}
The first-order Hermite component of $f_M$ is $c_{1,M}h+c_{2,M}\bar h$, and its orthogonal complement has variance $v_{z,M}$. Hence
\al{
	\EE\big[|f_M(h)|^2\big] 
	= v_h \big( |c_{1,M}|^2+|c_{2,M}|^2 \big)+v_{z,M}
	= \EE\big[\big|(\hat Q_M)_{i\alpha} \big|^2\big]\,,
}
where we used the independence of $Z_M$ and $\HH$. It follows that
\al{
	\frac{1}{p}\EE\bigg[
	\Tr\left( \frac{1}{n}f_M(\HH)^*f_M(\HH) \right)
	\bigg]
	= \frac{1}{p}\EE\bigg[
	\Tr\left( \frac{1}{n}\hat Q_M^*\hat Q_M \right)
	\bigg]
}
for every $n,d,p$. 

\begin{figure}[tb]
	\centering
	\includegraphics[width=.73\textwidth]{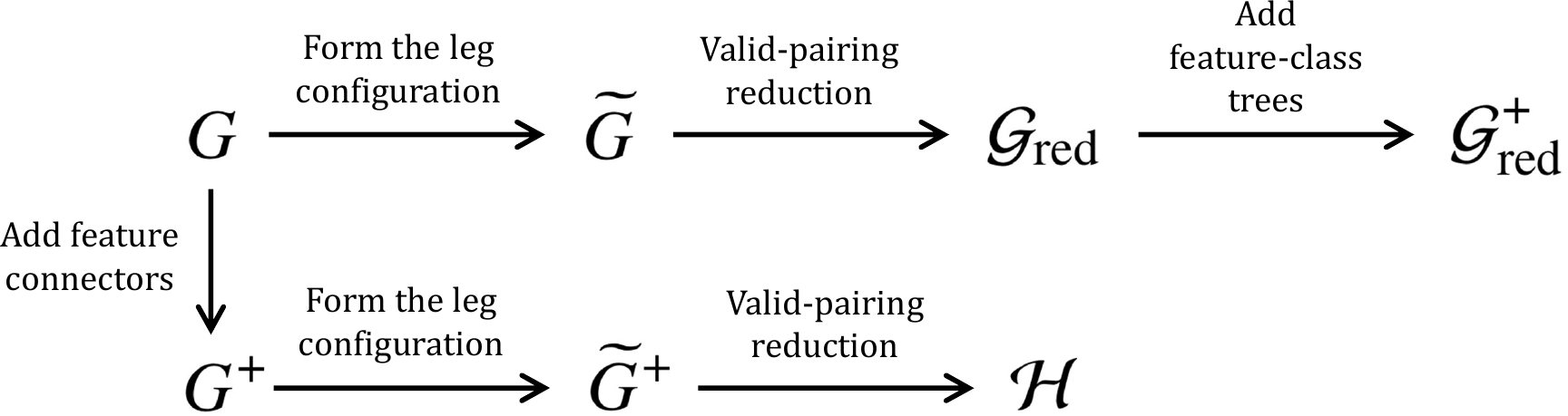}
	\caption{\label{fg:graph_relation}Relations among the graphs used in the proof of moment equivalence. With sample- and feature-index equality patterns fixed, the expectation over $X$ yields the Hermite contraction graph $G$. The subsequent normalized-column moment expansion, for a fixed data-index equality pattern, yields $\GG_{\rm red}$ and $\mathcal H$ through a Wick-pairing reduction. See the main text for details.}
\end{figure}

For higher moments, the joint dependence among the entries of $\HH$ becomes essential. We now introduce the graphical setup and notation needed to analyze this dependence. The relations between all graphical objects introduced below are summarized in Figure~\ref{fg:graph_relation} on page~\pageref{fg:graph_relation}. For any graph $K$, we write $C(K)$ for the number of its connected components. Throughout the proof, $X$ is proper. Recall from \eqref{eq:feature-matrices} that $\HH_{i\alpha}=(X\WW)_{i\alpha}/\sqrt d$. For fixed $W$, the family $\{\HH_{i\alpha}\}_{\alpha\in[p]}$ is jointly proper complex Gaussian for each fixed $i$, while different rows are independent. Its cross-moments are
\al{
	\EE_X\!\left[
	\HH_{i\alpha}\bar{\HH_{j\beta}}\,\middle|\,W
	\right]
	=
	\frac{v_x\delta_{ij}}{d} \sum_{\mu=1}^d
	\WW_{\mu\alpha}\bar{\WW_{\mu\beta}}
	\label{eq:Hhatcorr}
}
and $\EE_X[\HH_{i\alpha}\HH_{j\beta}\mid W]=0$. In particular, \eqref{eq:column-normalization} gives $\EE_X[|\HH_{i\alpha}|^2\mid W]=v_xv_w=v_h$.

\medskip
\noindent
$\rule{0.6em}{0.6em}$~\textbf{The Hermite contraction graph $G$ and the auxiliary graph $G^+$}%
\par\nopagebreak\smallskip\nopagebreak

\noindent
The first term in \eqref{eq:finite_degree_moment_equivalence} involves the product $\mkakko{f_M(\HH)^*f_M(\HH)}^k$. Substituting the Hermite expansion of $f_M$ defined in \eqref{eq:complex_hermite_truncation} gives a sum of products of $2k$ complex Hermite factors, which we represent graphically by $2k$ distinct \emph{factor nodes}. For brevity, let $\MM$ denote any such Hermite factor, either an ordinary factor $\He_{s,t}(\HH)$ or a conjugate factor $\bar{\He_{s,t}(\HH)}$. The product of the $\MM$s inherits a natural cyclic index structure from the single trace. For instance, when $k=4$, a typical index sum has the form
\al{
	\sum_{i,j,\ell,m = 1}^{n}
	\sum_{\alpha,\beta,\gamma,\delta = 1}^{p}
	\MM_{i\alpha}\MM_{i\beta}
	\MM_{j\beta}\MM_{j\gamma}
	\MM_{\ell\gamma}\MM_{\ell\delta}
	\MM_{m\delta}\MM_{m\alpha}\,.
	\label{eq:products_of_M}
}
Thus, every pair of cyclically neighboring factors shares either a sample index or a feature index. Additional index coincidences may also occur among nonneighboring factors because all sample and feature indices are summed over $[n]$ and $[p]$, respectively. If $i\neq j$, the $X$-variables appearing in $\MM_{i\alpha}$ and $\MM_{j\beta}$ are independent. The expectation over $X$ can therefore be taken separately within groups of factors carrying the same sample index. Similarly, distinct columns of $\WW$ are independent, so the expectation over $\WW$ can be taken separately within groups of factors carrying the same feature index.

Motivated by these observations, we decompose the sample- and feature-index sums according to equality patterns on the $2k$ factor nodes. A \emph{sample-index equality pattern} is a partition of the factor nodes into groups such that nodes in the same group carry the same sample-index value, whereas nodes in distinct groups carry distinct sample-index values. A \emph{feature-index equality pattern} is defined analogously. Both patterns must respect the index equalities imposed by the single-trace structure. Fix a sample-index equality pattern and a feature-index equality pattern, and define
\al{
	I_s = \text{the number of sample-index groups},
	\qquad
	I_f = \text{the number of feature-index groups}.
}
Because every group contains at least two factor nodes under the single-trace constraints, we have $I_s\leq k$ and $I_f\leq k$. The exact numbers of index assignments compatible with the fixed patterns are $(n)_{I_s}$ and $(p)_{I_f}$, where
\al{
	(m)_q = m(m-1)\cdots(m-q+1)\,.
}
Since $k$ and $M$ are fixed,
\al{
	(n)_{I_s} = n^{I_s} + O\bigl(n^{I_s-1}\bigr),
	\qquad
	(p)_{I_f} = p^{I_f} + O\bigl(p^{I_f-1}\bigr)\,.
}
Thus, the fixed equality patterns contribute the powers $n^{I_s}$ and $p^{I_f}$ to the leading power count.

As a concrete example, consider $k=4$ and label the eight factor-node positions by
\al{
	\MM_{i_1\alpha_1}\MM_{i_2\alpha_2}
	\MM_{i_3\alpha_3}\MM_{i_4\alpha_4}
	\MM_{i_5\alpha_5}\MM_{i_6\alpha_6}
	\MM_{i_7\alpha_7}\MM_{i_8\alpha_8}\,.
}
Here, $i_1,\ldots,i_8$ and $\alpha_1,\ldots,\alpha_8$ label index positions. Position labels remain distinct even when the corresponding index values are equal. The only admissible sample-index equality pattern with $I_s=4$ is
\al{
	\bigl\{\{i_1,i_2\},\{i_3,i_4\},\{i_5,i_6\},\{i_7,i_8\}\bigr\},
}
which is consistent with the cyclic constraints in \eqref{eq:products_of_M}. An admissible pattern with $I_s=3$ is
\al{
	\bigl\{\{i_1,i_2\},\{i_3,i_4\},\{i_5,i_6,i_7,i_8\}\bigr\},
}
while an admissible pattern with $I_s=2$ is
\al{
	\bigl\{\{i_1,i_2,i_5,i_6\},\{i_3,i_4,i_7,i_8\}\bigr\}.
}
Similarly, the only feature-index equality pattern with $I_f=4$ is
\al{
	\bigl\{\{\alpha_2,\alpha_3\},\{\alpha_4,\alpha_5\},
	\{\alpha_6,\alpha_7\},\{\alpha_8,\alpha_1\}\bigr\},
}
whereas an admissible pattern with $I_f=2$ is
\al{
	\bigl\{\{\alpha_4,\alpha_5\},
	\{\alpha_2,\alpha_3,\alpha_6,\alpha_7,\alpha_8,\alpha_1\}\bigr\}.
}

We now take the expectation over $X$. By Proposition~\ref{prop:he_M} in Appendix~\ref{app:hermite}, this expectation is a finite sum over directed contraction multigraphs compatible with the fixed sample-index equality pattern. A factor node representing $\He_{s,t}(\HH_{i\alpha})$ has $s$ outgoing edges and $t$ incoming edges. Self-loops, namely edges whose endpoints are the same factor node, are forbidden by Proposition~\ref{prop:he_M}. Each directed edge represents one Hermite contraction and contributes one cross-moment \eqref{eq:Hhatcorr}. Because this cross-moment contains $\delta_{ij}$, an edge can join only factor nodes in the same sample-index group. We denote such a graph on the original $2k$ factor nodes by $G$. Define
\al{
	E = \text{the number of edges in }G,
	\qquad C = C(G)\,.
}
A connected component of $G$ is called a \emph{dimer} if it contains exactly two factor nodes. All factor nodes in a connected component of $G$ belong to the same sample-index group. Conversely, one sample-index group may contain several connected components of $G$. Hence $I_s\leq C$. The quantities $I_s$ and $I_f$ are determined by the fixed equality patterns, whereas $E$ and $C$ are determined after a Hermite contraction graph $G$ has been fixed.

For example, let $k=4$ and fix the sample-index equality pattern
\al{
	\bigl\{\{i_1,i_2,i_5,i_6\},\{i_3,i_4,i_7,i_8\}\bigr\},
}
for which $I_s=2$. A compatible Hermite contraction may produce only two connected components, one within each sample-index group, or more than two components. For instance, within the first group it may produce two dimers whose factor-node positions are $\{i_1,i_5\}$ and $\{i_2,i_6\}$. Although all four positions carry the same sample-index value, they remain distinct factor-node positions in the graph.

For later use, we also define $G^+$. As observed above, the single-trace structure specifies $k$ pairs of neighboring factor nodes whose members carry the same feature index. We join the two nodes in each prescribed pair by an auxiliary line called a \emph{feature connector}. The graph obtained by adding the $k$ feature connectors to $G$ is denoted by $G^+$. For instance, when $k=4$, the feature connectors join the prescribed pairs $\{\alpha_2,\alpha_3\}$, $\{\alpha_4,\alpha_5\}$, $\{\alpha_6,\alpha_7\}$ and $\{\alpha_8,\alpha_1\}$. Feature connectors are auxiliary and are not counted among the $E$ Hermite-contraction edges. Thus, $G^+$ has $2k$ factor nodes, $E$ edges and $k$ feature connectors.

\medskip
\noindent
$\rule{0.6em}{0.6em}$~\textbf{The leg configuration $\tilde G$ and its auxiliary extension $\tilde G^+$}%
\par\nopagebreak\smallskip\nopagebreak

\begin{figure}[tb]
	\centering
	\includegraphics[width=.7\textwidth]{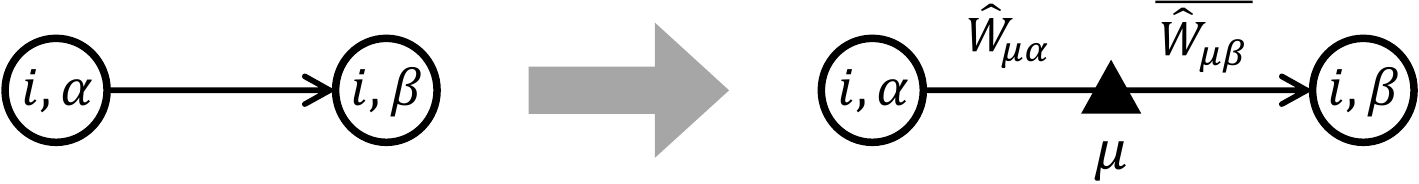}
	\caption{\label{fg:split_edge}An edge in $G$ (left) is split into two legs in $\tilde G$ (right) by inserting a data node, shown as a black triangle. The two factor nodes joined by the edge carry the same sample index ($i$ in this example); otherwise, the corresponding Hermite contraction vanishes. Each edge carries a factor $1/d$, originating from \eqref{eq:Hhatcorr}.}
\end{figure}

\noindent
After integrating out $X$, it remains to evaluate the resulting products of entries of $\WW$. For this purpose, we convert $G$ into a bipartite leg configuration, denoted by $\tilde G$. This construction is illustrated in Figure~\ref{fg:split_edge}. For each edge, we introduce a distinct \emph{data node} representing a data-index position $\mu\in[d]$, and replace the edge by two \emph{legs}, each joining one of the factor nodes to this data node. Thus, $\tilde G$ has $2k$ factor nodes, $E$ data nodes and $2E$ legs. Splitting edges into legs does not change the connectivity structure of the graph.

We further fix a \emph{data-index equality pattern}, namely a partition of the $E$ data nodes of $\tilde G$. The groups specify which data nodes carry equal data-index values: nodes in the same group carry the same value, and nodes in different groups carry distinct values. Define
\al{
	I_d = \text{the number of data-index groups}.
}
The exact number of admissible data-index assignments is $(d)_{I_d} = d^{I_d} + O\big(d^{I_d-1}\big)$.

If the feature connectors in $G^+$ are retained during the same splitting operation, the resulting graph is denoted by $\tilde G^+$. Feature connectors remain auxiliary objects: they are not split into legs and do not participate in the integration over $W$.

\medskip
\noindent
$\rule{0.6em}{0.6em}$~\textbf{Integration over $W$}%
\par\nopagebreak\smallskip\nopagebreak

\noindent
The $W$-dependent part of the contribution of $\tilde G$ is a product of entries of $\WW$ and their conjugates. For instance, consider a square graph ($E=4$) whose four factor nodes carry the same sample-index value $i$ and the same feature-index value $\alpha$, and label its four data nodes by $\mu_1,\mu_2,\mu_3,\mu_4$. Here, $\WW^{(\pm)}$ denotes either $\WW$ or $\bar{\WW}$, as appropriate. The contribution of the graph can then be written schematically as
\al{
	\prod_{j=1}^{4}
	\bigl(\WW_{\mu_j\alpha}^{(\pm)}\bigr)^2
	\ \sim\
	\bigl(\WW_{\square\alpha}^{(\pm)}\bigr)^8 \,.
}
Similarly, suppose that the feature-index values along the square are $\alpha,\alpha,\beta,\beta$, with $\alpha\neq\beta$. The contribution then takes the form
\al{
	\bigl(\WW_{\mu_1\alpha}^{(\pm)}\bigr)^2
	\WW_{\mu_2\alpha}^{(\pm)}
	\WW_{\mu_2\beta}^{(\pm)}
	\bigl(\WW_{\mu_3\beta}^{(\pm)}\bigr)^2
	\WW_{\mu_4\beta}^{(\pm)}
	\WW_{\mu_4\alpha}^{(\pm)}
	\ \sim\
	\bigl(\WW_{\square\alpha}^{(\pm)}\bigr)^4
	\bigl(\WW_{\square\beta}^{(\pm)}\bigr)^4 \,.
	\label{eq:Wmanyex}
}
Here, the symbols involving $\square$ suppress the data indices and record only how many factors belong to each feature column. The same structure extends to graphs of arbitrary fixed order. We are therefore led to evaluate fixed-degree moments of the columns of $\WW$.

A technical complication arises because the entries within a column of $\WW$ are dependent and are not jointly Gaussian owing to the column normalization. Nevertheless, Lemma~\ref{lem:normalized_column_wick} in Appendix~\ref{app:normalized_column_wick} shows that every fixed-degree joint moment of a normalized column of $\WW$ admits an expansion over the same pairings and the same data-index identifications as the Gaussian Wick expansion for the corresponding column of $W$. For each fixed pairing pattern, the scalar coefficient differs from the corresponding Gaussian Wick coefficient only by $O(d^{-1})$. Thus, at leading order as $d\to\infty$, the column normalization constraint can be ignored.

For the fixed feature-index equality pattern, the column-moment formula is applied separately to each feature-index group. For instance, \eqref{eq:Wmanyex} gives
\al{
	\EE\left[ \bigl(\WW_{\square\alpha}^{(\pm)}\bigr)^4 \right] 
	\EE\left[ \bigl(\WW_{\square\beta}^{(\pm)}\bigr)^4 \right].
}
Within each feature-index group, two legs can be paired only if their data-side endpoints belong to the same data-index group.

\medskip
\noindent
$\rule{0.6em}{0.6em}$~\textbf{The reduced graph $\GG_{\rm red}$}%
\par\nopagebreak\smallskip\nopagebreak

\noindent
Fix a Wick pairing of the $2E$ legs that is compatible with the fixed feature- and data-index equality patterns. The sample-index assignment of the $2k$ factor nodes plays no further role below, although its contribution has already been counted through $I_s$.

We construct a reduced bipartite multigraph $\GG_{\rm red}$ from $\tilde G$ by replacing each pair of legs with a single \emph{link}. For each pair whose legs emanate from distinct factor nodes, we merge their factor-side endpoints. The equivalence classes generated by all such factor-side merges are called \emph{reduced factor nodes}.

On the data side, recall that the $E$ data nodes are partitioned into $I_d$ data-index groups, with distinct groups carrying distinct data-index values. We represent each group by a single \emph{reduced data node}, giving exactly $I_d$ reduced data nodes in $\GG_{\rm red}$. A Wick pair contributes only if its two legs carry the same data-index value, so the data-side endpoints of every surviving pair are represented by the same reduced data node. Each link thus joins a reduced factor node to a reduced data node.

Define
\al{
	N_{\rm red}
	&= \text{the number of reduced factor nodes},
	\notag
	\\
	C_{\rm red} & = C(\GG_{\rm red})\,,
	\notag
	\\
	L_{\rm red}
	&= \text{the number of independent cycles in }\GG_{\rm red}
	   \text{ when the link directions are ignored}.
	\notag
}
The graph $\GG_{\rm red}$ has $N_{\rm red}+I_d$ nodes and $E$ links. Hence the Euler relation \autocite{Bondy1976book,Diestel2025book} reads
\al{
	N_{\rm red} + I_d - E = C_{\rm red} - L_{\rm red}\,.
	\label{eq:euler_relation_for_Gred}
}

\begin{figure}[tb]
	\centering
	\includegraphics[width=.37\textwidth]{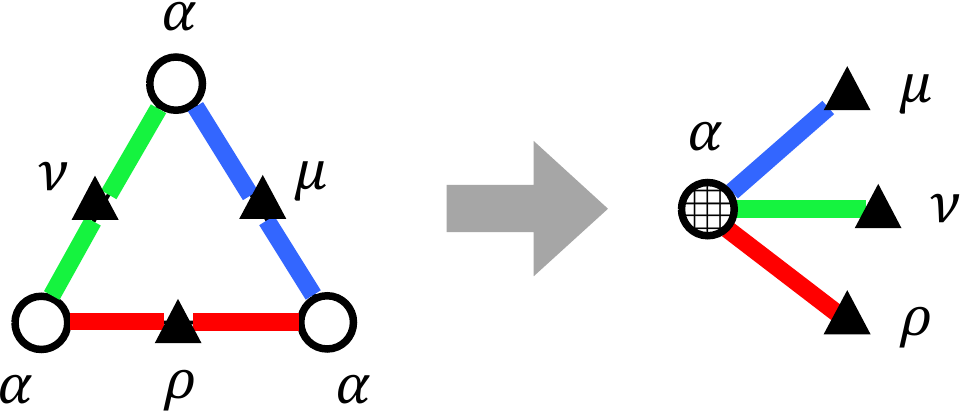}
	\qquad \qquad \quad 
	\includegraphics[width=.37\textwidth]{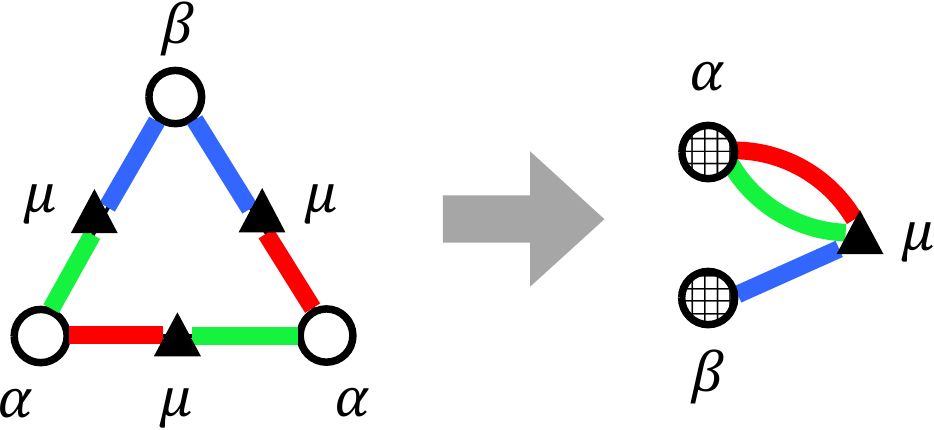}
	\put(-395, 60){\LARGE (a)}
	\put(-179, 60){\LARGE (b)}
	\vspace{1\baselineskip}
	\\
	\includegraphics[width=.42\textwidth]{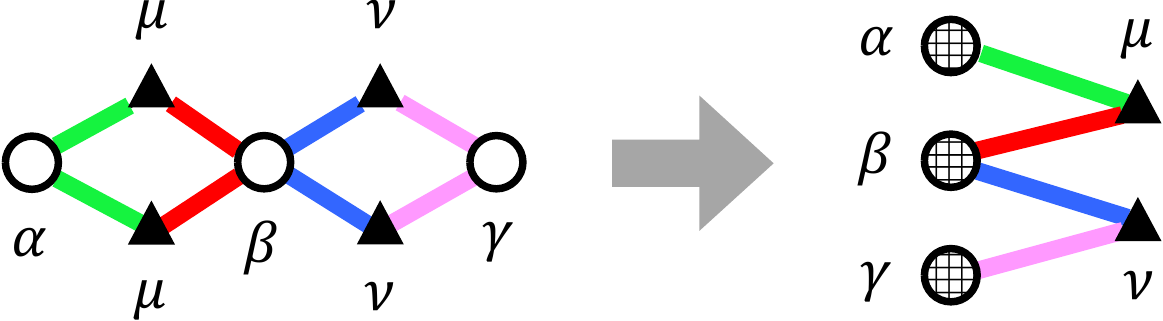}
	\qquad \quad
	\includegraphics[width=.42\textwidth]{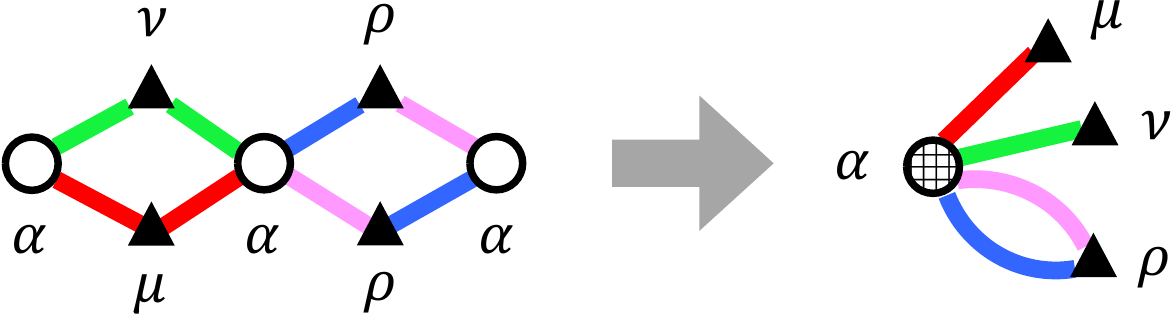}
	\put(-414, 50){\LARGE (c)}
	\put(-198, 50){\LARGE (d)}
	\vspace{.4\baselineskip}
	\caption{\label{fg:wick_graphs}Examples of the Wick-pairing reduction from $\tilde G$ to $\GG_{\rm red}$. In each panel, the graph on the left is $\tilde G$, with the two legs in each Wick pair shown in the same color, and the graph on the right is the resulting reduced graph $\GG_{\rm red}$. The hatched blobs represent reduced factor nodes, i.e., the equivalence classes generated by the factor-side mergers associated with the Wick pairing. Data nodes carrying the same data-index value are represented by a single reduced data node, while the links generated by distinct Wick pairs remain distinct and may therefore become parallel. Here, $\alpha,\beta,\gamma\in[p]$ denote distinct feature-index values, and $\mu,\nu,\rho\in[d]$ denote distinct data-index values. Sample-index labels are omitted because all factor nodes in each panel carry the same sample-index value.}
\end{figure}

Figure~\ref{fg:wick_graphs} illustrates how Wick pairings generate reduced factor nodes and links, while data nodes carrying the same data-index value are represented by a single reduced data node.

\medskip
\noindent
$\rule{0.6em}{0.6em}$~\textbf{The auxiliary reduced graphs $\GG_{\rm red}^+$ and $\mathcal H$}%
\par\nopagebreak\smallskip\nopagebreak

\noindent
We next define $\GG_{\rm red}^+$. For each feature-index group, add an arbitrary spanning tree of feature connectors on the reduced factor nodes belonging to that group. Applying this construction to all $I_f$ groups gives $\GG_{\rm red}^+$. Although the spanning trees are not unique, the connectivity of $\GG_{\rm red}^+$ is independent of their choice. Indeed, every such tree connects all reduced factor nodes in the corresponding feature-index group. Hence any path that uses a segment of one choice of spanning tree can be replaced by a path in any other spanning tree for the same group. It follows that any two reduced nodes are connected for one choice of the spanning trees if and only if they are connected for every other choice. Consequently, $C_{\rm red}^+$ is well defined.

The added trees contain $N_{\rm red}-I_f$ feature connectors in total. Define
\al{
	C_{\rm red}^+ = 
	C\bigl(\GG_{\rm red}^+\bigr) \,.
}
Since $\GG_{\rm red}^+$ is obtained by adding connectors to $\GG_{\rm red}$,
\al{
	C_{\rm red}^+ 	\leq C_{\rm red} \,.
}

The feature connectors in $\GG_{\rm red}^+$ need not be the literal images of those in $G^+$. What matters is their induced connectivity. The endpoints of the reduced image of every feature connector in $G^+$ belong to the same feature-index group and are therefore connected by the corresponding tree in $\GG_{\rm red}^+$.

Let $\mathcal H$ be obtained by applying the same pairing reduction to $\tilde G^+$ while retaining the images of its feature connectors. The factor-side mergers map the original factor nodes to reduced factor nodes, and the endpoints of every feature connector are transported by this map. Removing the feature connectors from $\mathcal H$ gives $\GG_{\rm red}$.

Figure~\ref{fg:graph_relation} summarizes the relations among these graphs and the corresponding operations. To illustrate the structures of $\tilde G^+$, $\mathcal H$ and $\GG_{\rm red}$ more explicitly, consider the following index assignment for $k=4$:
\begin{center}
	\begin{tabular}{lcccccccc}
		\toprule
		Factor node No.&1&2&3&4&5&6&7&8
		\\\midrule
		Sample index &$i$&$i$&$i$&$i$&$i$&$i$&$j$&$j$
		\\
		Feature index & $\beta$ & $\alpha$ & $\alpha$ 
		& $\beta$ & $\beta$ & $\alpha$ & $\alpha$ & $\beta$
		\\\bottomrule
	\end{tabular}
\end{center}
Here, $i,j\in[n]$ with $i\ne j$, and $\alpha,\beta\in[p]$ with $\alpha\ne\beta$. This index assignment is consistent with the single-trace structure.
\begin{figure}[tb]
	\centering
	\includegraphics[height=7.5em]{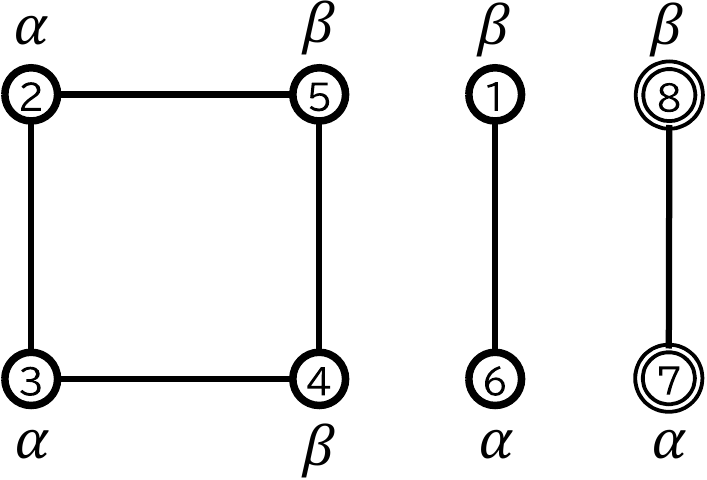}
	\qquad \qquad \qquad
	\includegraphics[height=7.5em]{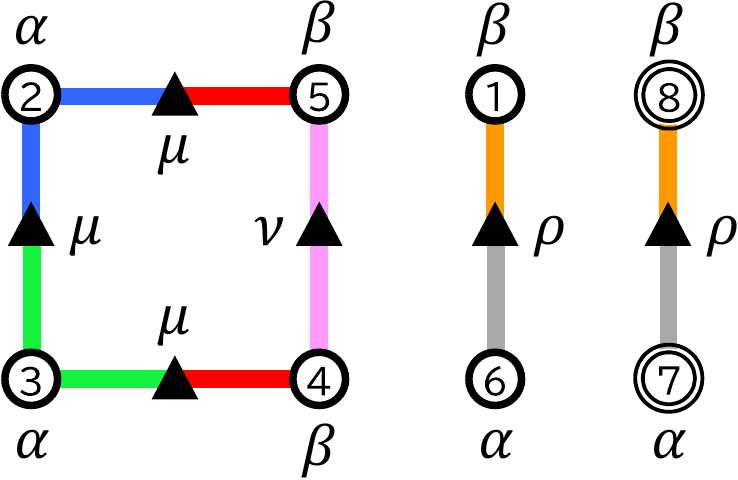}
	\put(-350, 71){\LARGE (a)}
	\put(-160, 71){\LARGE (b)}
	\vspace{1.2\baselineskip}
	\\
	\includegraphics[height=7.5em]{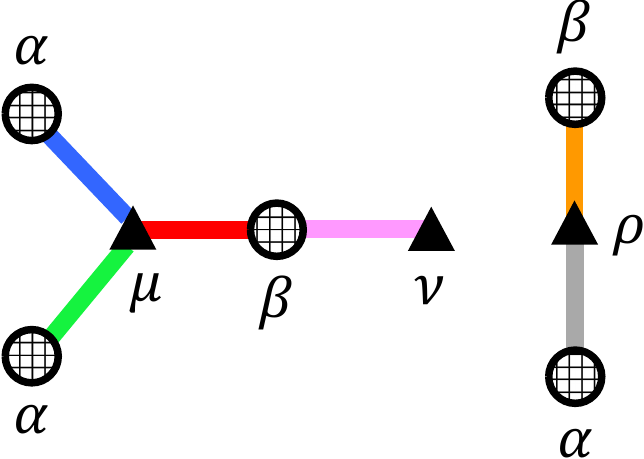}
	\qquad \qquad \quad
	\includegraphics[height=7.5em]{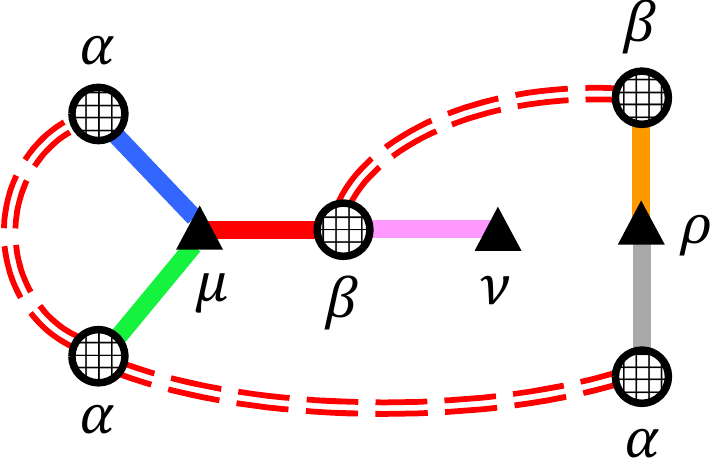}
	\put(-334, 69){\LARGE (c)}
	\put(-156, 69){\LARGE (d)}
	\vspace{1.2\baselineskip}
	\\
	\includegraphics[height=7.5em]{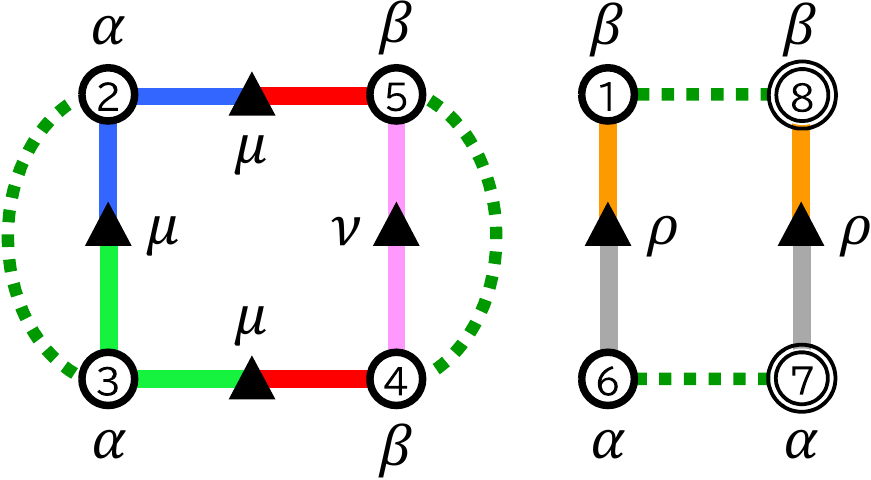}
	\qquad \qquad \quad 
	\includegraphics[height=7.5em]{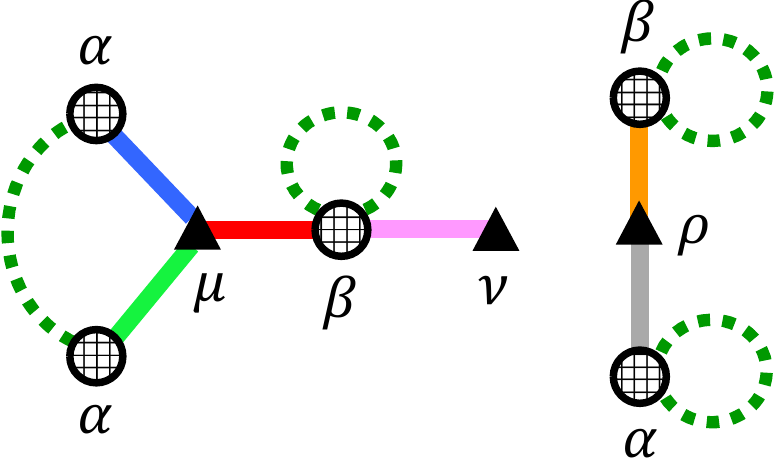}
	\put(-373, 69){\LARGE (e)}
	\put(-164, 69){\LARGE (f)}
	\caption{\label{fg:reduction_examples}Graphical constructions for the index assignment given in the preceding table: (a) $G$, (b) $\tilde G$, (c) $\GG_{\rm red}$, (d) $\GG_{\rm red}^+$, (e) $\tilde G^+$ and (f) $\mathcal H$. In panel~(a), the factor-node numbers correspond to the columns of the table. The directions of the edges are omitted for brevity. In panels~(a), (b) and~(e), the double circles mark factor nodes~$7$ and~$8$, which carry the sample-index value $j$; the remaining factor nodes carry $i$. In panels~(b) and~(e), the two legs in each Wick pair are drawn in the same color. The symbols $\mu,\nu,\rho\in[d]$ denote mutually distinct data-index values. The hatched blobs in panels~(c), (d) and~(f) represent reduced factor nodes. Green dotted lines denote the feature connectors prescribed by the single-trace structure. The red double-dashed lines in panel~(d) denote the feature-group spanning trees added to form $\GG_{\rm red}^+$.}
\end{figure}
Figure~\ref{fg:reduction_examples} illustrates the corresponding graphical constructions for one Hermite contraction graph, one data-index equality pattern and one Wick pairing. Panel~(a) shows the Hermite contraction graph $G$, with
\al{
	I_s=I_f=2\,, \qquad E=6\,, 
	\qquad C=3 \,.
}
Panel~(b) shows the bipartite graph $\tilde G$, obtained by splitting the six edges into legs and introducing one data node for each edge. The data-index equality pattern in this example has $I_d=3$. Panel~(c) shows the reduced graph $\GG_{\rm red}$ obtained from $\tilde G$ through the Wick-pairing reduction. Here,
\al{
	N_{\rm red}=5\,, \qquad C_{\rm red}=2\,,
	\qquad L_{\rm red}=0\,,
}
which indeed satisfies \eqref{eq:euler_relation_for_Gred}. Panel~(d) shows $\GG_{\rm red}^+$, obtained by adding a spanning tree to the reduced factor nodes in each feature-index group. These trees join the two connected components of $\GG_{\rm red}$, giving $C_{\rm red}^+ = 1 < C_{\rm red}$. 

By contrast, panel~(e) shows $\tilde G^+$, obtained by adding the four feature connectors prescribed by the single-trace structure to $\tilde G$ in panel~(b). Panel~(f) shows $\mathcal H$, obtained by applying the same Wick-pairing reduction to $\tilde G^+$ while retaining the images of these connectors. Three of the feature connectors become self-loops because their endpoints are merged by the reduction. The remaining connector does not join the two connected components, so that
\al{
	C(\mathcal H) = 2 > C_{\rm red}^+ = 1 \,.
}

\medskip
\noindent
$\rule{0.6em}{0.6em}$~\textbf{The graphical workflow}%
\par\nopagebreak\smallskip\nopagebreak

\noindent
The transformations above are generally one-to-many rather than one-to-one. The following schematic display records the sums required to reconstruct the original moment:
\al{
	\Tr\bigl[f_M(\HH)\dots\bigr]
	\begin{array}{l}
	\nearrow \cdots \\
	\to \Tr\bigl[\He_{s,t}(\HH)\dots\bigr]
	\\ \searrow \cdots
	\end{array}
	\!\! \to \EE_X
	\begin{array}{l}
	\nearrow \cdots \\
	\to ~G ~\to~ \tilde G
	\\ \searrow \cdots
	\end{array}
	\!\! \to \EE_W
	\begin{array}{l}
	\nearrow \cdots \\
	\to ~\GG_{\rm red}
	\\ \searrow \cdots
	\end{array}
}
Starting from a single trace of a product of the nonlinearly transformed matrices $f_M(\HH)$ and $f_M(\HH)^*$, we substitute the Hermite expansion of $f_M$ and obtain a sum of traces of products of complex Hermite factors. We first classify the sample- and feature-index assignments by their equality patterns. For each fixed pair of such patterns, averaging over $X$ gives a sum over compatible Hermite contraction graphs $G$. Each edge of $G$ is then represented by its two-leg configuration in $\tilde G$, after which we classify the data-index assignments by their equality patterns. Finally, the normalized-column moment formula gives a sum over valid columnwise pairings, each of which produces a reduced graph $\GG_{\rm red}$. At this stage, no random factor remains.

\medskip
\noindent
$\rule{0.6em}{0.6em}$~\textbf{Power counting}%
\par\nopagebreak\smallskip\nopagebreak

\noindent
Fix a Hermite contraction graph $G$, compatible sample-, feature- and data-index equality patterns, and a valid pairing. The exact number of index assignments is
\al{
	(n)_{I_s}(p)_{I_f}(d)_{I_d}\,.
}
Since $k$ and $M$ are fixed, this quantity has leading order $n^{I_s}p^{I_f}d^{I_d}$. Each Hermite-contraction edge contributes a factor $d^{-1}$. The Wick pairing coefficients are $O(1)$ by Lemma~\ref{lem:normalized_column_wick}. Hence the corresponding normalized moment term has the dimensional order
\al{
	\frac{n^{I_s}p^{I_f}d^{I_d-E}}{pn^k}
	=
	O\mkakko{p^{I_s+I_f+I_d-E-k-1}}
	\label{eq:leading_power_count}
}
in the proportional regime \eqref{eq:lmlm}.

\medskip
With this setup in place, we now prove Proposition~\ref{prop:finite_degree_equivalence}.

\begin{proof}[Proof of Proposition~\ref{prop:finite_degree_equivalence}]
We divide the proof into three parts. We first establish the limiting measure of the truncated Gaussian surrogate, then prove the moment equivalence by a combinatorial expansion, and finally transfer the limiting measure to the truncated nonlinear model.

\medskip
\noindent\textbf{$\bigcirc$\;Part~1:~Limiting measure of the truncated surrogate.}%
\par\nopagebreak\smallskip\nopagebreak
\noindent
Fix $M\geq1$ and introduce the corresponding surrogate before column normalization, 
\al{
	Q_M \equiv c_{1,M}H + c_{2,M}\bar H + Z_M\,.
}
For fixed $M$, the parameters $c_{1,M}$, $c_{2,M}$, $v_{z,M}$ and $\tau_{z,M}$ are independent of $n,d,p$, while $Z_M$ has \mbox{i.i.d.} complex Gaussian entries and is independent of $H$. Lemma~\ref{lem:surrogate_limit} therefore applies to $Q_M$ and shows that the expected empirical spectral distribution of $Q_M^*Q_M/n$ converges weakly to a compactly supported probability measure $\nu_M$ on $\RR_{\geq0}$.

We next transfer this convergence to the column-normalized surrogate $\hat Q_M$ defined in \eqref{eq:Qhat_M_definition}. By the definitions of $Q_M$ and $\hat Q_M$,
\al{
	Q_M-\hat Q_M = c_{1,M}(H-\HH)
	+ c_{2,M}(\bar H-\bar{\HH})\,.
}
Consequently,
\al{
	\frac{1}{np}\EE\big[ \|Q_M-\hat Q_M\|_{\rm F}^2 \big] 
	& \leq \mkakko{|c_{1,M}|+|c_{2,M}|}^2
	\frac{1}{np}\EE\big[ \|H-\HH\|_{\rm F}^2 \big]
	\longrightarrow 0 \,,
}
where the convergence follows from \eqref{eq:prop2_H_Hhat_frobenius}, established in Step~1 of Proposition~\ref{prop:approximation_steps}. Furthermore, for fixed $M$, the estimate leading to \eqref{eq:Q0Qh}, with the surrogate parameters replaced by their truncated counterparts, gives 
\al{
	\frac{1}{np}\EE\big[
	\|Q_M\|_{\rm F}^2+\|\hat Q_M\|_{\rm F}^2
	\big] =O(1)\,.
}
Lemma~\ref{lem:stieltjes_comparison} thus gives, for every fixed $z\in\CC_+$,
\al{
	\EE\big[ \big| m_{Q_M}(z)-m_{\hat Q_M}(z) \big| \big]
	\longrightarrow 0\,.
	\label{eq:QM_QhatM_stieltjes_difference}
}
Since the expected empirical spectral distribution of $Q_M^*Q_M/n$ converges weakly to $\nu_M$, \eqref{eq:QM_QhatM_stieltjes_difference} and the characterization of weak convergence in expectation in terms of Stieltjes transforms \autocite[Theorem~5.8(3a)]{FleermannKirsch2023} imply that the expected empirical spectral distribution of
$\hat Q_M^*\hat Q_M/n$ also converges weakly to $\nu_M$. At this stage, only weak convergence is asserted; convergence of the spectral moments will follow from the moment bounds established in Part~2.

\medskip
\noindent\textbf{$\bigcirc$\;Part~2:~Moment equivalence.}%
\par\nopagebreak\smallskip\nopagebreak
\noindent
We first apply the graphical expansion to the $k$-th expected moment of the nonlinear model,
\al{
	\frac{1}{p}\EE\bigg[ \Tr\bigg\{
	\bigg( \frac{1}{n}f_M(\HH)^*f_M(\HH) \bigg)^k
	\bigg\} \bigg] \,.
}
Our first objective is to identify which terms in this expansion can survive at leading order. We show that this requires every connected component of the Hermite contraction graph to be a dimer.

Fix a Hermite contraction graph $G$, compatible sample-, feature- and data-index equality patterns and a Wick pairing of the $2E$ legs. Since $k$ and $M$ are fixed, only finitely many such combinatorial choices occur, and all Hermite coefficients and normalized-column pairing coefficients are $O(1)$. It therefore suffices to bound one fixed contribution.

Every connected component of $G$ lies within a single sample-index group, whereas one sample-index group may contain several connected components. Suppose that the $a$-th sample-index group contains $c_a$ connected components of $G$. These components can be connected by adding $c_a-1$ auxiliary lines within that group. The total number of lines required for all sample-index groups is
\al{
	\sum_{a=1}^{I_s}(c_a-1) = C-I_s \,.
}

After these lines are added to $G^+$, the resulting graph is connected. Indeed, the added lines connect all factor nodes within each sample-index group, while the feature connectors join the prescribed pairs with common feature indices. Alternating between these sample- and feature-index connections traces the cyclic index structure of the original single trace and reaches all $2k$ factor nodes. Since adding one line can decrease the number of connected components by at most one,
\al{
	C(G^+) - 1 \leq C - I_s \,.
	\label{eq:Gplus_component_bound}
}

The passage from $G^+$ to $\mathcal H$ cannot increase the number of connected components. Splitting the edges into legs preserves connectivity, and the subsequent reduction maps every connected path to a connected walk in $\mathcal H$. Hence
\al{
	C(\mathcal H) \leq C(G^+) \,.
	\label{eq:H_Gplus_component_bound}
}
We next compare the connectivity of $\mathcal H$ and $\GG_{\rm red}^+$. Every feature connector retained in $\mathcal H$ joins two reduced factor nodes belonging to the same feature-index group. In $\GG_{\rm red}^+$, these two nodes are connected, possibly through several intermediate nodes, by a path in the corresponding feature-group tree. Thus, every path in $\mathcal H$ can be replaced by a path in $\GG_{\rm red}^+$, and hence
\al{
	C_{\rm red}^+ \leq C(\mathcal H) \,.
}
Together with \eqref{eq:Gplus_component_bound} and \eqref{eq:H_Gplus_component_bound}, this yields
\al{
	C_{\rm red}^+ \leq C(\mathcal H)
	\leq C(G^+) \leq C-I_s+1 \,.
	\label{eq:Gredplus_component_bound}
}

We next compare $\GG_{\rm red}$ and $\GG_{\rm red}^+$. By construction, the feature-group trees added to $\GG_{\rm red}$ contain $N_{\rm red}-I_f$ connectors in total. Each connector decreases the number of connected components by at most one, and hence
\al{
	C_{\rm red}-C_{\rm red}^+ \leq
	N_{\rm red}-I_f \,.
	\label{eq:inequality_on_spanning_trees}
}
Combined with \eqref{eq:Gredplus_component_bound}, this yields
\al{
	N_{\rm red}-C_{\rm red}
	\geq I_f-C_{\rm red}^+ \geq I_f-C+I_s-1 \,.
	\label{eq:factor_component_bound}
}
Combining \eqref{eq:factor_component_bound} with \eqref{eq:euler_relation_for_Gred}, we obtain
\al{
	E-I_d & =
	N_{\rm red}-C_{\rm red}+L_{\rm red}
	\\
	& \geq I_f-C+I_s-1+L_{\rm red} \,.
	\label{eq:master_graph_inequality}
}

The power count \eqref{eq:leading_power_count} bounds the fixed contribution by $\displaystyle O\mkakko{p^{I_s+I_f+I_d-E-k-1}}$. Using \eqref{eq:master_graph_inequality},
\al{
	I_s+I_f+I_d-E-k-1
	& = I_s+I_f-k-1-(E-I_d)
	\\
	&\leq C-k-L_{\rm red} \,.
}
Thus, every fixed contribution is
\al{
	O\mkakko{p^{C-k-L_{\rm red}}} .
	\label{eq:fixed_graph_power_bound}
}

Since the expansion \eqref{eq:complex_hermite_truncation} of $f_M$ has no constant term, isolated factor nodes are forbidden in $G$. Namely, each connected component of $G$ must contain at least two nodes. Since $G$ has $2k$ factor nodes, we conclude
\al{
	C \leq k \,.
}
Since $L_{\rm red}\geq0$, \eqref{eq:fixed_graph_power_bound} shows that every contribution is $O(1)$. A contribution can survive at leading order, i.e., have a nonzero high-dimensional limit, only if
\al{
	C=k \qquad \text{and} \qquad L_{\rm red}=0 \,.
	\label{eq:leading_graph_conditions}
}
Therefore, $G$ must be a collection of $k$ dimers. This condition is necessary but not sufficient.

A dimer containing exactly one edge is called \emph{thin}, while a dimer containing two or more parallel edges is called \emph{thick}. A Wick pair is called \emph{cross-dimer} if its two legs belong to different dimers of $G$. Figure~\ref{fg:inter_dimer_pairings} illustrates the three possible types of cross-dimer pairing: thin--thin, thin--thick and thick--thick. We next show that \emph{no leg of a thick dimer can participate in a cross-dimer Wick pairing} at leading order. This implies that the cross-dimer pairings in panels~(b) and~(c) of Figure~\ref{fg:inter_dimer_pairings} are excluded at leading order.

Let $D$ be a thick dimer with factor nodes $x$ and $y$, and suppose that one of its legs participates in a cross-dimer Wick pair. Choose an edge of $D$ containing this leg, and denote its data node by $\mu$. Thus, this edge is represented in $\tilde G$ by the path
\al{
	x\text{---}\mu\text{---}y \,.
}
Let $\pi$ denote the reduction map from the factor and data nodes of $\tilde G$ to the corresponding reduced nodes of $\GG_{\rm red}$.

\begin{figure}[tb]
	\centering
	\includegraphics[height=.09\textheight]{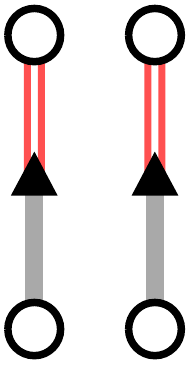}
	\qquad \qquad \quad
	\includegraphics[height=.09\textheight]{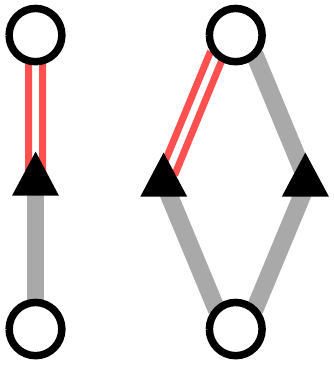}
	\qquad \qquad \quad
	\includegraphics[height=.09\textheight]{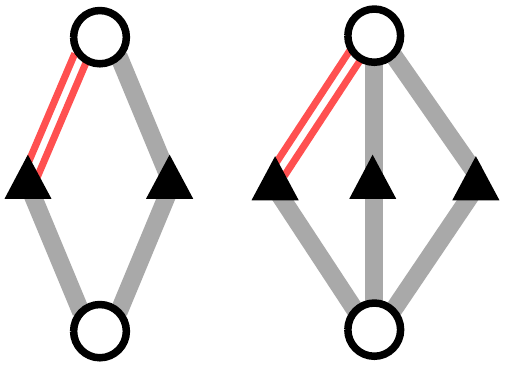}
	\put(-318,55){\LARGE (a)}
	\put(-229,55){\LARGE (b)}
	\put(-109,55){\LARGE (c)}
	\caption{\label{fg:inter_dimer_pairings}Schematic illustrations of cross-dimer Wick pairings. As before, open circles represent factor nodes and black triangles represent data nodes. In each panel, the two legs drawn with red double lines represent a cross-dimer Wick pair. Panels~(a)--(c) show pairings between two thin dimers, a thin and a thick dimer and two thick dimers, respectively.}
\end{figure}

First suppose that $\pi(x)=\pi(y)$. The two legs incident on $\mu$ belong to distinct Wick pairs, since one of them participates in a cross-dimer pairing. The two Wick pairs therefore produce distinct links with the same endpoints $\pi(x)=\pi(y)$ and $\pi(\mu)$. These parallel links form a cycle of length two, and hence
\al{
	L_{\rm red}\geq1 \,.
}

Now suppose that $\pi(x)\neq\pi(y)$. Since $D$ is thick, it contains another edge between $x$ and $y$. Denote its data node by $\nu$, so that this edge is represented by $x\text{---}\nu\text{---}y$.

First consider the case $\pi(\mu)=\pi(\nu)$. The Wick pairs containing $x\text{---}\mu$ and $x\text{---}\nu$ produce two links between $\pi(x)$ and the common reduced data node $\pi(\mu)=\pi(\nu)$. These links are distinct because $x\text{---}\mu$ already belongs to a cross-dimer Wick pair and hence cannot be paired with $x\text{---}\nu$. They therefore form a cycle of length two, so that
\al{
	L_{\rm red}\geq1 \,.
}

It remains to consider the case $\pi(\mu)\neq\pi(\nu)$. Since $\pi(x)\neq\pi(y)$, the two legs of each edge belong to distinct Wick pairs. The edges through $\mu$ and $\nu$ therefore produce the two paths
\al{
	\pi(x)\text{---}\pi(\mu)\text{---}\pi(y)\,, 
	\qquad \pi(x)\text{---}\pi(\nu)\text{---}\pi(y) \,.
}
Their intermediate nodes are distinct, so their union forms a cycle. Hence,
\al{
	L_{\rm red}\geq1 \,.
}

We have proved that if a thick dimer participates in a cross-dimer Wick pairing, then $L_{\rm red}\geq1$. Since every leading-order contribution must satisfy $L_{\rm red}=0$ by \eqref{eq:leading_graph_conditions}, no thick dimer can participate in a cross-dimer Wick pairing; every leg of a thick dimer must be paired with another leg in the same dimer.

This rule alone does not ensure that distinct thick dimers remain separate components in $\GG_{\rm red}$. Data nodes from distinct components of $\widetilde G$ are merged in the reduction whenever they belong to the same data-index group, and this can connect the corresponding components in $\GG_{\rm red}$. We show that such cross-dimer coincidences are subleading whenever a thick dimer is involved.

Consider an index equality pattern in which a thick dimer shares a data-index group with another dimer. We perform the {\it separation surgery} that preserves data-index equalities within each individual dimer and among thin dimers while completely removing every thin--thick or thick--thick coincidences. For example, if two data nodes in thick dimer $A$ and two in thick dimer $B$ belong to the same data-index group $\mu$, we replace
\al{
	\{\mu,\mu\}\subset\text{thick-$A$}, \qquad 
	\{\mu,\mu\}\subset\text{thick-$B$}
}
by
\al{
	\{\mu,\mu\}\subset\text{thick-$A$}, \qquad
	\{\mu',\mu'\}\subset\text{thick-$B$}
}
with $\mu$ and $\mu'$ denoting distinct data-index groups. As another example, suppose that five data nodes in the group $\mu$ lie in four dimers $(\text{thin-A},\text{thin-B},\text{thick-C},\text{thick-D})$ with multiplicities $(1,1,2,1)$. The surgery gives
\al{
	\{\mu\}\subset\text{thin-A},
	\qquad
	\{\mu\}\subset\text{thin-B},
	\qquad
	\{\mu',\mu'\}\subset\text{thick-C},
	\qquad
	\{\mu''\}\subset\text{thick-D},
}
where $\mu'$ and $\mu''$ are new data-index groups that occur in no other dimer.

The separation surgery does not alter the Wick pairing because no cross-dimer Wick pairing involves a thick dimer. If the surgery creates $r\geq1$ additional data-index groups, then $I_d$ increases by $r$, while $I_s$, $I_f$, $E$, and $k$ remain unchanged. Thus, by \eqref{eq:leading_power_count}, the exponent of $p$ for the original pattern is smaller by $r$ than that for the new pattern and is subleading; at leading order, no thick dimer shares a data-index group with another dimer. Hence the reduced image of a thick dimer forms a connected component of $\GG_{\rm red}$ by itself.

We next show that the two factor nodes of every thick dimer have the same feature index at leading order. We call a dimer \emph{homo} if its two factor nodes belong to the same feature-index group, and \emph{hetero} otherwise.

Suppose that at least one hetero thick dimer is present. Let $\widetilde\Omega$ be the collection of all hetero thick dimers before Wick reduction, and let $k'\geq1$ be their number. After Wick reduction, the reduced image of each of these dimers is a connected component of $\GG_{\rm red}$ by itself, as just shown. Since $L_{\rm red}=0$ at leading order, each such component is a tree. We denote the union of these $k'$ tree components by $\Omega$.

Each dimer in $\widetilde\Omega$ has two factor nodes with different feature indices. No Wick pair can join their legs, so the two factor nodes remain distinct after Wick reduction. Thus, $\Omega$ contains exactly $2k'$ reduced factor nodes, two from each hetero thick dimer.

We first count feature-index multiplicities before Wick reduction, among the factor nodes of $\widetilde\Omega$. Fix a feature index $\alpha$. The single-trace structure of the moment makes its total multiplicity among all factor nodes of $\widetilde G$ even. Its multiplicity among the factor nodes in homo dimers is also even, because each homo dimer contains either zero or two $\alpha$-nodes. Consequently, its total multiplicity in the hetero-dimer sector of $\widetilde G$ is even.

Suppose that $\alpha$ occurred an odd number of times among the factor nodes of the hetero thick dimers, i.e., in $\widetilde\Omega$. It would then occur an odd number of times among the hetero thin dimers as well. The homo thin dimers contribute an even number of $\alpha$-nodes, so the entire thin-dimer sector would contain an odd number of legs carrying $\alpha$. The thin-dimer legs carrying $\alpha$ must be paired among themselves: they cannot be paired with thick-dimer legs, and a Wick pair cannot join legs carrying different feature indices. A complete pairing requires an even number of such legs, contradicting the odd multiplicity obtained above. Thus, every feature index represented in $\widetilde\Omega$ must occur there an even number of times.

For the hetero thick dimers, Wick reduction does not merge factor nodes carrying different feature indices. Hence, by the parity argument above, every feature index represented in $\Omega$ occurs there an even number of times. Let $I_f'$ be the number of feature-index groups represented in $\Omega$. Each such group therefore contains at least two of the $2k'$ reduced factor nodes, so
\al{
	I_f'\leq k'\,.
	\label{eq:hetero_feature_group_bound}
}

We now simplify each tree component in $\Omega$. Only its two reduced factor nodes are relevant when feature connectors are added. Delete all pendant branches and replace the unique path between these two nodes by a single edge. After this simplification, $\Omega$ consists of $k'$ disjoint edges on $2k'$ nodes.

The next step is to add the feature connectors to construct $\GG_{\rm red}^+$. Its connected-component partition does not depend on the spanning trees chosen within the feature-index groups. For each group represented in $\Omega$, choose a tree that first connects its nodes in $\Omega$ and then, if necessary, extends to nodes outside $\Omega$. If $u_\alpha$ is the number of reduced factor nodes in $\Omega$ carrying $\alpha$, this adds
\al{
	\sum_{\alpha\,{\rm in}\,\Omega}(u_\alpha - 1)
	= 2 k' - I_f'
}
feature connectors within $\Omega$.

The resulting graph on $\Omega$ has $2k'$ nodes and
\al{
	E'\equiv k'+(2k'-I_f')=3k'-I_f'
}
edges. Let $C'$ and $L'$ be its numbers of connected components and independent cycles. Euler's relation gives
\al{
	C'-L' = 2k'-E' = I_f'-k'
	\leq0\,,
}
where the last inequality is \eqref{eq:hetero_feature_group_bound}. Since $\Omega$ is nonempty, $C'\geq1$, so
\al{
	L'\geq C'\geq1\,.
}
This means that the feature connectors have created a cycle within $\Omega$. Before they were added, the simplified hetero thick dimers were a disjoint union of edges and had no cycle. The cycle must contain an added feature connector. Removing that connector leaves its endpoints connected by the rest of the cycle, so the connector does not reduce the number of connected components. At least one of the $N_{\rm red}-I_f$ feature connectors added to $\GG_{\rm red}$ is redundant for connectivity. As a result, the inequality \eqref{eq:inequality_on_spanning_trees} is strict:
\al{
	C_{\rm red}-C_{\rm red}^+ \leq
	N_{\rm red}-I_f-1\,.
}
Then the bound \eqref{eq:fixed_graph_power_bound} loses at least one power of $p$, implying that every configuration containing a hetero thick dimer is subleading. We conclude that every thick dimer is homo at leading order.

Furthermore, for a homo thick dimer with $m$ Hermite edges, the leading contribution is obtained when its $m$ data nodes carry distinct indices, since this maximizes $I_d$. Legs belonging to different edges then cannot be paired, so the only possible Wick reduction pairs the two legs on each edge. The two factor nodes consequently merge into a single reduced factor node, and the reduced image of the dimer is a star-shaped tree with $m$ emanating links.

We now compare the leading contributions in the $f_M$ moment expansion with the Wick expansion of the surrogate $\hat Q_M$. Write
\al{
	\begin{gathered}
		f_M(h)=\ell_M(h)+r_M(h)\,,
		\\
		\ell_M(h)\equiv c_{1,M}h+c_{2,M}\bar h\,,
		\qquad
		r_M(h)\equiv
		\sum_{\substack{s,t\geq0\\2\leq s+t\leq M}}
		a_{s,t}\He_{s,t}(h,\bar h)\,.
	\end{gathered}
}
By Hermite orthogonality and the definitions of $v_{z,M}$ and $\tau_{z,M}$,
\al{
	v_{z,M}=\EE\big[|r_M(h)|^2\big]\,,
	\qquad
	\tau_{z,M}=\EE\big[r_M(h)^2\big]
	\qquad
	\text{for}~~h\sim\CN(0,v_h,0)\,.
	\label{eq:residual_second_moments}
}

We first consider the thin dimers. Their factor nodes are associated with either $\He_{1,0}(\HH,\bar{\HH})=\HH$ or $\He_{0,1}(\HH,\bar{\HH})=\bar{\HH}$ and therefore arise solely from $\ell_M(\HH)$ and its conjugate. Since $\ell_M(\HH)$ agrees with the widely linear part of $\hat Q_M$ defined in \eqref{eq:Qhat_M_definition}, the thin-dimer configurations are in a coefficient-preserving bijection with the corresponding Wick configurations of the widely linear part of $\hat Q_M$. This bijection preserves the factor positions, the sample- and feature-index groups, and the Wick pairings between thin dimers.

As shown above, every thick dimer contributing at leading order is homo, has no Wick pairing with another dimer, shares no data-index group with another dimer, and has distinct indices on its internal data nodes. Its reduced image is therefore a connected component of $\GG_{\rm red}$, graphically separated from the remaining dimers.

Fix a homo thick dimer with $m$ Hermite edges, sample index $i$ and feature index $\alpha$. To complete its internal sum, we restore all Hermite modes and all internal data-index assignments excluded by the preceding leading-order analysis, including those in which two or more edges carry the same data index. The preceding power bounds show that the configurations restored in this way contribute only $o(1)$ to the normalized global moment.

We also remove, at an $o(1)$ cost, the distinctness constraints between data indices belonging to different dimers. The added terms are precisely the cross-dimer data-index coincidences shown to be subleading by the separation-surgery argument. The data-index sums of the thick dimers can then be taken independently.

After this completion, the $m$ data-index sum on each edge gives (cf.~\eqref{eq:Hhatcorr} and Figure~\ref{fg:split_edge})
\al{
	\frac{v_x}{d} \sum_{\mu=1}^{d}
	\bar{\WW_{\mu\alpha}}\WW_{\mu\alpha}
	= v_xv_w = v_h
}
exactly. Thus the completed thick-dimer contribution no longer depends on $\WW$. Summing over its internal Hermite modes then gives exactly the corresponding second moment of $r_M$ in \eqref{eq:residual_second_moments}. In particular, when the two factor nodes arise from $r_M(\HH)$ and $\overline{r_M(\HH)}$, respectively, the completed contribution is
\al{
	\EE\big[r_M(h)\bar{r_M(h)}\big]	= v_{z,M} \,.
}
When both factor nodes arise from $r_M(\HH)$, it is
\al{
	\EE\big[ r_M(h)^2 \big] = \tau_{z,M} \,,
}
and the conjugate case gives $\bar{\tau_{z,M}}$. These are precisely the possible Wick contractions
\al{
	\begin{gathered}
		\EE\big[ (Z_M)_{i\alpha}\bar{(Z_M)_{j\beta}} \big]
		= v_{z,M}\delta_{ij}\delta_{\alpha\beta} \,,
		\qquad \EE\big[
		(Z_M)_{i\alpha}(Z_M)_{j\beta}
		\big] = \tau_{z,M}\delta_{ij}\delta_{\alpha\beta} \,,
		\\
		\EE\big[
		\bar{(Z_M)_{i\alpha}} \bar{(Z_M)_{j\beta}}
		\big] = \bar{\tau_{z,M}}\delta_{ij}\delta_{\alpha\beta} \,.
	\end{gathered}
}
The two factor nodes of every leading thick dimer already have the same sample and feature indices, so their index identifications agree with those imposed by the corresponding pair of entries of $Z_M$.

After these completions, the thin-dimer sector is in the coefficient-preserving correspondence described above with the widely linear part of the surrogate, while each completed thick dimer has exactly the coefficient and index identifications of one corresponding $Z_M$-pair. Conversely, every Wick configuration of $\hat Q_M=\ell_M(\HH)+Z_M$ determines a nonlinear dimer configuration of the same form: the $Z_M$-pairs specify the completed thick dimers, while the remaining widely linear factors and their Wick pairings specify the thin-dimer sector. This correspondence preserves the dimer partition, the sample- and feature-index equality patterns, and the Wick pairings within the thin-dimer sector. It therefore also matches configurations that are subleading because their equality patterns leave fewer free indices; the corresponding configurations have the same power of $p$ on both sides.

The completions above change the normalized moment only by $o(1)$. Since only finitely many Hermite degrees, equality patterns, contraction graphs, and Wick pairings occur for fixed $k$ and $M$, summing over all configurations preserves this error bound. Therefore,
\al{
	\frac{1}{p}\EE\bigg[\Tr \bigg\{\left(
	\frac{1}{n}f_M(\HH)^*f_M(\HH)
	\right)^k \bigg\}\bigg] 
	- \frac{1}{p}\EE\bigg[\Tr\bigg\{ \left( 
	\frac{1}{n}\hat Q_M^*\hat Q_M
	\right)^k \bigg\}\bigg] \longrightarrow 0\,. 
	\label{eq:truncated_moment_comparison}
}

The weak convergence established in Part~1 does not by itself imply convergence against the unbounded function $\lambda^k$. The fixed-order moment argument in the proof of Lemma~\ref{lem:surrogate_limit} applies equally to $\hat Q_M$, with Lemma~\ref{lem:normalized_column_wick} used for the normalized columns. At order $k+1$, it gives, along the proportional regime,
\al{
	\limsup_{n,d,p\to\infty} \frac{1}{p}
	\EE\bigg[\Tr\bigg\{ \left(
	\frac{1}{n}\hat Q_M^*\hat Q_M
	\right)^{k+1} \bigg\}\bigg] <\infty\,.
}
As in the proof of Lemma~\ref{lem:surrogate_limit}, this higher-moment bound gives uniform integrability of the $k$-th spectral moments. Together with the weak convergence proved in Part~1, it follows that
\al{
	\frac{1}{p}\EE\bigg[\Tr\bigg\{ \left(
	\frac{1}{n}\hat Q_M^*\hat Q_M
	\right)^k \bigg\}\bigg] 
	\longrightarrow \int_{\RR_{\geq0}}
	\lambda^k\,\rmd\nu_M(\lambda)\,.
}
Together with \eqref{eq:truncated_moment_comparison}, this proves 
\al{
	\frac{1}{p}\EE\bigg[\Tr\bigg\{ \left(
	\frac{1}{n}f_M(\HH)^*f_M(\HH)
	\right)^k \bigg\}\bigg] 
	\longrightarrow \int_{\RR_{\geq0}}
	\lambda^k\,\rmd\nu_M(\lambda)\,. 
	\label{eq:nonlinear_truncated_moment_limit}
}

\medskip
\noindent\textbf{$\bigcirc$\;Part~3:~Spectral convergence of the truncated nonlinear model.}%
\par\nopagebreak\smallskip\nopagebreak
\noindent
We now use the moment equivalence proved in Part~2 to show that the expected empirical spectral distribution of $f_M(\HH)^*f_M(\HH)/n$ converges weakly to $\nu_M$. Let $\bar\nu_{f_M}^{(n,d,p)}$ denote the expected empirical spectral distribution of
\al{
	\frac{1}{n}f_M(\HH)^*f_M(\HH)\,.
}
This is a deterministic probability measure on $\RR_{\geq0}$, and
\al{
	\int_{\RR_{\geq0}}\lambda^k\,
	\rmd\bar\nu_{f_M}^{(n,d,p)}(\lambda)
	=
	\frac{1}{p}\EE\bigg[\Tr\bigg\{\bigg(
	\frac{1}{n}f_M(\HH)^*f_M(\HH)
	\bigg)^k\bigg\}\bigg] \,.
}
Hence \eqref{eq:nonlinear_truncated_moment_limit} gives, for every fixed $k\geq1$,
\al{
	\int_{\RR_{\geq0}}\lambda^k\,
	\rmd\bar\nu_{f_M}^{(n,d,p)}(\lambda)
	\longrightarrow
	\int_{\RR_{\geq0}}\lambda^k\,\rmd\nu_M(\lambda)\,.
	\label{eq:all_truncated_moments_converge}
}
Moment convergence alone does not immediately imply weak convergence. So, we first establish tightness and then use the higher-moment bounds from Part~2 to pass moments to subsequential weak limits. Taking $k=1$ in \eqref{eq:all_truncated_moments_converge}, we obtain
\al{
	\limsup_{n,d,p\to\infty}
	\int_{\RR_{\geq0}}\lambda\,
	\rmd\bar\nu_{f_M}^{(n,d,p)}(\lambda)
	<\infty \,.
}
It follows that, for every $R>0$,
\al{
	\bar\nu_{f_M}^{(n,d,p)}
	\bigl([R,\infty)\bigr) 
	\leq \frac{1}{R} \int_{\RR_{\geq0}}\lambda\, 
	\rmd\bar\nu_{f_M}^{(n,d,p)}(\lambda)= O \big(R^{-1} \big)
}
uniformly for all sufficiently large $n,d,p$. The right-hand side can be made arbitrarily small by taking $R$ sufficiently large. Hence the expected empirical spectral distributions are tight.

Consider any subsequence along the proportional regime. By Prokhorov's theorem \autocite{Billingsley1999}, tightness gives a further subsequence, denoted by $(n_j,d_j,p_j)$, such that
\al{
	\bar\nu_{f_M}^{(n_j,d_j,p_j)}
	\longrightarrow \tilde\nu_M
}
weakly for some probability measure $\tilde\nu_M$ on $\RR_{\geq0}$. Fix $k\geq1$. Part~2, applied at order $k+1$, gives
\al{
	\sup_j
	\int_{\RR_{\geq0}}\lambda^{k+1}\,
	\rmd\bar\nu_{f_M}^{(n_j,d_j,p_j)}(\lambda)
	<\infty \,.
}
Consequently,
\al{
	\sup_j \int_{\lambda > R}\lambda^k\,
	\rmd\bar\nu_{f_M}^{(n_j,d_j,p_j)}(\lambda)
	\leq \frac{1}{R} \sup_j
	\int_{\RR_{\geq0}}\lambda^{k+1}\,
	\rmd\bar\nu_{f_M}^{(n_j,d_j,p_j)}(\lambda)
	\longrightarrow 0
}
as $R\to\infty$. The functions $\lambda^k$ are therefore uniformly integrable along the subsequence, so weak convergence yields, using \eqref{eq:all_truncated_moments_converge},
\al{
	\int_{\RR_{\geq0}}\lambda^k\,\rmd\tilde\nu_M(\lambda)
	=
	\lim_{j\to\infty} \int_{\RR_{\geq0}}\lambda^k\,
	\rmd\bar\nu_{f_M}^{(n_j,d_j,p_j)}(\lambda)
	=
	\int_{\RR_{\geq0}}\lambda^k\,\rmd\nu_M(\lambda) \,.
}
That is, $\tilde\nu_M$ and $\nu_M$ have the same moments of every order. Since $\nu_M$ is compactly supported, it is uniquely determined by its moments \autocite{Lin2017review}. Hence
\al{
	\tilde\nu_M = \nu_M \,.
}
We have thus shown that every subsequence along the proportional regime has a further subsequence converging weakly to $\nu_M$, which proves that
\al{
	\bar\nu_{f_M}^{(n,d,p)}
	\longrightarrow \nu_M
}
weakly along the proportional regime. This completes the proof of Proposition~\ref{prop:finite_degree_equivalence}.
\end{proof}

\subsection{Proof of Theorem~\ref{thm:complex_gep}}
\label{sc:spectral_gep_proof}

\begin{proof}
By the transposition argument in Section~\ref{sc:spectral_gep_overview}, it is enough to consider the case in which $X$ is proper. We work with the centered activation introduced in that section.

Let $m_\nu$ denote the Stieltjes transform of the compactly supported probability measure $\nu$ in Lemma~\ref{lem:surrogate_limit}. Since $\EE[m_Q(z)]$ is the Stieltjes transform of the expected empirical spectral distribution of $Q^*Q/n$, Lemma~\ref{lem:surrogate_limit} gives, for every fixed $z\in\CC_+$,
\al{
	\EE\big[m_Q(z)\big]
	\longrightarrow m_\nu(z)\,.
	\label{eq:surrogate_stieltjes_limit}
}

We now compare $Q$ and $f(H)$ along \eqref{eq:proof_chain}. Using $|\EE[Y]|\leq\EE[|Y|]$ at each approximation step, the triangle inequality gives, for every fixed $M$,
\al{
	\big| \EE\big[m_{f(H)}(z)\big] - \EE\big[m_Q(z)\big] \big|
	\leq{} &
	\EE\big[ |m_{f(H)}(z)-m_{f(\HH)}(z)| \big]
	\notag
	\\
	& + \EE\big[
	|m_{f(\HH)}(z)-m_{f_M(\HH)}(z)|
	\big]
	\notag
	\\
	& + \big| \EE\big[m_{f_M(\HH)}(z)\big]
	- \EE\big[m_{\hat Q_M}(z)\big] \big|
	\notag
	\\
	& + \EE\big[ |m_{\hat Q_M}(z)-m_{\hat Q}(z)| \big]
	\notag
	\\
	& + \EE\big[ |m_{\hat Q}(z)-m_Q(z)| \big]\,.
	\label{eq:main_proof_stieltjes_chain}
}

For fixed $M$, Proposition~\ref{prop:finite_degree_equivalence} shows that the expected empirical spectral distributions of $f_M(\HH)^*f_M(\HH)/n$ and $\hat Q_M^*\hat Q_M/n$ converge weakly to the same probability measure $\nu_M$. Since $\lambda\mapsto(\lambda-z)^{-1}$ is bounded and continuous on $\RR_{\geq0}$ for every fixed $z\in\CC_+$, it follows that
\al{
	\EE\big[m_{f_M(\HH)}(z)\big] -
	\EE\big[m_{\hat Q_M}(z)\big] 	\longrightarrow 0
}
as $n,d,p\to\infty$ with $M$ fixed.

The first and fifth terms on the right-hand side of \eqref{eq:main_proof_stieltjes_chain} also converge to zero by \eqref{eq:closure_Hhat_H} and \eqref{eq:closure_Q_Qhat} of Proposition~\ref{prop:approximation_steps}, respectively. Hence, for every fixed $M$,
\al{
	\limsup_{n,d,p\to\infty} \big|\EE\big[m_{f(H)}(z)\big]
	- \EE\big[m_Q(z)\big] \big|
	\leq{}&
	\limsup_{n,d,p\to\infty} \EE\big[
	|m_{f(\HH)}(z)-m_{f_M(\HH)}(z)| \big]
	\notag
	\\
	& + \limsup_{n,d,p\to\infty}\EE\big[
	|m_{\hat Q_M}(z)-m_{\hat Q}(z)| \big]\,.
	\label{eq:main_proof_stieltjes_limsup}
}

The right-hand side of \eqref{eq:main_proof_stieltjes_limsup} converges to zero as $M\to\infty$ by \eqref{eq:closure_Q_QhatM} and \eqref{eq:closure_fM_f} of Proposition~\ref{prop:approximation_steps}. Since the left-hand side is independent of $M$, we conclude that 
\al{
	\EE\big[m_{f(H)}(z)\big] - \EE\big[m_Q(z)\big]
	\longrightarrow 0
	\label{eq:main_proof_stieltjes_equivalence}
}
as $n,d,p\to\infty$. Combining \eqref{eq:surrogate_stieltjes_limit} and \eqref{eq:main_proof_stieltjes_equivalence} yields
\al{
	\EE\big[m_{f(H)}(z)\big] \longrightarrow m_\nu(z)
	\label{eq:nonlinear_stieltjes_limit}
}
for every fixed $z\in\CC_+$. By the characterization of weak convergence in expectation in terms of Stieltjes transforms \autocite[Theorem~5.8(3a)]{FleermannKirsch2023}, \eqref{eq:nonlinear_stieltjes_limit} implies that the expected empirical spectral distribution of $f(H)^*f(H)/n$ converges weakly to $\nu$. Together with Lemma~\ref{lem:surrogate_limit}, this proves \eqref{eq:expected_esd_equivalence}. This completes the proof.
\end{proof}

\section{Numerical Experiments}
\label{sc:numerics}

In this section, we numerically test the spectral equivalence established in Theorem~\ref{thm:complex_gep}. To examine the scope of this equivalence, we consider not only proper $X$ and activation functions $f$ satisfying \eqref{eq:activation_increment_bound}, but also improper $X$ and activation functions outside this regularity condition. Although Theorem~\ref{thm:complex_gep} is stated in terms of the eigenvalues of the feature covariance matrices, we instead plot the nonzero singular values of $f(H)/\sqrt{n}$, whose squares are exactly the nonzero eigenvalues of $f(H)^*f(H)/n$. Taking square roots spreads out the values near the origin and makes the spectral comparison clearer.

The plotting range is chosen to display the bulk spectrum. Isolated values outside this range are not shown, since they do not affect the bulk comparison addressed by Theorem~\ref{thm:complex_gep} and would obscure the displayed spectrum. The number of independent matrix realizations used to construct each histogram is denoted by $N_{\mathrm{MC}}$ and displayed in the title of each figure. Without loss of generality, we set $v_x=1$ throughout this section. The code used to perform all experiments and generate all figures in this section is publicly available on GitHub.\footnote{\url{https://github.com/TaTKSM/complex-gep}}

\subsection{Even Functions}

As shown in Lemma~\ref{lem:bothczero} of Section~\ref{sc:vavasd}, $c_1$ and $c_2$ vanish whenever $f$ is even, reducing the Gaussian surrogate to a matrix with independent Gaussian entries and yielding a scaled \MP law. We numerically test this prediction using three activation functions: $f(z)=|z|$ (Modulus), $f(z)=|\Re z|+i|\Im z|$ (Split Modulus), and $f(z)=\cos z$ (Cosine). The last choice lies outside the regularity assumptions of Theorem~\ref{thm:complex_gep}.

\begin{figure}[tb]
	\centering
	\mbox{
	\includegraphics[width=.33\textwidth]{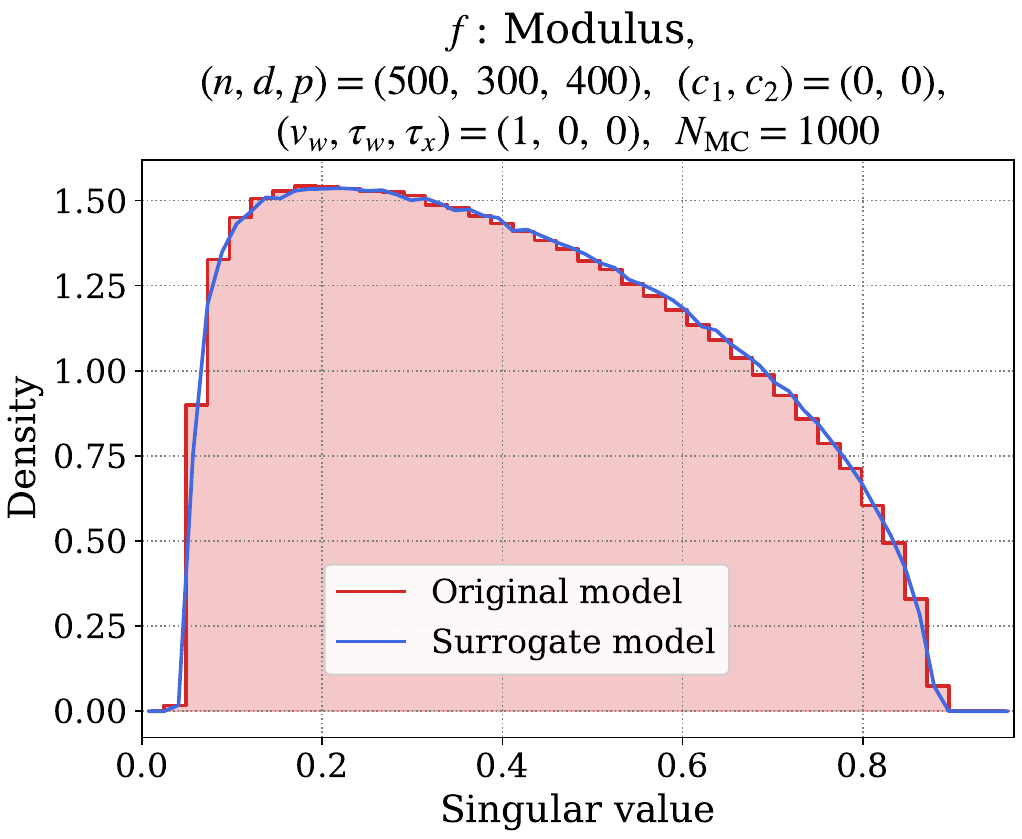}
	\includegraphics[width=.33\textwidth]{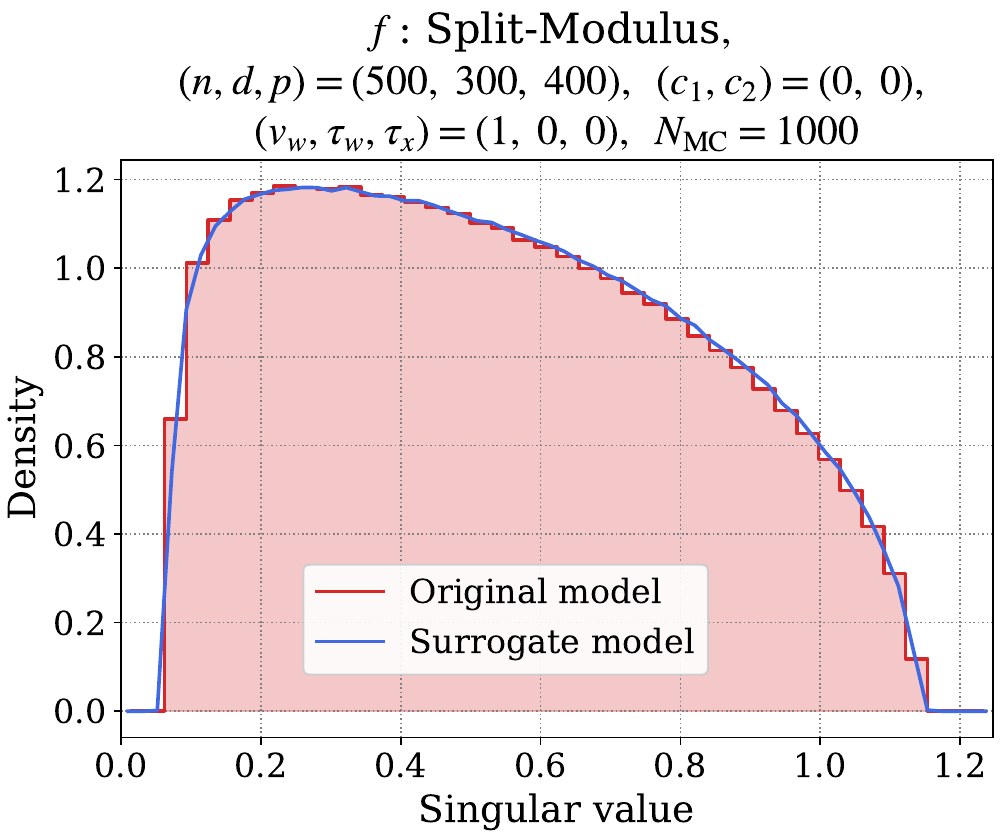}
	\includegraphics[width=.33\textwidth]{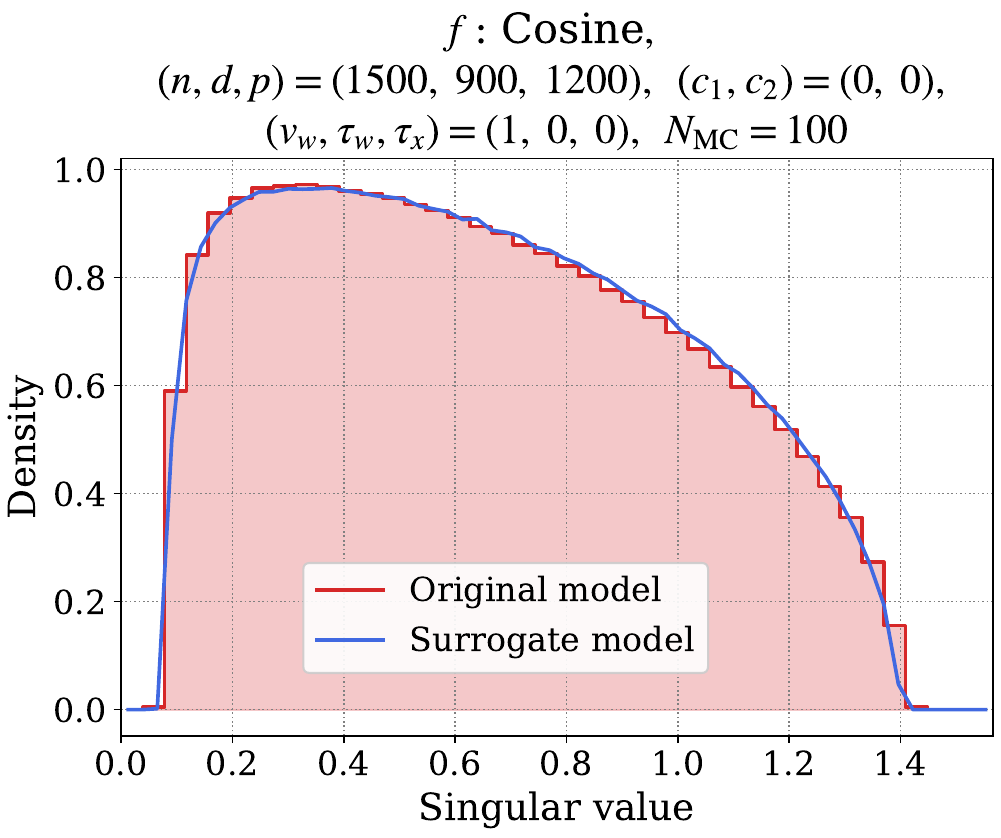}
	}
	\vspace{-6pt}
	\\
	\mbox{
	\includegraphics[width=.33\textwidth]{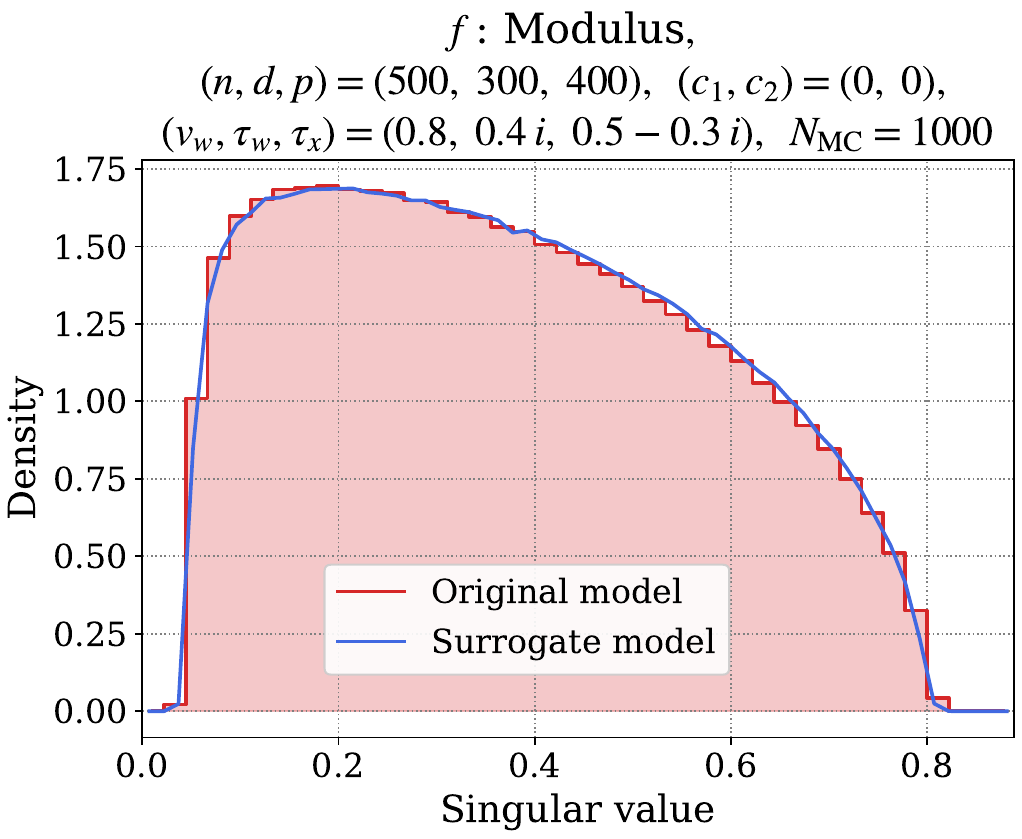}
	\includegraphics[width=.33\textwidth]{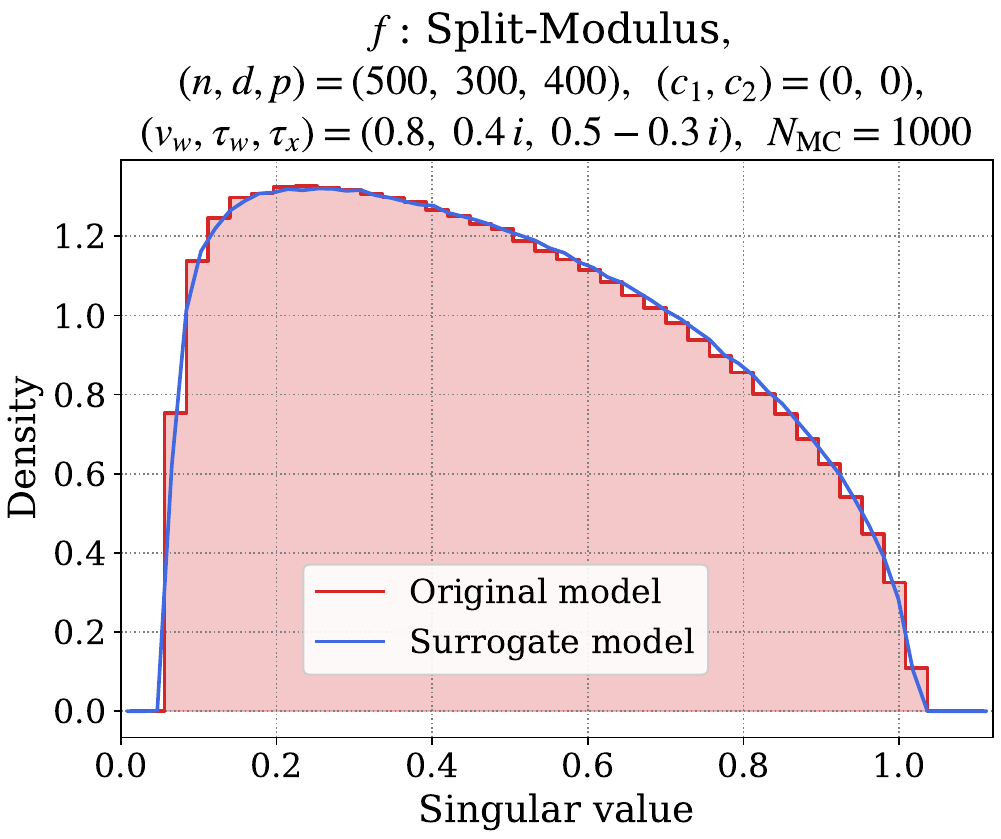}
	\includegraphics[width=.33\textwidth]{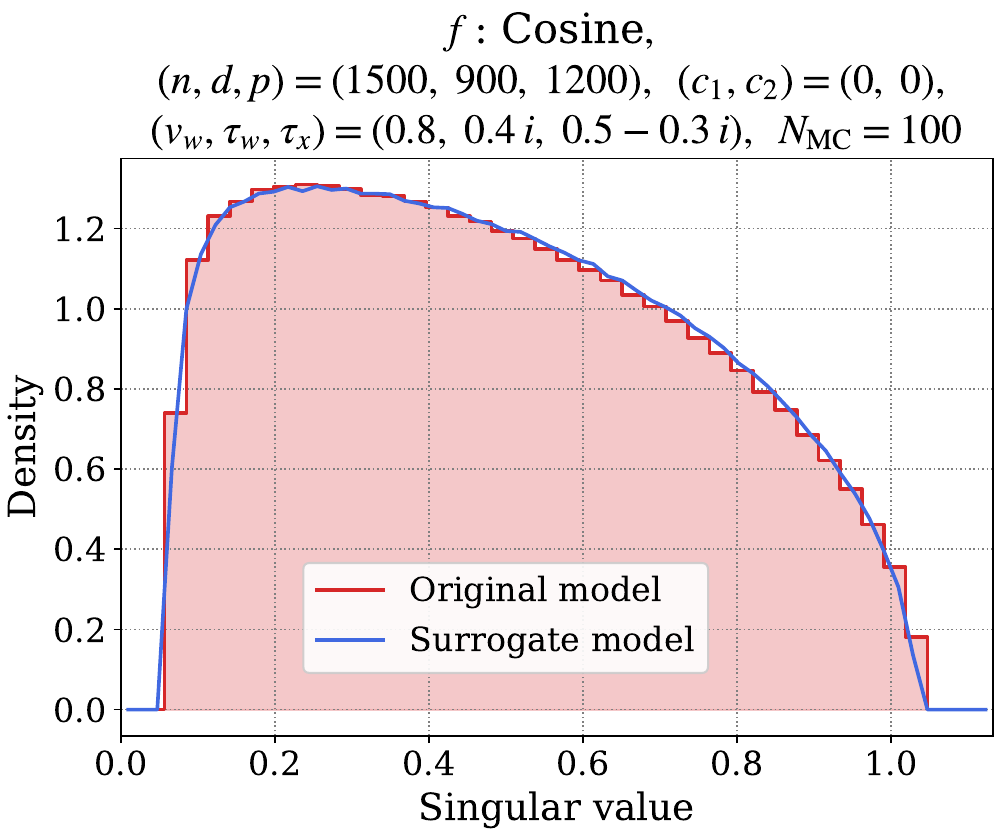}
	}
	\caption{\label{fg:even_f}Comparison of the empirical singular-value distributions of the original model (red histograms) and the Gaussian surrogate (blue curves) for $\tau_x=\tau_w=0$ (top row) and $\tau_x\ne0$ and $\tau_w\ne0$ (bottom row). The three columns correspond, from left to right, to Modulus, Split Modulus, and Cosine.}
\end{figure}

The results are shown in Figure~\ref{fg:even_f}. The agreement is excellent in both the proper and improper settings. Convergence is somewhat slower for the cosine activation, so we triple the matrix dimensions in the corresponding experiments to obtain a clearer finite-dimensional agreement.

\begin{figure}[tb]
	\centering
	\mbox{
	\includegraphics[width=.33\textwidth]{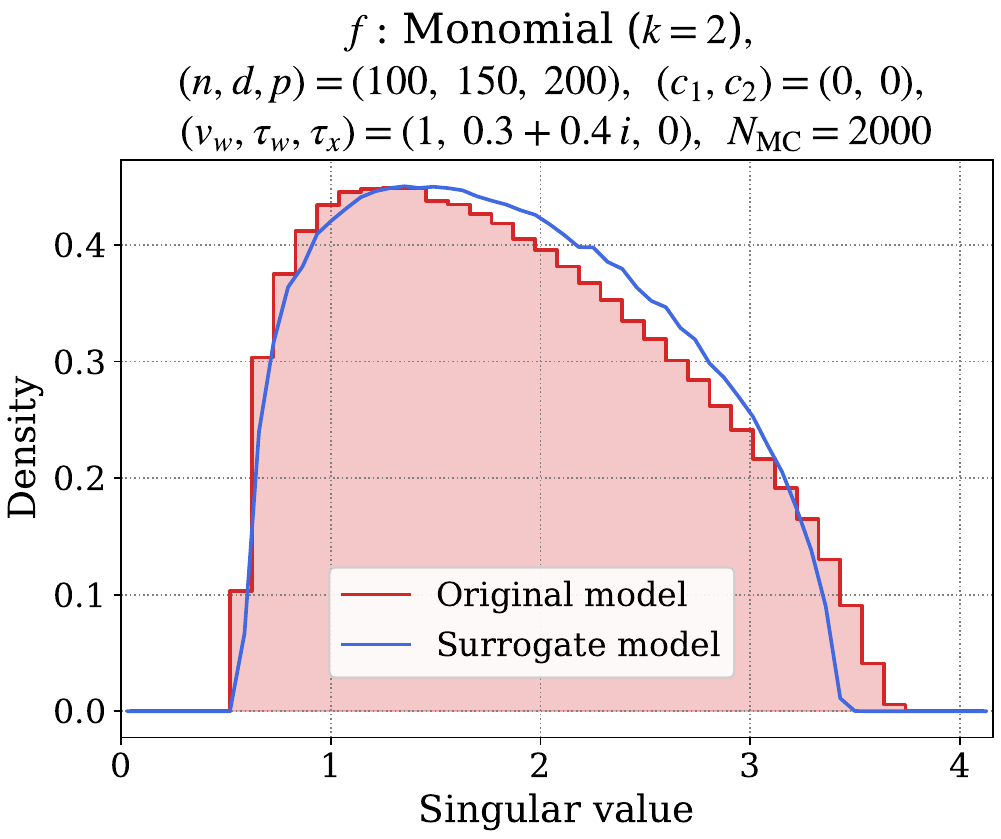}
	\includegraphics[width=.33\textwidth]{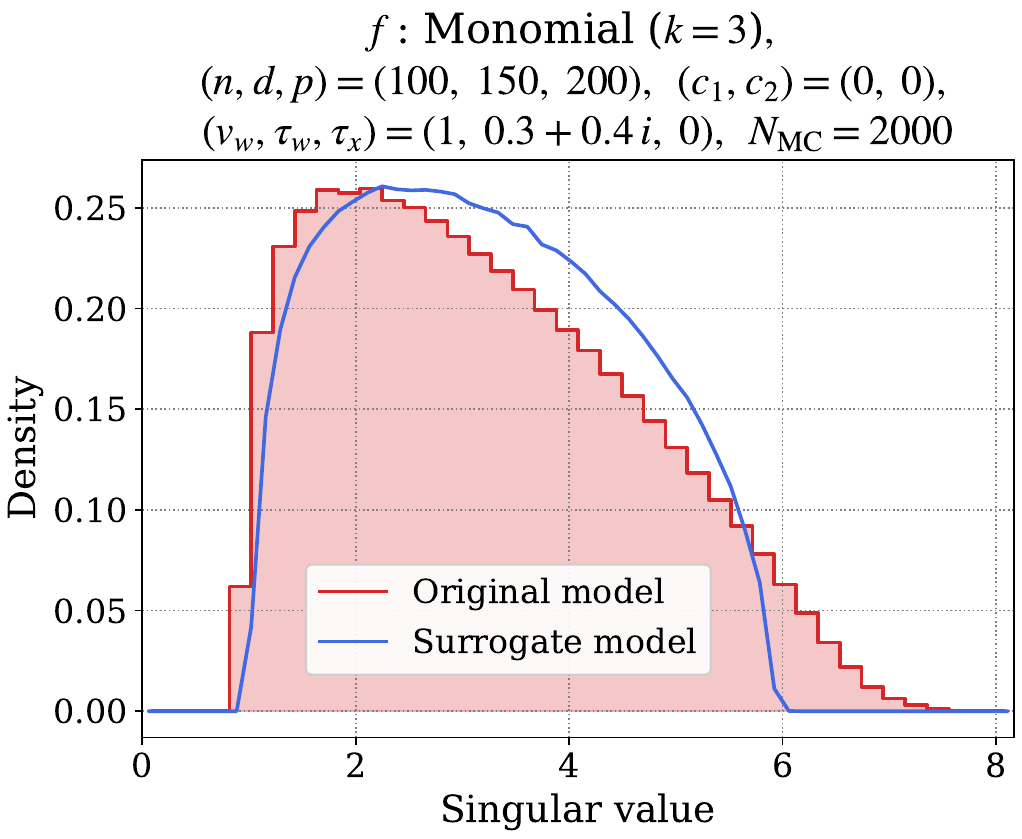}
	\includegraphics[width=.33\textwidth]{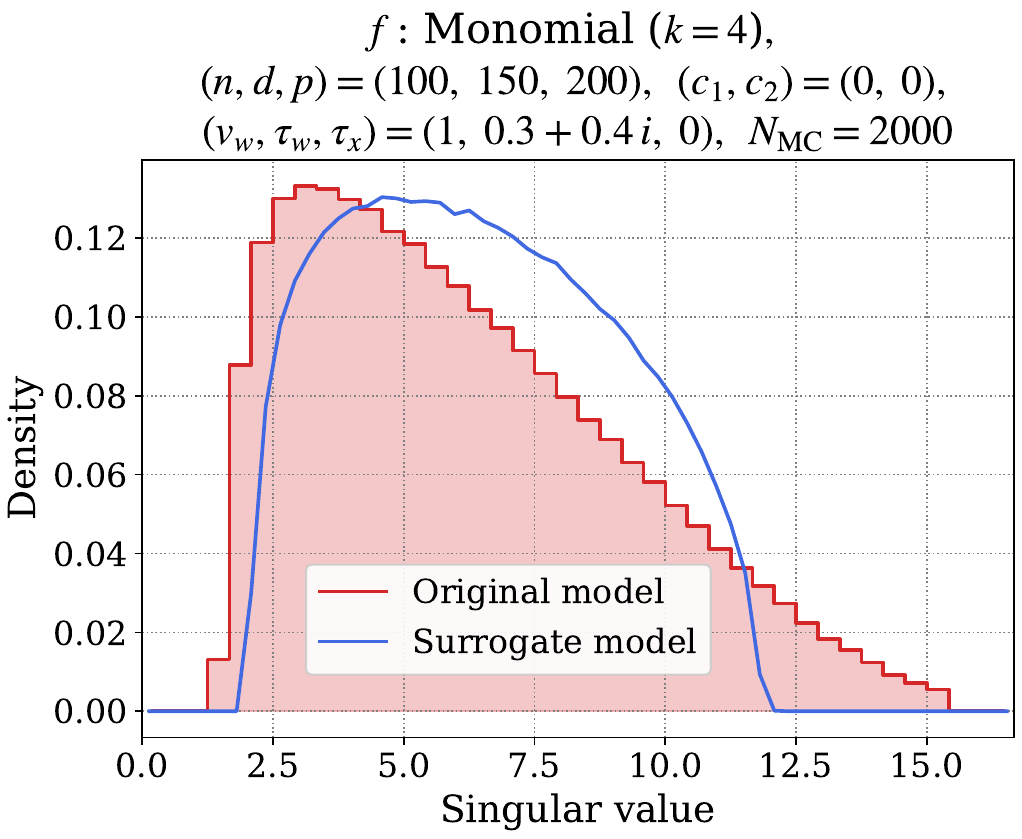}
	}
	\vspace{-4pt}
	\\
	\mbox{
	\includegraphics[width=.33\textwidth]{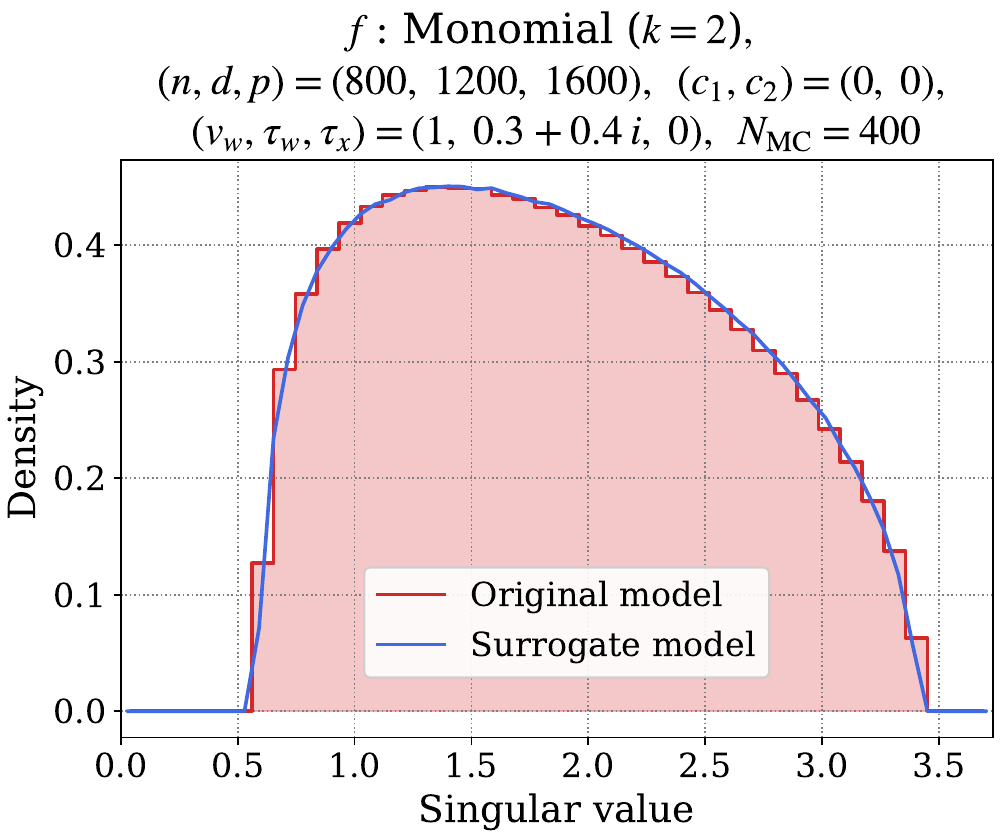}
	\includegraphics[width=.33\textwidth]{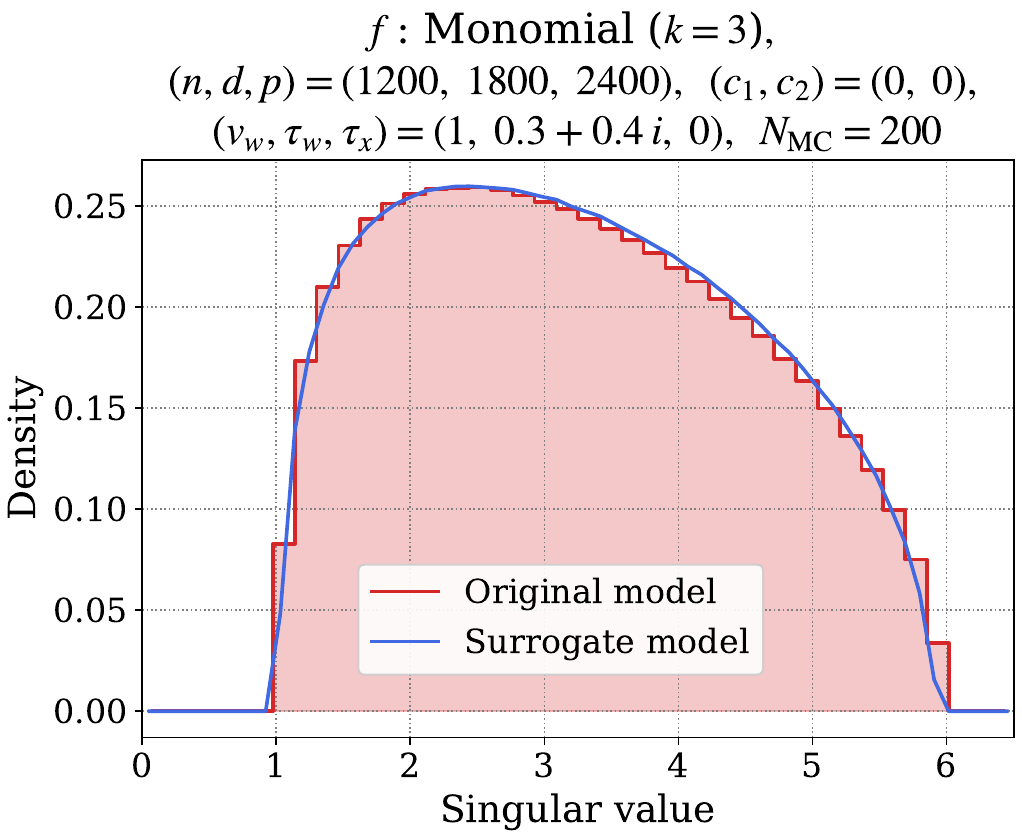}
	\includegraphics[width=.33\textwidth]{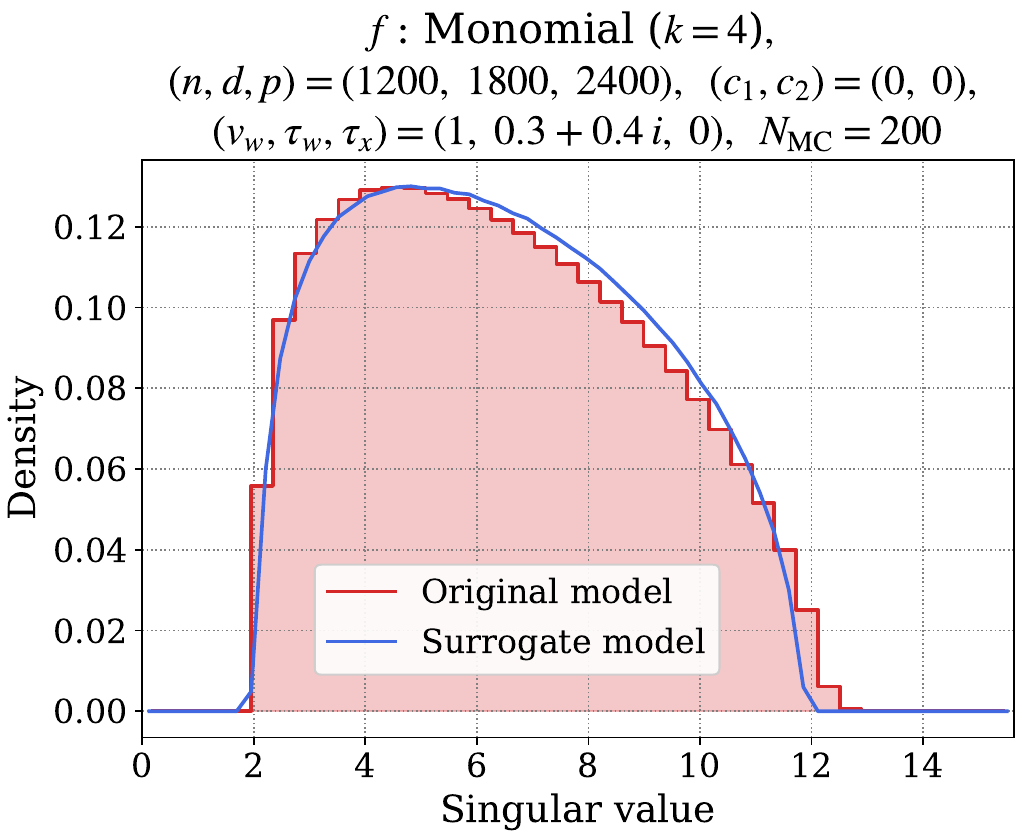}
	}
	\caption{\label{fg:monomials}Comparison for the monomial activations $z^k$, with $k=2,3,4$ from left to right. The matrix dimensions are increased from the top row to the bottom row while their ratios are preserved.}
\end{figure}

\subsection{Monomials}

We next consider the monomial activations $f(z)=z^k$ for $k=2,3,4$. For $k=2$ and $4$, evenness gives $c_1=c_2=0$, whereas for $k=3$, we have $c_1=3\tau_h$ and $c_2=0$ (recall Lemma~\ref{lem:c2zero}). Hence both coefficients vanish for all $k\in\{2,3,4\}$ when $\tau_h=\tau_x\tau_w=0$. Figure~\ref{fg:monomials} compares the original and surrogate models at two different matrix sizes. The top row shows that the finite-dimensional discrepancy becomes more pronounced as $k$ increases. The bottom row shows that this discrepancy decreases substantially when the dimensions are increased with their ratios fixed.

\begin{figure}[tb]
	\centering
	\includegraphics[width=.33\textwidth]{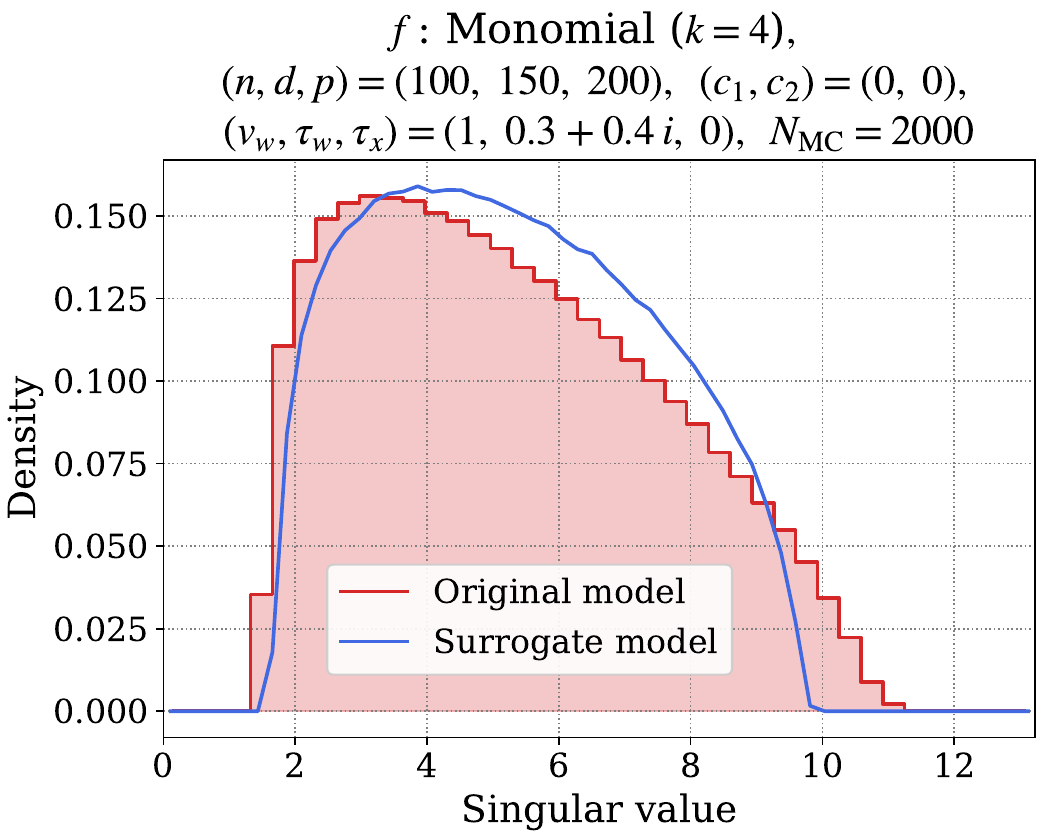}
	\qquad \quad
	\includegraphics[width=.33\textwidth]{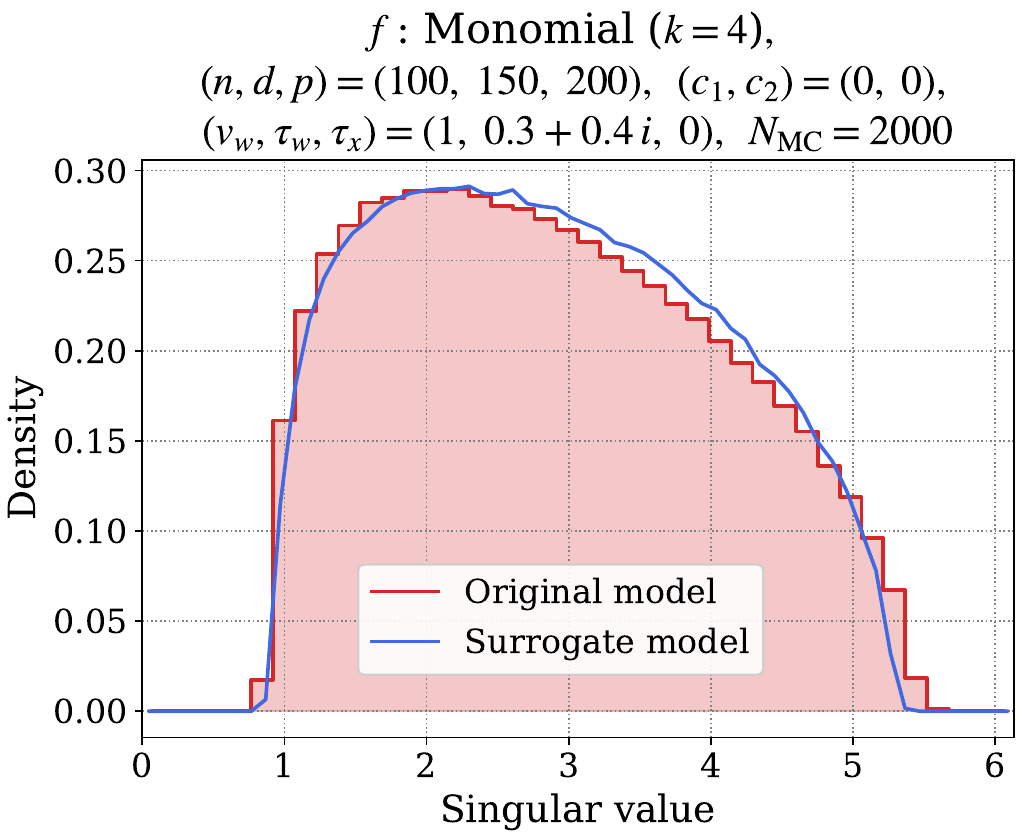}
	\caption{\label{fg:monomials_clipped}Comparison for the quartic activation after radial clipping of the preactivation matrix $H$. The upper $1\%$ and $10\%$ of the preactivation moduli are clipped in the left and right panels, respectively.}
\end{figure}

The larger discrepancy observed for higher-degree monomials may be attributed to the tails of the monomial outputs. Indeed, for $h\sim\CN(0,v_h,0)$, we have $\PP(|h^k|>r) = \exp\left(- r^{2/k} / v_h \right)$, showing slower decay for higher $k$, which may amplify rare large entries of the preactivation matrix $H$. To examine this effect, we radially clip the entries of $H$ before applying the quartic activation, leaving their phases unchanged. The clipping thresholds are estimated empirically as the $99$th and $90$th percentiles of the moduli of $10^7$ independent samples from $\CN(0,v_h,0)$, corresponding to clipping the upper $1\%$ and $10\%$, respectively. The Gaussian-surrogate parameters are then recomputed separately for the two clipped models. Figure~\ref{fg:monomials_clipped} shows the results. Clipping only the upper $1\%$ already substantially reduces the discrepancy seen in the top-right panel of Figure~\ref{fg:monomials}. When the upper $10\%$ is clipped, the agreement with the Gaussian surrogate is satisfactory despite the relatively small dimensions used here.

\subsection{Further Activation Functions}

\begin{figure}[tb]
	\centering
	\mbox{
	\includegraphics[width=.33\textwidth]{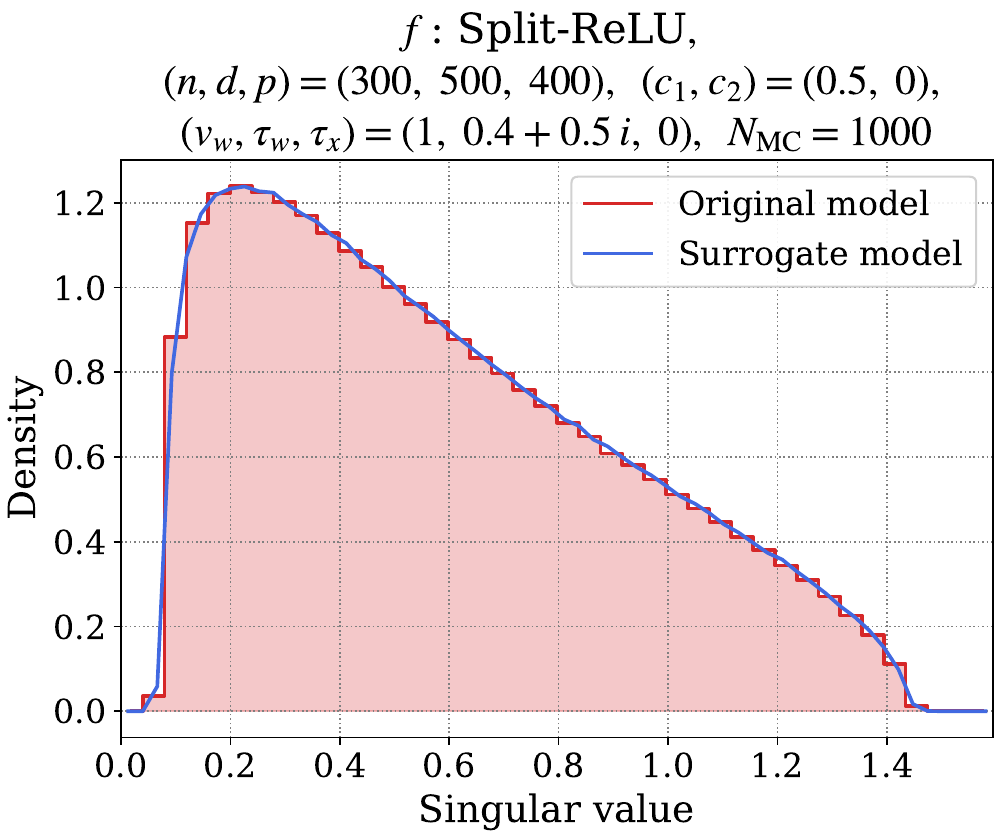}
	\includegraphics[width=.33\textwidth]{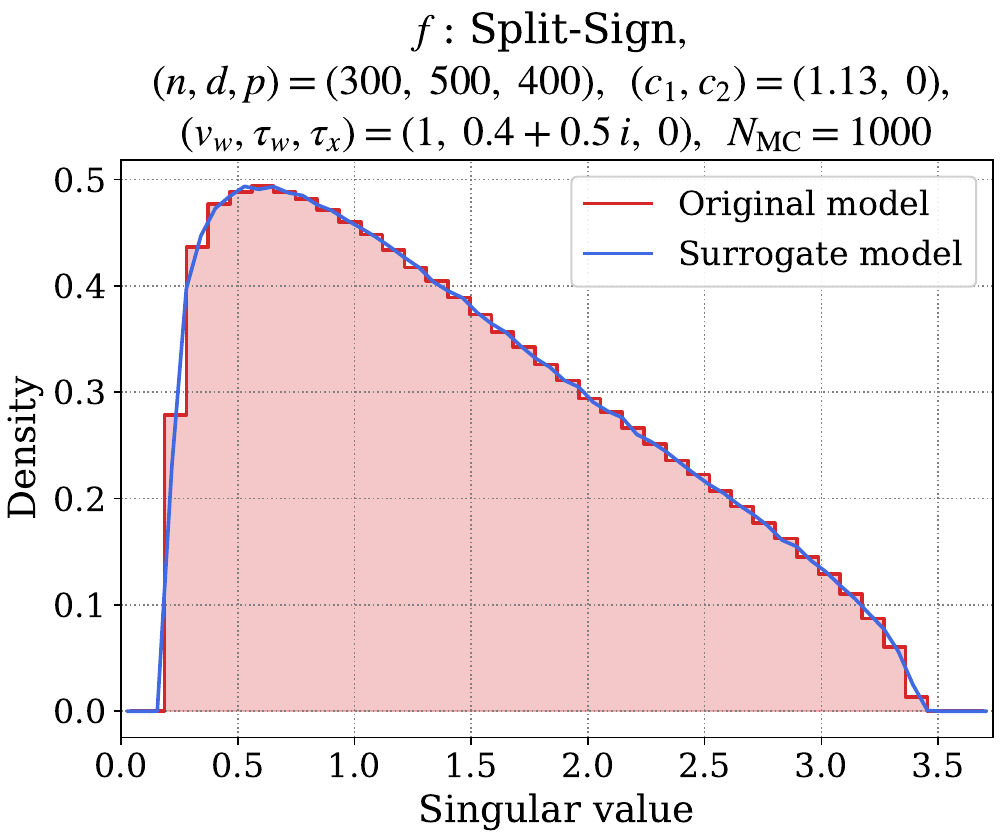}
	\includegraphics[width=.335\textwidth]{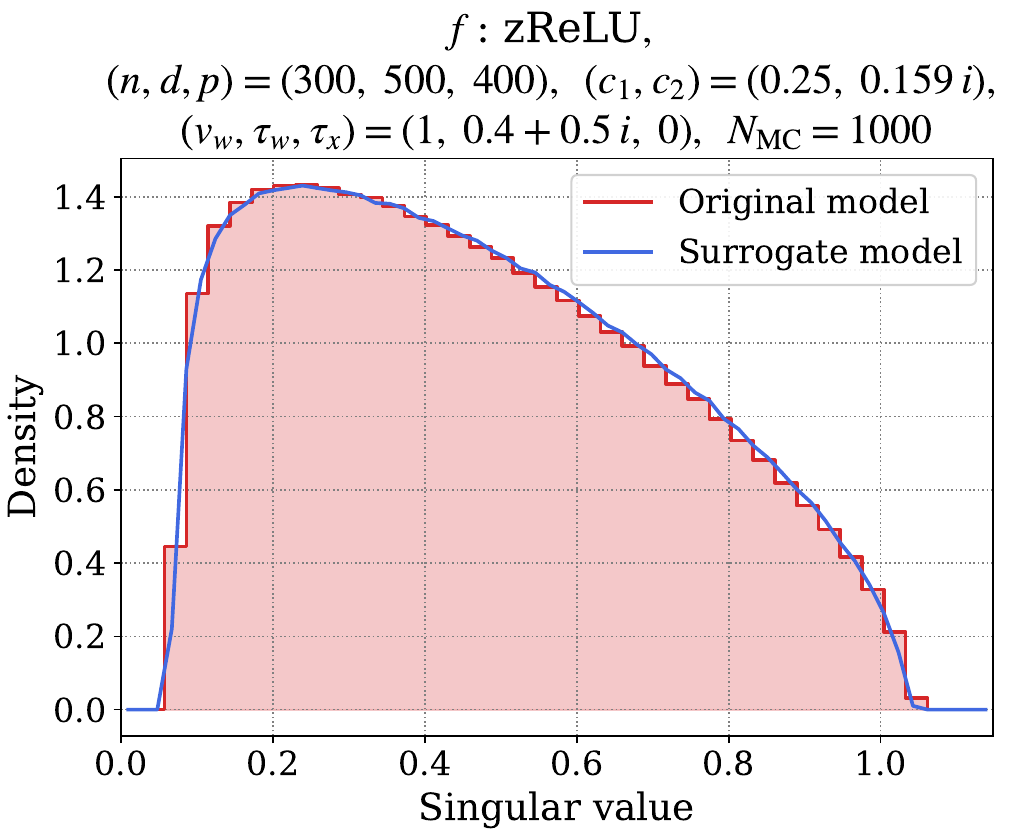}
	}
	\vspace{-4pt}
	\\
	\mbox{
	\includegraphics[width=.33\textwidth]{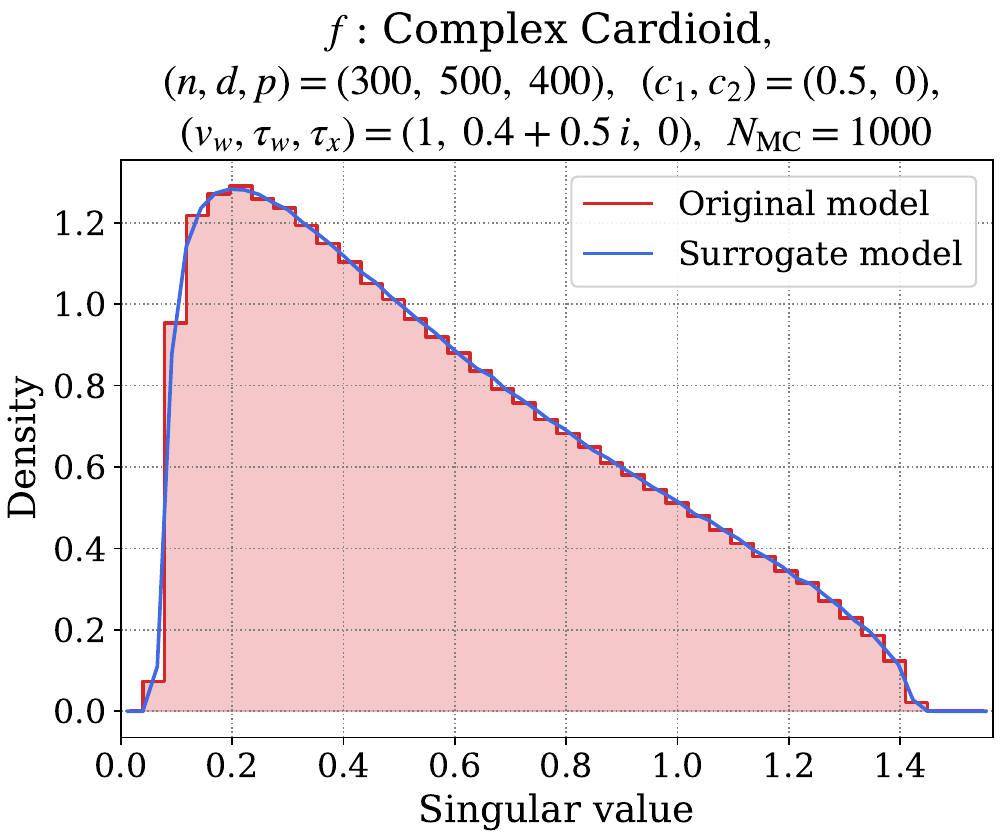}
	\includegraphics[width=.33\textwidth]{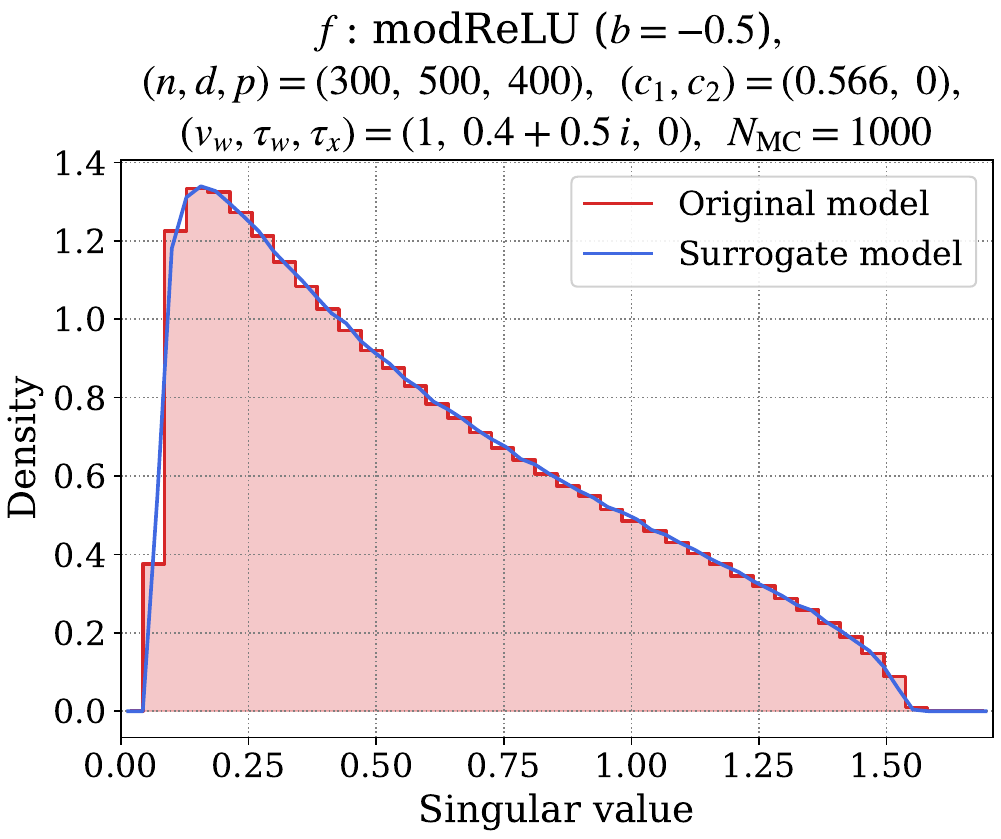}
	\includegraphics[width=.33\textwidth]{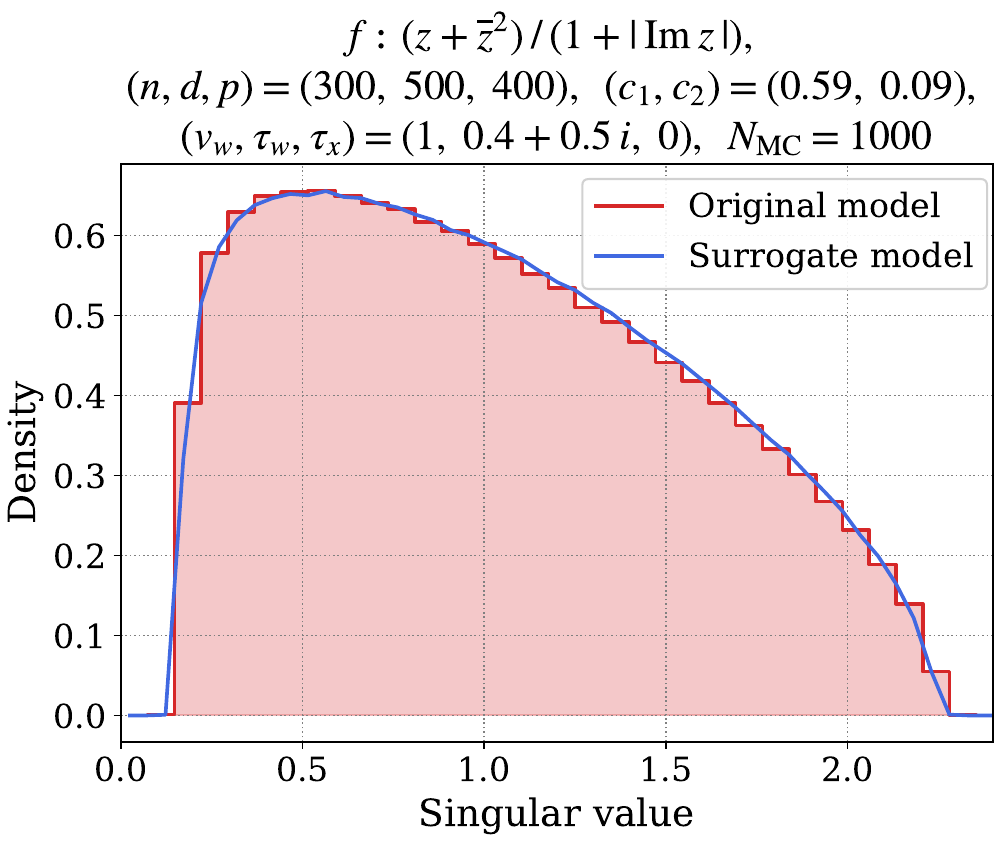}
	}
	\caption{\label{fg:various_f_proper_X}Comparison in the proper setting $\tau_x=0$ for Split-ReLU, Split-Sign, and zReLU in the top row, and Complex Cardioid, modReLU with $b=-0.5$, and $f(z)=(z+\bar{z}^2)/(1+|\Im z|)$ in the bottom row.}
\end{figure}
\begin{figure}[tb]
	\centering
	\hspace*{-10pt}
	\mbox{
	\includegraphics[width=.335\textwidth]{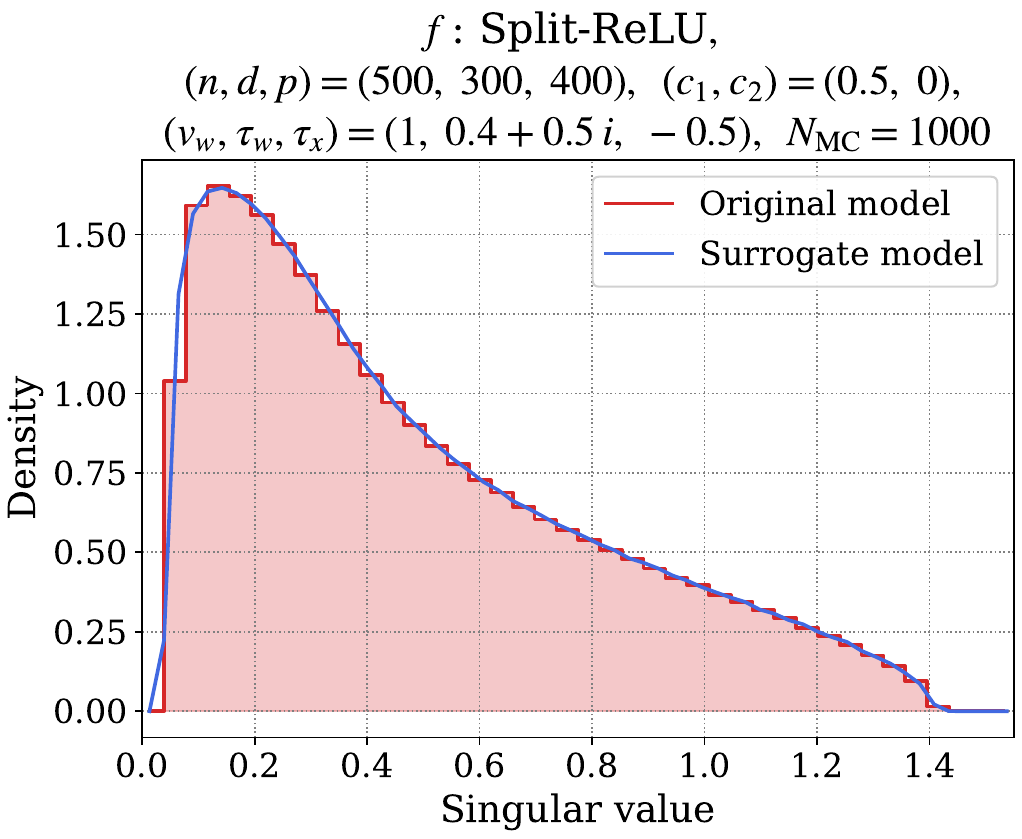}
	\includegraphics[width=.33\textwidth]{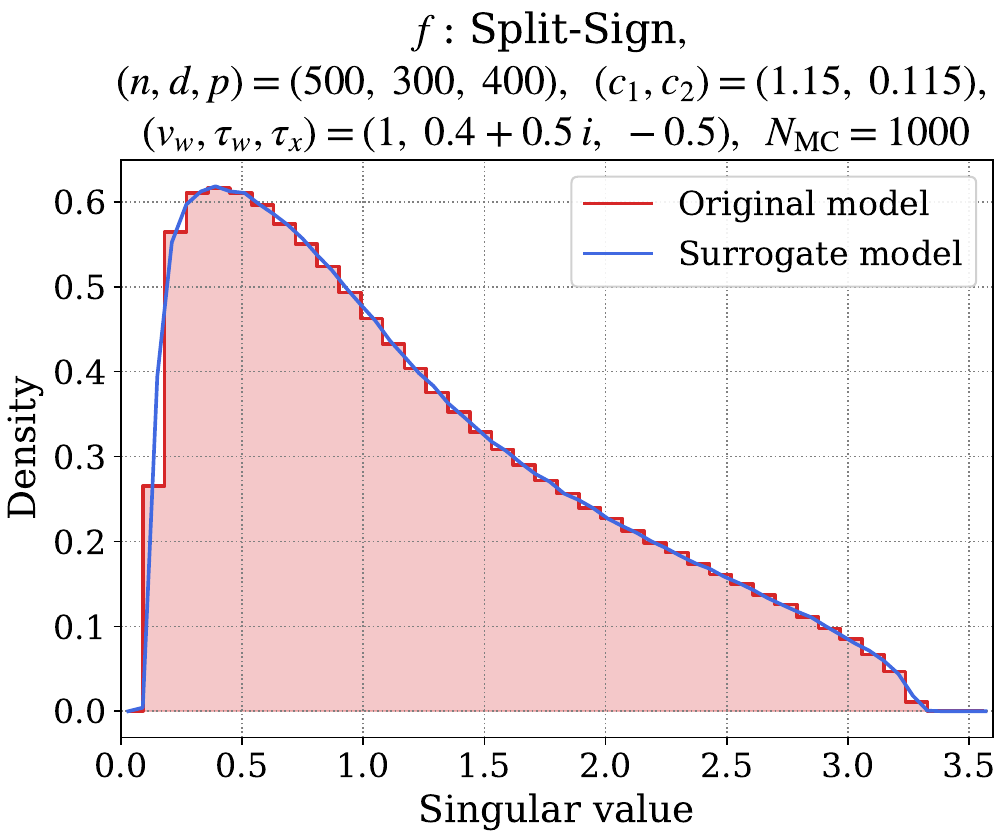}
	\includegraphics[width=.36\textwidth]{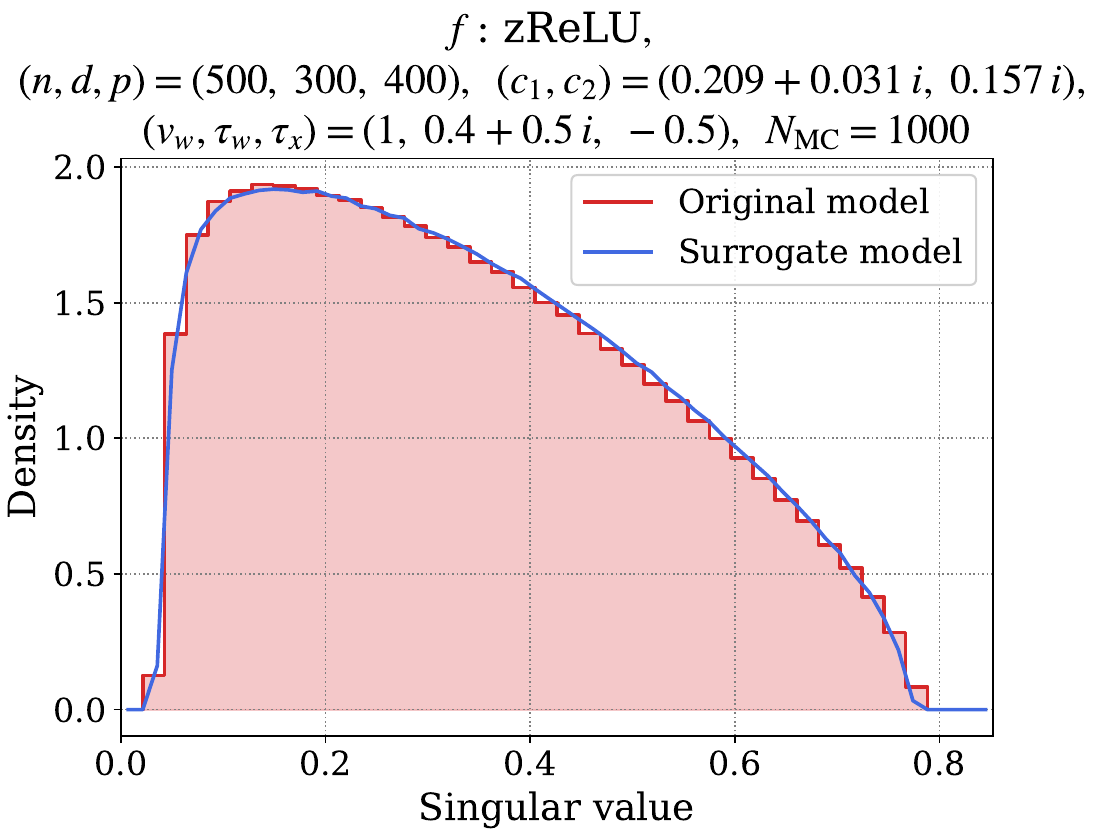}
	}
	\vspace{-4pt}
	\\
	\hspace*{-18pt}
	\mbox{
	\includegraphics[width=.33\textwidth]{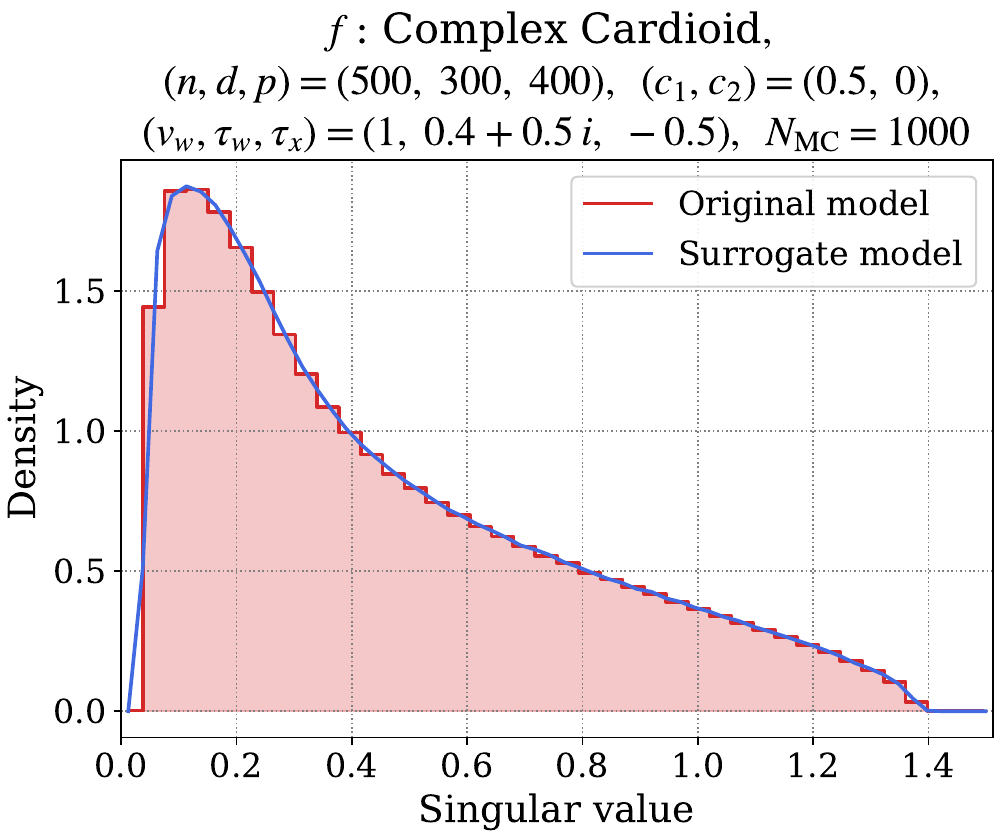}
	\includegraphics[width=.37\textwidth]{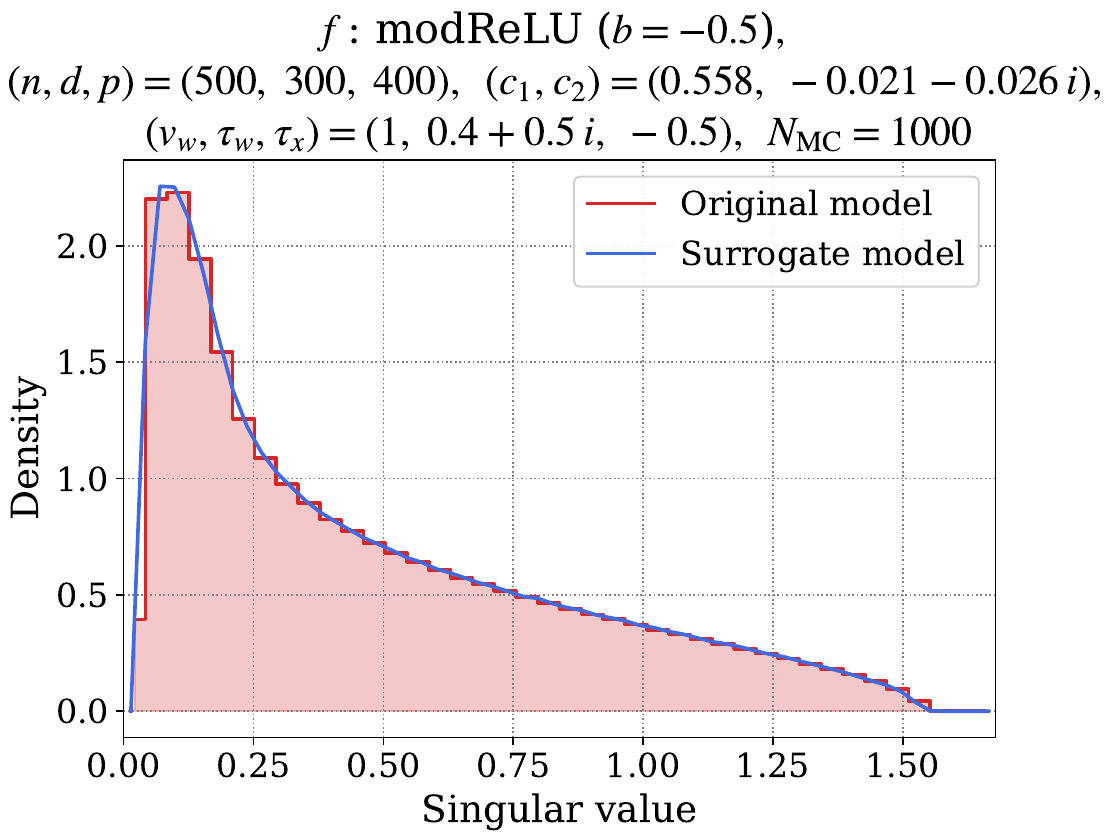}
	\includegraphics[width=.392\textwidth]{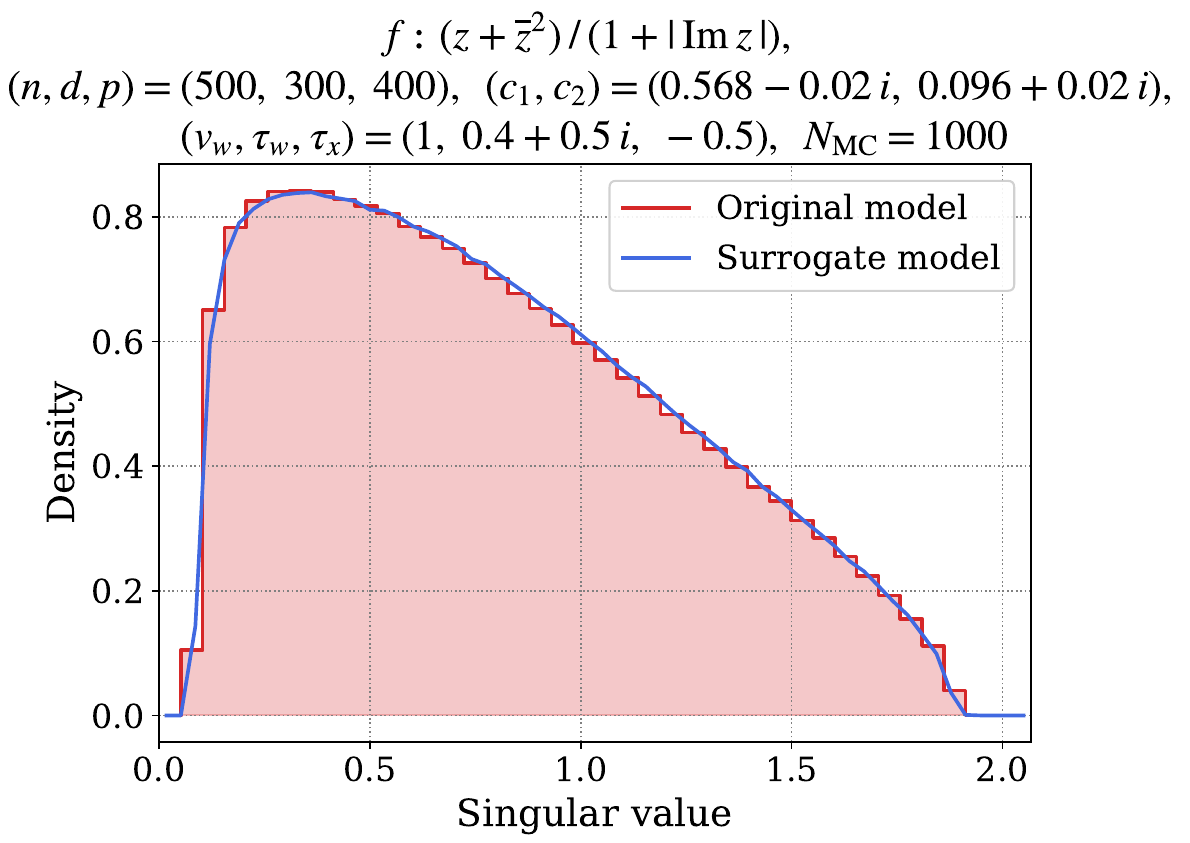}
	}
	\caption{\label{fg:various_f_improper_X}Comparison in the improper setting $\tau_x\tau_w\ne0$ for Split-ReLU, Split-Sign, and zReLU in the top row, and Complex Cardioid, modReLU with $b=-0.5$, and $f(z)=(z+\bar{z}^2)/(1+|\Im z|)$ in the bottom row.}
\end{figure}

We next examine a broader range of activation functions: Split-ReLU, Split-Sign, zReLU, Complex Cardioid, modReLU, and the nonlinear activation $f(z)=(z+\bar{z}^2)/(1+|\Im z|)$, introduced as an additional comparison case. Here, zReLU \autocite{Guberman2016} is defined by
\al{
	\mathrm{zReLU}(z) = 
	\begin{cases}
		z, & z\ne 0 \text{ and } 0 \leq\arg z\leq\pi/2\,,
		\\
		0, & \text{otherwise}.
	\end{cases}
}
These activation functions differ substantially in their regularity, symmetry, and treatment of the modulus and phase. In particular, Split-Sign and zReLU are discontinuous and therefore fall outside the regularity condition \eqref{eq:activation_increment_bound}.

We first consider the proper setting $\tau_x=0$. In this case, Lemma~\ref{lem:ooqwQOWc} gives $c_2=0$ for Split-ReLU and Split-Sign, Lemma~\ref{lem:c2zero} gives $c_2=0$ for Complex Cardioid, and Lemma~\ref{lem:c2_0_for_f_with_phase} gives $c_2=0$ for modReLU. As shown in Figure~\ref{fg:various_f_proper_X}, the agreement with the Gaussian surrogate is excellent for all six activation functions. 

Figure~\ref{fg:various_f_improper_X} shows the corresponding results in the improper setting $\tau_x\tau_w\ne0$. We interchange $n$ and $d$ relative to Figure~\ref{fg:various_f_proper_X} to examine a different dimensional configuration. Close agreement is again observed for all six activation functions.

\section{Conclusions}
\label{sc:conclusions}

In this work, we formulated a general complex GEP in terms of a widely linear Gaussian surrogate that retains both linear and conjugate-linear response channels, together with possibly improper Gaussian residual noise. We then proved a spectral version of this principle for the expected empirical spectral distribution of the feature covariance matrix. The proof combines a complex Hermite expansion with a graphical moment analysis that identifies the leading contraction structures, while the remaining approximation steps are controlled through Stieltjes-transform comparisons. We also characterized broad classes of complex activation functions for which one or both GEP coefficients vanish. The numerical experiments provide strong support for the spectral complex GEP across a diverse range of activation functions and covariance structures, including several settings beyond the assumptions of the theorem. Combined with the leave-one-out analysis of the widely linear Gaussian surrogate developed in the companion paper \autocite{kanazawa_rmt_2026}, the present result provides a route from the original nonlinear feature model to an explicit characterization of its limiting spectrum.

The present proof requires either the data matrix $X$ or the weight matrix $W$ to be proper and the activation function to satisfy the polynomially weighted Lipschitz condition \eqref{eq:activation_increment_bound}. The excellent numerical agreement observed when both matrices are improper and for activation functions outside this regularity condition suggests that these restrictions may be technical rather than intrinsic. Moreover, the rigorous result established here concerns an ensemble-averaged spectral distribution. An important next step is to determine whether the same equivalence holds at the level of individual matrix realizations, as suggested by the self-averaging behavior commonly observed in random matrix models. A further direction is to extend the general complex GEP beyond spectral observables and establish its consequences for training and test errors in high-dimensional supervised learning problems.

\section*{Acknowledgments}

The author used GPT-5.6 to assist with English-language editing and to improve the clarity and presentation of the manuscript. The tool was not used as a source of mathematical results. All AI-assisted revisions were reviewed by the author, who takes full responsibility for the final content.

\begin{appendices}
\section{Properties of Complex Hermite Polynomials}
\label{app:hermite}

In this appendix, we summarize the definitions and fundamental properties of complex Hermite polynomials \autocite{Ito1952,Ghanmi2008,Ghanmi2013,Ismail2015,Ismail2016}. Multivariate product formulas and their combinatorial interpretations were developed by \textcite{Ismail2015}. Motivated by this framework, we derive the corresponding formula for a jointly proper complex Gaussian vector in the covariance-based form needed for our graphical expansion.

Let $z \in \CC$ be a proper complex Gaussian random variable distributed as $z \sim \CN(0, v, 0)$ with variance $v > 0$.
The complex Hermite polynomials $\He_{p,q}(z, \bar{z})$ for nonnegative integers $p, q \ge 0$ are defined via the generating function $G(s, t, z)$ as
\al{
	G(s, t, z) \equiv 
	\exp\left( s z + t \bar{z} - st v \right) 
	= \sum_{p,q=0}^\infty 
	\frac{s^p t^q}{p! q!} \He_{p,q}(z, \bar{z})\,.
	\label{eq:gen_func_def}
}
These polynomials form a complete orthogonal basis, satisfying the orthogonality relation
\al{
	\EE\kkakko{\He_{p,q}(z, \bar{z}) \bar{\He_{m,n}(z, \bar{z})}} 
	= \delta_{pm} \delta_{qn} p! q! v^{p+q}\,.
	\label{eq:heorthog}
}
The explicit forms of the polynomials up to total degree three are given by
\al{
	\begin{array}{c}
		\He_{0,0} = 1\,, \quad  
		\He_{1,0} = z\,, \quad 
		\He_{0,1} = \bar{z}\,,
		\vspace{0.5em}
		\\
		\He_{2,0} = z^2 \,, \quad 
		\He_{1,1} = |z|^2-v \,, \quad 
		\He_{0,2} = \bar{z}^2 \,,
		\vspace{0.5em}
		\\
		\He_{3,0} = z^3 \,,~~~ 
		\He_{2,1} = z^2\bar{z} - 2vz \,,~~~
		\He_{1,2} = z\bar{z}^2 - 2v\bar{z} \,,~~~
		\He_{0,3} = \bar{z}^3 \,,
	\end{array}
	\hspace{-1em}
}
where the arguments $(z,\bar{z})$ were omitted for brevity. 

We begin with a generating-function identity for correlated complex Gaussian variables. It will be used to derive the higher-order cross-moments later.
\begin{lemma}
	\label{lem:prod_gen_func}
	Let $z_1,z_2,\dots,z_M$ be jointly proper zero-mean complex Gaussian
	random variables satisfying
	\al{
		\EE\bigl[|z_a|^2\bigr]
		& = v > 0
		\quad (a\in[M]),
		\\
		\EE[z_a\bar{z_b}]
		& = C_{ab} = \bar{C_{ba}}
		\quad (a,b\in[M],\ a\ne b),
		\\
		\EE[z_a z_b]
		& = 0 \quad (a,b\in[M]).
		\label{eq:4wf2563}
	}
	Then the expectation of the product of the $M$ generating functions is
	given by
	\al{
		\EE\left[
		\prod_{a=1}^{M}G(s_a,t_a,z_a)
		\right]
		=
		\exp\left(
		\sum_{a\ne b}s_a t_b C_{ab}
		\right).
	}
\end{lemma}
\begin{proof}
From the definition of $G$, we have
\al{
    \EE \left[ \prod_{a=1}^M G(s_a, t_a, z_a) \right] 
    &= \EE \left[ \prod_{a=1}^M \exp(s_a z_a + t_a \bar{z_a} - s_a t_a v) \right] \\
    &= \exp\left( -v \sum_{a=1}^M s_a t_a \right) \EE \left[ \exp\left( \sum_{a=1}^M (s_a z_a + t_a \bar{z_a}) \right) \right] \,. \label{eq:proof_step1}
}
Since $\{z_a\}_{a=1}^M$ are jointly Gaussian random variables with zero mean, their linear combination $J \equiv \sum_{a=1}^M (s_a z_a + t_a \bar{z_a})$ is also a zero-mean Gaussian random variable. Thus, we have the relation $\EE[\exp(J)] = \exp\big(\frac{1}{2}\EE[J^2]\big)$. We evaluate the second moment of $J$ explicitly by expanding the square:
\al{
    \EE[J^2] &= \EE \left[ \left( \sum_{a=1}^M 
    (s_a z_a + t_a \bar{z_a}) \right)^2 \right] 
    \\
    & = \sum_{a=1}^M \sum_{b=1}^M \EE \left[ 
    	s_a s_b z_a z_b + t_a t_b 
    	\bar{z_a} \bar{z_b} 
    	+ 2 s_a t_b z_a \bar{z_b} 
    \right]
    \\
    & = 2 \sum_{a=1}^M \sum_{b=1}^M 
    s_a t_b \EE[z_a \bar{z_b}] 
    \\
    & = 2 v \sum_{a=1}^M s_a t_a 
    + 2 \sum_{a \ne b} s_a t_b C_{ab} \,, 
    \label{eq:proof_step2}
}
where \eqref{eq:4wf2563} was used. Finally, substituting \eqref{eq:proof_step2} back into \eqref{eq:proof_step1}, we observe that the self-variance terms in the two exponentials cancel exactly:
\al{
    \EE \left[ \prod_{a=1}^M G(s_a, t_a, z_a) \right] 
    & = \exp\left( -v \sum_{a=1}^M s_a t_a \right) 
    \exp\left( v \sum_{a=1}^M s_a t_a + 
    \sum_{a \ne b} s_a t_b C_{ab} \right) 
    \\
    & = \exp\left( \sum_{a \ne b} s_a t_b C_{ab} \right) \,.
}
This completes the proof.
\end{proof}

\begin{lemma}\label{lem:HHrho}
Let $z_1, z_2$ be complex Gaussian variables satisfying the conditions in Lemma~\ref{lem:prod_gen_func}. Then
\al{
	\EE \kkakko{\He_{p,q}(z_1, \bar{z_1}) \He_{m,n}(z_2, \bar{z_2})} 
	& = \delta_{pn} \delta_{qm} p! q! C_{12}^p C_{21}^q \,.
	\label{eq:399t4sd}
	\\
	\EE \kkakko{\He_{p,q}(z_1, \bar{z_1}) \bar{\He_{m,n}(z_2, \bar{z_2})}} 
	& = \delta_{pm} \delta_{qn} p! q! C_{12}^p C_{21}^q \,,
	\label{eq:t34Uas}
}
\end{lemma}
\begin{proof}
Setting $M=2$ in Lemma~\ref{lem:prod_gen_func}, the expectation of the generating functions simplifies to
\al{
    \EE [ G(s_1, t_1, z_1) G(s_2, t_2, z_2) ] 
    = \exp \big( s_1 t_2 C_{12} + s_2 t_1 C_{21} \big) \,.
}
Expanding the right-hand side in a Taylor series and matching the coefficients of the monomials $s_1^p t_1^q s_2^m t_2^n$ against the definition \eqref{eq:gen_func_def}, we immediately obtain \eqref{eq:399t4sd}. Equation \eqref{eq:t34Uas} follows from \eqref{eq:399t4sd} by invoking the conjugate symmetry of the complex Hermite polynomials, $\bar{\He_{m,n}(z_2, \bar{z_2})}=\He_{n,m}(z_2, \bar{z_2})$.
\end{proof}

Notice that taking the limit $C_{12},C_{21}\to v$ in \eqref{eq:t34Uas} yields the previously stated orthogonality condition \eqref{eq:heorthog}, as it should. 

The generating-function argument in Lemma~\ref{lem:HHrho} extends to a product of an arbitrary number of complex Hermite polynomials and yields the following cross-moment formula.

\begin{proposition}[Cross-moments of complex Hermite polynomials] \label{prop:he_M}
Let $z_1,z_2,\dots,z_M$ be complex Gaussian variables satisfying the conditions in Lemma~\ref{lem:prod_gen_func}. There exists a set of nonnegative integers $\{m_{ab}\}_{a,b\in[M]}$ satisfying
\al{
	p_a = \sum_{b=1}^M m_{ab}\,, 
	\qquad 
	q_a=\sum_{b=1}^M m_{ba}\,,
	\qquad m_{aa}=0 \quad\text{for all }a\in[M]	
	\label{eq:degree_matching} 
}
if and only if
\al{
	\sum_{a=1}^M p_a=\sum_{a=1}^M q_a\equiv S\,,
	\qquad p_a+q_a\leq S \quad\text{for all }a\in[M]\,.
}
If these conditions fail, then
\al{
	\EE\left[ \prod_{a=1}^M \He_{p_a,q_a}(z_a,\bar{z_a}) \right] =0\,.
}
If they hold, then
\al{
	\EE\left[ \prod_{a=1}^M \He_{p_a,q_a}(z_a,\bar{z_a}) \right] 
	= \left(\prod_{a=1}^M p_a!q_a!\right) \sum_{\{m_{ab}\}} 
	\prod_{a\ne b} \frac{C_{ab}^{m_{ab}}}{m_{ab}!}\,, 
	\label{eq:three_point_exp}
}
where the sum runs over all configurations satisfying \eqref{eq:degree_matching}.
\end{proposition}

\begin{proof}
Applying Lemma~\ref{lem:prod_gen_func} yields
\al{
    \EE \left[ \prod_{a=1}^{M} G(s_a, t_a, z_a) \right]
    & = \exp\left( \sum_{a \neq b} s_a t_b C_{ab} \right)
    \\
    &= \prod_{a \neq b} \left(
    \sum_{m_{ab}=0}^\infty \frac{(s_a t_b C_{ab})^{m_{ab}}}{m_{ab}!}
    \right)
    \\
    &= \sum_{\{m_{ab} \ge 0\}}
    \prod_{a \neq b} \frac{(s_a t_b C_{ab})^{m_{ab}}}{m_{ab}!} \,.
    \label{eq:G_M_product_expanded}
}
To find the expectation of the specific polynomials, we extract the coefficient of the monomial $\prod_{a=1}^{M} s_a^{p_a} t_a^{q_a}$. Geometrically, this coefficient is a sum over directed multigraphs with $M$ nodes, where node $a$ has out-degree $p_a$ and in-degree $q_a$, and self-loops are excluded ($m_{aa}=0$).

Summing the degree-matching conditions $p_a = \sum_{b \neq a} m_{ab}$ and $q_a = \sum_{b \neq a} m_{ba}$ over $a \in [M]$ gives the global conservation constraint $\sum_{a=1}^M p_a = \sum_{a=1}^M q_a \equiv S$.

The existence of a configuration $\{m_{ab}\}$ can be characterized by Hall's marriage theorem \autocite{Hall1935,Bondy1976book,Diestel2025book}. For each node $a$, place $p_a$ outgoing copies of $a$ on the left and $q_a$ incoming copies of $a$ on the right. We first draw all candidate edges: an outgoing copy of node $a$ is joined to an incoming copy of node $b$ whenever $a\ne b$. These candidate edges do not yet form a matching. A valid configuration $\{m_{ab}\}$ corresponds to a perfect matching selected from them, where $m_{ab}$ counts the selected edges from copies of node $a$ to copies of node $b$.

Hall's marriage theorem states that such a matching exists if and only if every subset $X$ on the left has at least $|X|$ possible partners on the right. We consider two cases:
\bi
	\item
	If $X$ consists entirely of outgoing copies of one node $a$, its possible partners are all incoming copies except those of node $a$. Their number is $S-q_a$. Since $|X|$ can be as large as $p_a$, the required condition is $p_a+q_a\leq S$.
	\item
	If $X$ contains outgoing copies of at least two distinct nodes, every incoming copy is a possible partner of at least one element of $X$. All $S$ incoming copies are therefore available, and the condition is automatically satisfied.
\ei
Hence, a valid configuration exists if and only if $p_a+q_a\leq S$ for every $a\in[M]$. If this condition fails, the coefficient of the monomial is zero. Otherwise, the admissible configurations are precisely the collections $\{m_{ab}\}$ satisfying the degree-matching conditions above. Returning to the coefficient extraction in \eqref{eq:G_M_product_expanded} and comparing it with \eqref{eq:gen_func_def} yields \eqref{eq:three_point_exp}.
\end{proof}

\begin{figure}[tb]
	\centering
	\includegraphics[width=2.2cm]{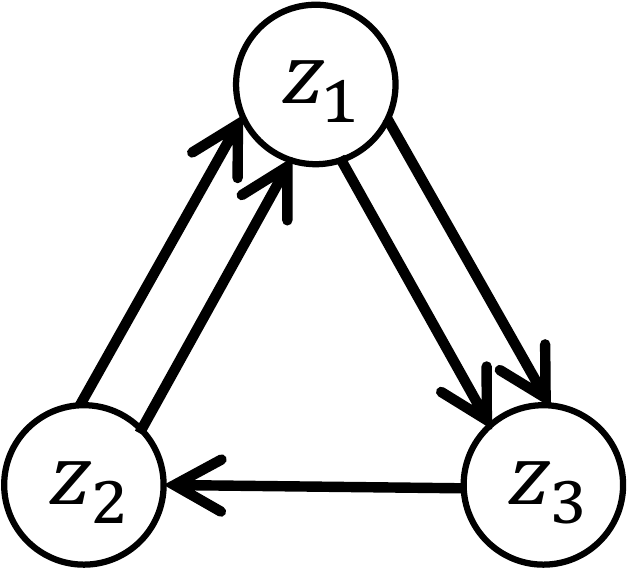}
	\hspace{15mm}
	\includegraphics[width=2.2cm]{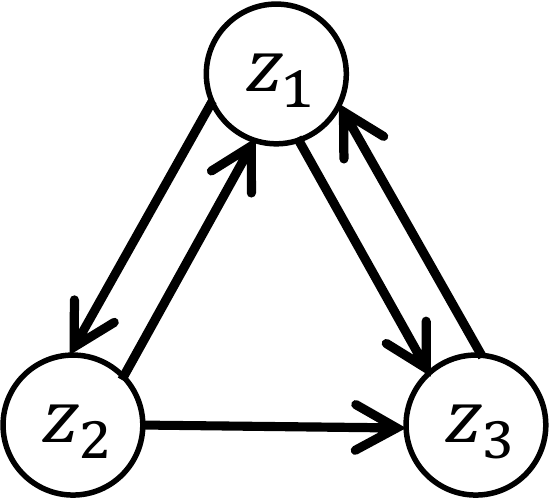}
	\put(-95, 21){\scalebox{3.0}{$\mathbf{+}$}}
	\caption{\label{fg:tripleheex}Graphs contributing to the triple product of complex Hermite polynomials \eqref{eq:34gdAPGOEIR}.}
\end{figure}

To illustrate Proposition~\ref{prop:he_M}, we note that the following expectations all vanish:
\al{
	& \EE[\He_{2,2}(z_1) \He_{1,0}(z_2) \He_{0,1}(z_3)]
	= \EE[\He_{3,2}(z_1) \He_{1,1}(z_2) \He_{0,1}(z_3)]
	\notag
	\\
	=~ & \EE[\He_{3,3}(z_1) \He_{1,2}(z_2) \He_{1,0}(z_3)]
	= \EE[\He_{4,2}(z_1) \He_{0,2}(z_2) \He_{1,1}(z_3)]
	= 0\,.
}
In each case, although the global constraint $\sum_a p_a = \sum_a q_a \equiv S$ is satisfied, the local constraint $p_a + q_a \le S$ is violated. Conversely, an example of a nonvanishing expectation is given by 
\al{
	\EE[\He_{2,2}(z_1)\He_{2,1}(z_2)\He_{1,2}(z_3)]
	= 4 C_{21}^2 C_{13}^2 C_{32} 
	+ 16 C_{12} C_{21} C_{13} C_{31} C_{23} \,,
	\label{eq:34gdAPGOEIR}
}
where these two terms correspond to the sum of contributions from two distinct valid graphs, as illustrated in Figure~\ref{fg:tripleheex}.

\section{Frobenius Bound for Stieltjes Transforms}
\label{app:frobenius}

In this appendix, we derive a general bound on the difference between the normalized Stieltjes transforms of $A^*A/n$ and $B^*B/n$ in terms of the Frobenius norms of $A$, $B$, and $A-B$. For $A\in\CC^{n\times p}$ and $z\in\CC_+$, define
\al{
	K_A \equiv \frac{1}{n}A^*A\,, \qquad
	G_A(z) \equiv \frac{1}{K_A-z\1_p}\,,
}
and define $K_B$ and $G_B(z)$ analogously.

\begin{lemma}[Comparison bound for Stieltjes transforms]
\label{lem:stieltjes_comparison}
For arbitrary matrices $A,B\in\CC^{n\times p}$ and every $z\in\CC_+$,
\al{
	\left| m_A(z)-m_B(z) \right|
	& \leq \frac{
	\left(\|A\|_{\rm F}+\|B\|_{\rm F}\right)
	\|A-B\|_{\rm F}
	}{ (\Im z)^2 np }\,.
	\label{eq:stieltjes_comparison_bound}
}
Furthermore, if $A$ and $B$ are random matrices satisfying $\EE\big[\|A\|_{\rm F}^2\big] < \infty$ and $\EE\big[\|B\|_{\rm F}^2 \big] < \infty$, then
\al{
	\EE[ \left|m_A(z)-m_B(z)\right| ]
	& \leq \frac{\sqrt{2}}{(\Im z)^2}
	\sqrt{
		\frac{\EE\big[ \|A\|_{\rm F}^2 + \|B\|_{\rm F}^2 \big]}{np}
	}
	\sqrt{ \frac{ \EE\big[ \|A-B\|_{\rm F}^2 \big]}{np}}\,.
	\label{eq:expected_stieltjes_comparison_bound}
}
\end{lemma}

\begin{proof}
The resolvent identity gives
\al{
	G_A(z)-G_B(z) = G_A(z)(K_B-K_A)G_B(z)\,.
	\label{eq:resolvent_identity_AB}
}
Since $K_A$ and $K_B$ are Hermitian and $z\in\CC_+$, the spectral theorem yields
\al{
	\|G_A(z)\|_{\rm op} \leq \frac{1}{\Im z}\,,
	\qquad
	\|G_B(z)\|_{\rm op} \leq \frac{1}{\Im z}\,.
	\label{eq:resolvent_operator_bounds}
}
Using \eqref{eq:resolvent_identity_AB}, \eqref{eq:resolvent_operator_bounds} and the inequalities
\al{
	|\Tr(CD)| \leq \|C\|_{\rm op}\|D\|_{\star}\,,
	\qquad \|CD\|_{\rm op} \leq
	\|C\|_{\rm op}\|D\|_{\rm op}
}
for arbitrary matrices $C$ and $D$ of compatible dimensions, we obtain
\al{
	\left| m_A(z)-m_B(z) \right|
	& = \frac{1}{p}
	\left| \Tr\left[
		G_A(z)(K_B-K_A)G_B(z)
	\right] \right|
	\\
	& = \frac{1}{p} \left|
	\Tr\left[ G_B(z)G_A(z)(K_B-K_A) \right]
	\right|
	\\
	& \leq \frac{1}{p}
	\|G_B(z)G_A(z)\|_{\rm op} \|K_B-K_A\|_{\star}
	\\
	& \leq \frac{1}{p}\frac{1}{(\Im z)^2}
	\|K_A-K_B\|_{\star}\,.
	\label{eq:stieltjes_nuclear_intermediate}
}

The difference above may be represented as
\al{
	K_A-K_B & = \frac{1}{n}( A^*A-B^*B )
	\\
	& = \frac{1}{n} \big\{ A^*(A-B) + (A-B)^*B \big\}\,.
	\label{eq:covariance_difference_factorization}
}
The nuclear norm satisfies $\|C^*D\|_{\star} \leq \|C\|_{\rm F}\|D\|_{\rm F}$. Applying this to the two terms in \eqref{eq:covariance_difference_factorization} gives
\al{
	\|K_A-K_B\|_{\star}
	& \leq \frac{1}{n} \big( \|A^*(A-B)\|_{\star}
	+ \|(A-B)^*B\|_{\star} \big)
	\\
	& \leq \frac{1}{n}
	( \|A\|_{\rm F} + \|B\|_{\rm F} )\|A-B\|_{\rm F}\,.
	\label{eq:covariance_nuclear_bound}
}
Combining \eqref{eq:stieltjes_nuclear_intermediate} and \eqref{eq:covariance_nuclear_bound} proves \eqref{eq:stieltjes_comparison_bound}.

If $A$ and $B$ are random, taking the expectation in \eqref{eq:stieltjes_comparison_bound} and applying the Cauchy--Schwarz inequality gives
\al{
	\EE[ | m_A(z)-m_B(z) | ] 
	& \leq \frac{1}{(\Im z)^2 np}
	\sqrt{
		\EE\big[ ( \|A\|_{\rm F}+\|B\|_{\rm F} )^2 \big]
	}\sqrt{
		\EE\big[ \|A-B\|_{\rm F}^2 \big]
	}
	\\
	& \leq \frac{\sqrt{2}}{(\Im z)^2}
	\sqrt{
		\frac{\EE\big[\|A\|_{\rm F}^2+\|B\|_{\rm F}^2 \big]}{np}
	}
	\sqrt{
		\frac{\EE\big[\|A-B\|_{\rm F}^2 \big]}{np}
	}\,,
}
where the last inequality follows from $( \|A\|_{\rm F}+\|B\|_{\rm F} )^2 \leq 2 \big(\|A\|_{\rm F}^2+\|B\|_{\rm F}^2 \big)$. This proves \eqref{eq:expected_stieltjes_comparison_bound}.
\end{proof}

\section{Moments of Normalized Gaussian Columns}
\label{app:normalized_column_wick}

This appendix derives an exact pairwise-contraction formula for the fixed-degree moments of the column-normalized Gaussian matrix $\WW$ introduced in Section~\ref{sc:notation}. Although column normalization makes the entries within each column dependent and non-Gaussian, their moments involve the same pairings and data-index identifications as the Wick expansion for the corresponding column of $W$. The normalization changes only the scalar coefficient associated with each pairing, and this coefficient differs from its Gaussian Wick counterpart by $O(d^{-1})$.

Recall from \eqref{eq:column-normalization} that
\al{
	\WW_{\mu\alpha} = \sqrt{dv_w}
	\frac{W_{\mu\alpha}}{
		\sqrt{\sum_{\nu=1}^d |W_{\nu\alpha}|^2}
	} \,.
}
For $\varepsilon\in\{+,-\}$, write
\al{
	W_{\mu\alpha}^{(+)} = W_{\mu\alpha}\,,
	\qquad W_{\mu\alpha}^{(-)}
	= \bar{W_{\mu\alpha}}\,,
}
and define $\WW_{\mu\alpha}^{(+)}$ and $\WW_{\mu\alpha}^{(-)}$ analogously. We also set
\al{
	\kappa_{+-} = \kappa_{-+} = v_w\,,
	\qquad \kappa_{++} = \tau_w\,, 
	\qquad \kappa_{--} = \bar{\tau_w} \,.
}
Thus, $\kappa_{\varepsilon\varepsilon'}$ is the Gaussian second moment associated with a pair of types $\varepsilon$ and $\varepsilon'$.

\begin{lemma}[Pairing formula for normalized Gaussian columns]
\label{lem:normalized_column_wick}

Fix an integer $r\geq1$, a feature index $\alpha\in[p]$, data indices $\mu_1,\ldots,\mu_{2r}\in[d]$ and signs $\varepsilon_1,\ldots,\varepsilon_{2r}\in\{+,-\}$. Then
\al{
	\EE\bigg[
		\prod_{j=1}^{2r}
		\WW_{\mu_j\alpha}^{(\varepsilon_j)}
	\bigg]
	= \sum_{\pi\in\mathcal P_2(2r)}
	c_{\pi,d} \prod_{\{j,\ell\}\in\pi}
	\delta_{\mu_j\mu_\ell} \,.
	\label{eq:normalized_column_wick}
}
Here, $\mathcal P_2(2r)$ is the set of pairings of $[2r]$. If a pairing $\pi$ contains $a$ pairs of type $(+,-)$ or $(-,+)$, $b$ pairs of type $(+,+)$ and $c$ pairs of type $(-,-)$, so that $a+b+c=r$, then its coefficient is given exactly by
\al{
	c_{\pi,d} = \frac{(dv_w)^r}{\Gamma(r)}
	\int_0^\infty
	\rmd t ~ t^{r-1}
	\bigl(1+2v_wt+D_wt^2\bigr)^{-d/2}
	v(t)^a\tau(t)^b\bar{\tau(t)}^{\,c} \, ,
	\label{eq:normalized_column_pairing_integral}
}
where
\al{
	D_w = v_w^2-|\tau_w|^2 \,,
	\qquad v(t) =
	\frac{v_w+D_wt}{1+2v_wt+D_wt^2} \,,
	\qquad \tau(t) =
	\frac{\tau_w}{1+2v_wt+D_wt^2} \,.
	\label{eq:tilted_column_parameters}
}
For every fixed $r$ and every $\pi\in\mathcal P_2(2r)$,
\al{
	c_{\pi,d} = \prod_{\{j,\ell\}\in\pi}
	\kappa_{\varepsilon_j\varepsilon_\ell}
	+ O_r(d^{-1}) \,.
	\label{eq:normalized_column_wick_coefficient}
}
The error is uniform over the data indices and the feature index. Its implicit constant may depend on $r$, $v_w$ and $\tau_w$, but not on $d$, $\alpha$ or $\mu_1,\ldots,\mu_{2r}$.

Moments of odd total degree vanish. Moreover, for any fixed product involving finitely many feature columns, independence of the columns gives the product of the corresponding columnwise expansions. For every fixed collection of columnwise pairings, the resulting coefficient differs from its Gaussian Wick counterpart by $O(d^{-1})$.
\end{lemma}

\begin{proof}
Fix one feature index $\alpha$ and abbreviate $w_\mu=W_{\mu\alpha}$. Let
\al{
	P(w) = \prod_{j=1}^{2r}
	w_{\mu_j}^{(\varepsilon_j)} .
}
Since $P$ is homogeneous of total degree $2r$,
\al{
	\EE \left[P(\WW_{\scalebox{0.6}{$\bullet$}\alpha})\right]
	= (dv_w)^r \EE\left[ \frac{P(w)}{\|w\|_2^{2r}} \right]. 	\label{eq:normalized_column_homogeneous}
}
The Laplace representation
\al{
	\frac{1}{x^r} = \frac{1}{\Gamma(r)}
	\int_0^\infty \rmd t~t^{r-1} \rme^{-tx}\,, 
	\qquad x>0
}
may be inserted into \eqref{eq:normalized_column_homogeneous}. Indeed, since 
\al{
	\frac{|P(w)|}{\|w\|_2^{2r}} =
	\prod_{j=1}^{2r} \frac{|w_{\mu_j}|}{\|w\|_2} 
	\leq 1\,,
}
Tonelli's theorem gives
\al{
	\int_0^\infty \rmd t ~ t^{r-1}
	\EE\left[ |P(w)| \rme^{-t\|w\|_2^2} \right]
	= \Gamma(r) \, \EE\left[ 
	\frac{|P(w)|}{\|w\|_2^{2r}}
	\right] \leq \Gamma(r)<\infty \,.
}
Thus, Fubini's theorem allows us to take the expectation inside the integral, yielding
\al{
	\EE\left[P(\WW_{\scalebox{0.6}{$\bullet$}\alpha})
	\right] = \frac{(dv_w)^r}{\Gamma(r)} 
	\int_0^\infty \rmd t ~ t^{r-1}
	\EE\left[ P(w) \rme^{-t\|w\|_2^2} 
	\right].
	\label{eq:normalized_column_laplace}
}

We next evaluate the expectation in \eqref{eq:normalized_column_laplace}. First suppose that $|\tau_w|<v_w$, and let $D_w \equiv v_w^2-|\tau_w|^2>0$. For one entry $z\sim\mathcal N_{\mathbb C}(0,v_w,\tau_w)$, multiplication of its density by $\rme^{-t|z|^2}$ only shifts the coefficient of $|z|^2$ in the Gaussian exponent. After normalization, the tilted law is again centered complex Gaussian. Its normalizing factor is
\al{
	\phi(t) \equiv \EE\big[\rme^{-t|z|^2}\big] =
	\bigl(1+2v_wt+D_wt^2\bigr)^{-1/2} \,.
}
The covariance and pseudo-covariance under the normalized tilted law are precisely $v(t)$ and $\tau(t)$ in \eqref{eq:tilted_column_parameters}. These formulas follow by inverting the covariance matrix $\bmat v_w & \tau_w \\ \bar{\tau_w} & v_w \emat$, adding $t\1_2$, and then inverting the resulting matrix.

The entries $w_1,\ldots,w_d$ remain independent under the tilt. Hence
\al{
	\EE\left[ \rme^{-t\|w\|_2^2} \right]
	= \phi(t)^d	=
	\bigl(1+2v_wt+D_wt^2\bigr)^{-d/2} \,.
}
Applying the improper complex Wick formula under the normalized tilted law yields
\al{
	\EE\left[
		P(w) \rme^{-t\|w\|_2^2}
	\right]
	= 
	\bigl(1+2v_wt+D_wt^2\bigr)^{-d/2}
	\sum_{\pi\in\mathcal P_2(2r)}
	\prod_{\{j,\ell\}\in\pi}
	\kappa_{\varepsilon_j\varepsilon_\ell}(t)
	\delta_{\mu_j\mu_\ell} \, ,
	\label{eq:kappawick}
}
where
\al{
	\kappa_{+-}(t) = \kappa_{-+}(t)
	= v(t),
	\qquad
	\kappa_{++}(t) = \tau(t),
	\qquad
	\kappa_{--}(t) = \bar{\tau(t)} \,.
}
Substituting \eqref{eq:kappawick} into \eqref{eq:normalized_column_laplace} proves \eqref{eq:normalized_column_wick} and the exact integral formula \eqref{eq:normalized_column_pairing_integral}. In particular, the data-index identifications are exactly the same as in the ordinary Gaussian Wick expansion.

It remains to expand the pairing coefficient as $d\to\infty$ and compare its leading term with the corresponding Gaussian Wick coefficient. Let
\al{
	\rho = \frac{|\tau_w|}{v_w},
	\qquad
	\delta = 1-\rho^2
	= \frac{D_w}{v_w^2}
	\in[0,1] \,.
}
After the change of variables $u=dv_wt$, \eqref{eq:normalized_column_pairing_integral} becomes
\al{
	c_{\pi,d} = {} &
	v_w^a\tau_w^b\bar{\tau_w}^{\,c}
	\frac{1}{\Gamma(r)}
	\int_0^\infty \rmd u~
	u^{r-1} \left(1 + \frac{\delta u}{d} \right)^a
	\left(
		1 + \frac{2u}{d} + \frac{\delta u^2}{d^2}
	\right)^{-d/2-r}.
}
We next show that, for fixed $r$, $a$ and $\delta\in[0,1]$, the integral above satisfies
\al{
	\frac{1}{\Gamma(r)} \int_0^\infty
	\rmd u ~ u^{r-1}
	\left(1 + \frac{\delta u}{d}\right)^a
	\left(
		1 + \frac{2u}{d} + \frac{\delta u^2}{d^2}
	\right)^{-d/2-r}
	= 1 + O_r(d^{-1}) \,.
	\label{eq:normalized_column_laplace_estimate}
}
To verify this estimate, we divide the integral into the three ranges $[0,d^{1/4}]$, $(d^{1/4},d/2]$ and $(d/2,\infty)$ and estimate their contributions separately.
\bnum
	\item
	For $0\leq u\leq d^{1/4}$, taking logarithms and using the Taylor expansion gives, uniformly over this range and $\delta\in[0,1]$,
	\al{
		\left(1+\frac{\delta u}{d}\right)^a
		\left(
			1+\frac{2u}{d}+\frac{\delta u^2}{d^2}
		\right)^{-d/2-r}
		= \rme^{-u}
		\left[
			1 + O_r \left(\frac{u+u^2}{d}\right)
		\right].
	}
	Since
	\al{
		\int_0^\infty \rmd u\,
		u^{r-1}(u+u^2)\rme^{-u} < \infty\,,
	}
	the integrated error over this range is $O_r(d^{-1})$.
	\item
	For $d^{1/4}<u\leq d/2$, the conditions $0\leq\delta\leq1$ and $a\leq r$ imply
	\al{
		\left(1 + \frac{\delta u}{d}\right)^a
		\leq \left(\frac{3}{2}\right)^r,
		\qquad
		\log\left(
		1 + \frac{2u}{d} + \frac{\delta u^2}{d^2}
		\right)
		\geq \log\left(1+\frac{2u}{d}\right)
		\geq \frac{u}{d}\,.
	}
	The integrand is therefore bounded by $C_r u^{r-1}\rme^{-u/2}$, whose integral over this range is $O_r\bigl(d^{(r-1)/4}\rme^{-d^{1/4}/2}\bigr)$ and hence $O_r\bigl(\rme^{-c d^{1/4}}\bigr)$ for some $c>0$.
	\item
	For $u>d/2$, set $x=u/d$. Since
	\al{
		1 + \frac{\delta u}{d} \leq 1+x \,,
		\qquad
		1 + \frac{2u}{d}+\frac{\delta u^2}{d^2} \geq 1+2x \,,
	}
	we obtain, for all sufficiently large $d$,
	\al{
		& \int_{d/2}^{\infty}\rmd u\, u^{r-1}
		\left(1+\frac{\delta u}{d}\right)^a
		\left(
			1 + \frac{2u}{d}+\frac{\delta u^2}{d^2}
		\right)^{-d/2-r}
		\notag
		\\
		& \qquad \leq
		d^r \int_{1/2}^{\infty}\rmd x\,
		x^{r-1}(1+2x)^{-d/2}
		= O_r\bigl(d^r2^{-d/2}\bigr)\,.
	}
\enum
Thus, the entire complementary integral is $o(d^{-1})$, uniformly in $\delta\in[0,1]$. This proves \eqref{eq:normalized_column_laplace_estimate} and hence \eqref{eq:normalized_column_wick_coefficient} when $|\tau_w|<v_w$.

If $|\tau_w|=v_w$, the complex Gaussian variable is supported on a real one-dimensional subspace after multiplication by a deterministic phase. The same calculation reduces to the corresponding real Gaussian integral. Equivalently, the formulas above extend continuously to $D_w=0$, and the estimate \eqref{eq:normalized_column_wick_coefficient} remains valid.

Odd-degree moments vanish because the distribution of each normalized column is symmetric under a sign change. Finally, distinct feature columns are independent. Applying the one-column formula to each distinct feature index gives the product expansion for finitely many columns. Since the number of columns and the total degree are fixed, a finite product of coefficients of the form $c_{\pi,d}=c_{\pi,\infty}+O(d^{-1})$ again differs from its Gaussian counterpart by $O(d^{-1})$.
\end{proof}

\end{appendices}

\begin{refcontext}[sorting=nyt]
\printbibliography
\end{refcontext}
\end{document}